\PassOptionsToPackage{table}{xcolor}

\documentclass[oneside,12pt,letterpaper]{article}

\usepackage{authblk} %

\usepackage{amsmath,amssymb,bm,bbold}
\usepackage{amsfonts,amsthm,mathrsfs,mathtools}
\usepackage{stmaryrd}
\SetSymbolFont{stmry}{bold}{U}{stmry}{m}{n}
\usepackage{upgreek}
\usepackage{stmaryrd}
\usepackage{econometrics}
\usepackage{flexisym}
\usepackage{breqn}

\usepackage{tikz}
\usetikzlibrary{positioning}
\usepackage{xr-hyper}
\usepackage{pdfpages}
\usepackage{pgf,interval}

\usepackage{float}
\usepackage{booktabs, array, longtable, ragged2e, pdflscape}
\usepackage{subcaption}

\usepackage{hyperref}
\hypersetup{
  colorlinks = false,
  urlbordercolor = {blue},
  linkbordercolor = {red},
  citebordercolor = {green}
}
\usepackage[page,toc,titletoc,title]{appendix}
\usepackage{makeidx}
\usepackage[nottoc,numbib]{tocbibind}
\hypersetup{bookmarksopen=true}

\usepackage{cleveref}
\renewcommand{\cref}{\Cref}
\renewcommand{\crefrange}{\Crefrange}
\AtBeginEnvironment{appendices}{\crefalias{section}{appendix}}

\usepackage{algorithm}
\usepackage{algpseudocode}
\algnewcommand\algorithmicinput{\textbf{Initialize:}}
\algnewcommand\Initialize{\item[\algorithmicinput]}

\makeatletter
\let\save@mathaccent\mathaccent
\newcommand*\if@single[3]{%
  \setbox0\hbox{${\mathaccent"0362{#1}}^H$}%
  \setbox2\hbox{${\mathaccent"0362{\kern0pt#1}}^H$}%
  \ifdim\ht0=\ht2 #3\else #2\fi
  }
\newcommand*\rel@kern[1]{\kern#1\dimexpr\macc@kerna}
\newcommand*\widebar[1]{\@ifnextchar^{{\wide@bar{#1}{0}}}{\wide@bar{#1}{1}}}
\newcommand*\wide@bar[2]{\if@single{#1}{\wide@bar@{#1}{#2}{1}}{\wide@bar@{#1}{#2}{2}}}
\newcommand*\wide@bar@[3]{%
  \begingroup
  \def\mathaccent##1##2{%
    \let\mathaccent\save@mathaccent
    \if#32 \let\macc@nucleus\first@char \fi
    \setbox\z@\hbox{$\macc@style{\macc@nucleus}_{}$}%
    \setbox\tw@\hbox{$\macc@style{\macc@nucleus}{}_{}$}%
    \dimen@\wd\tw@
    \advance\dimen@-\wd\z@
    \divide\dimen@ 3
    \@tempdima\wd\tw@
    \advance\@tempdima-\scriptspace
    \divide\@tempdima 10
    \advance\dimen@-\@tempdima
    \ifdim\dimen@>\z@ \dimen@0pt\fi
    \rel@kern{0.6}\kern-\dimen@
    \if#31
      \overline{\rel@kern{-0.6}\kern\dimen@\macc@nucleus\rel@kern{0.4}\kern\dimen@}%
      \advance\dimen@0.4\dimexpr\macc@kerna%
      \let\final@kern#2%
      \ifdim\dimen@<\z@ \let\final@kern1\fi
      \if\final@kern1 \kern-\dimen@\fi
    \else
      \overline{\rel@kern{-0.6}\kern\dimen@#1}%
    \fi
  }%
  \macc@depth\@ne%
  \let\math@bgroup\@empty \let\math@egroup\macc@set@skewchar%
  \mathsurround\z@ \frozen@everymath{\mathgroup\macc@group\relax}%
  \macc@set@skewchar\relax
  \let\mathaccentV\macc@nested@a%
  \if#31
    \macc@nested@a\relax111{#1}%
  \else
    \def\gobble@till@marker##1\endmarker{}%
    \futurelet\first@char\gobble@till@marker#1\endmarker%
    \ifcat\noexpand\first@char A\else%
      \def\first@char{}%
    \fi
    \macc@nested@a\relax111{\first@char}%
  \fi
  \endgroup
}
\makeatother

\usepackage[american]{babel}
\usepackage{csquotes}
\usepackage[
  style=authoryear,maxnames=4,natbib=true,uniquename=false,
  isbn=false,url=false,eprint=false,date=year,dashed=false
]{biblatex}
\bibliography{wig-r.bib}
\DefineBibliographyStrings{english}{bibliography = {References}}

\usepackage{xcolor}
\usepackage{cancel} %
\usepackage[left=2.5cm, right=2.5cm, top=3cm]{geometry}
\usepackage{indentfirst}
\usepackage{fancyhdr}
\newcommand{\ubar}[1]{\text{\b{$#1$}}} %
\usepackage{enumitem}

\usepackage{pifont}%
\newcommand{\interior}[1]{%
  {\kern0pt#1}^{\mathrm{o}}%
}

\DeclareMathOperator*{\argmin}{arg\,min}

\newtheorem{theorem}{Theorem}%

\newtheorem{proposition}[theorem]{Proposition}
\newtheorem{lemma}[theorem]{Lemma}

\newtheorem{definition}[theorem]{Definition}

\newtheorem*{remark}{Remark}

\newtheorem{problem}[theorem]{Problem}
\newtheorem{update}{Update}

\counterwithin{update}{section}
\counterwithin{algorithm}{section}

\title{Deriving the Gradients of Some Popular Optimal Transport Algorithms\footnote{
  This version has seen many major revisions of the first version \citep{xie2025},
  most notably in \cref{sec:sinkhorn-gradient,sec:wdl} regarding the differentiation of the 
  OT algorithms.
}}
\author{Fangzhou Xie}
\affil{Department of Economics, Rutgers University}
\date{\today}

\begin{document}
\maketitle

\begin{abstract}
  In this note, I review entropy-regularized Monge-Kantorovich problem in Optimal Transport,
  and derive the gradients of several popular algorithms popular in Computational Optimal Transport,
  including the Sinkhorn algorithms, Wasserstein Barycenter algorithms, and the Wasserstein Dictionary Learning
  algorithms.
\end{abstract}
\newpage

\tableofcontents
\newpage

\listofalgorithms
\addtocontents{loa}{\def\string\figurename{Algorithm}}

\vspace{8em}

\begin{figure}[H]
  \centering
  \begin{tikzpicture}[
      round/.style={circle, draw=black, very thick, minimum size=8mm},
    ]\centering
    \node[round] (review)                            {2};
    \node[round] (sinkhorn)    [below=of review]     {3};
    \node[round] (sinkgrad)    [right=of sinkhorn]   {4};
    \node[round] (barycenter)  [below=of sinkhorn]   {5};
    \node[round] (barygrad)    [right=of barycenter] {6};
    \node[round] (wdl)         [below=of barygrad]   {7};
    \draw[->,very thick] (review.south) -- (sinkhorn.north);
    \draw[->,very thick] (sinkhorn.east) -- (sinkgrad.west);
    \draw[->,very thick] (sinkhorn.south) -- (barycenter.north);
    \draw[->,very thick] (barycenter.east) -- (barygrad.west);
    \draw[->,very thick] (barygrad.south) -- (wdl.north);
    \draw[->,very thick,dashed] (sinkgrad.south) -- (barygrad.north);
  \end{tikzpicture}
  \caption{Diagram of the structure of the paper.}\label{fig:article-diagram}
\end{figure}

\newpage

\section{Introduction}

The French mathematician Gaspard Monge (1746 -- 1818) was allegedly interested in the question
of how to efficiently dig dirt from the earth and use it to patch the castle wall,
i.e.~moving dirt from one location to another\footnote{
  ``Lorsqu'on doit tranfporter des terres d'un lieu dans un autre,
  on a contume de donner le nom de \textit{D{\'e}blai} au volume des terres que l'on doit
  tranfporter, \& le nom de \textit{Remblai} {\`a} l'dfpace qu'elles doivent occuper apr{\'e}s le tranfport''.
  -- \fullcite{monge1781}.
}.
Therefore, the question of how to ``transporting'' dirt efficiently became the start of the study
of Optimal Transport.
The field has received much attention over the years,
because it is deeply connected to many interesting topics across different fields,
including but not limited to Linear Programming, PDE, statistics, image processing, natural language processing, etc.
Interested readers should refer to the now standard monographs by Villani \parencite*{villani2003,villani2008}
for the mathematical theories, \citet{santambrogio2015a} for applied mathematicians,
\citet{galichon2016} for economists, and \citet{peyre2019} for computational aspects of OT.

This note can be seen as a quick review of some basic algorithms on Computational Optimal Transport (COT),
and provide the manually-derived gradients for the said algorithms.
This will also form the foundation for the companion R package\footnote{
  \url{https://github.com/fangzhou-xie/wig}
} \textbf{\textit{wig}},
which provides an efficient implementation of the computation (OT algorithms and their gradients) laid out
in this note and also the WIG models \citep{xie2020,xie2020a}.
The algorithms are also checked numerically by the Julia Zygote\footnote{
  \url{https://github.com/FluxML/Zygote.jl} \citep{innes2019}
} and ForwardDiff\footnote{
  \url{https://github.com/JuliaDiff/ForwardDiff.jl} \citep{revels2016}
} library to ensure the correct derivation of the gradients.

\citet{xie2020} already provides a package, \textbf{\textit{wigpy}}, for the computation of WIG models,
but the automatic differentiation was built upon the PyTorch \citep{paszke2017} and hence the package was written in Python.
The new \textbf{\textit{wig}} package, on the other hand, is in R with all the manually derived gradients implemented in C++.
This ensures the high performance of the computation, and also gets rid of the dependency on any AD library\footnote{
  The idea of bypassing the AD libraries to reduce dependencies for packages is also applied in the R package
  \textbf{\textit{rethnicity}} \citet{xie2022,xie2021}.
  More broadly speaking, writing packages with less dependencies is usually
  a good practice (\url{https://www.tinyverse.org/}).
}.

\cref{fig:article-diagram} provides a diagram of the structure of the paper.
For readers interested in entropic regularized OT problems,
they can refer to \cref{sec:review,sec:sinkhorn-algorithm} with optional \cref{sec:sinkhorn-gradient};
for readers interested in Wasserstein Barycenter problems,
they can refer to \cref{sec:review,sec:sinkhorn-algorithm,sec:wasserstein-barycenter};
with optional \cref{sec:barycenter-gradient};
for readers interested in Wasserstein Dictionary Learning,
they can refer to \cref{sec:review,sec:sinkhorn-algorithm,sec:wasserstein-barycenter,sec:wdl}.
The remainder of this paper will be organized as follows.
\cref{sec:review} briefly reviews the Monge--Kantorivich problem, its entropic regularized problem and solution.
\cref{sec:sinkhorn-algorithm} introduces the Sinkhorn algorithm and its variants to solve the entropic regularized OT problem.
\cref{sec:sinkhorn-gradient} derives the gradients and Jacobians of the family of Sinkhorn algorithms.
\cref{sec:wasserstein-barycenter} discusses the Wasserstein Barycenter problem.
\cref{sec:barycenter-gradient} shows the gradient of Sinkhorn-like algorithms to solve the Wasserstein Barycenter problem.
\cref{sec:wdl} considers the Wasserstein Dictionary Learning problem.
\cref{appendix:math-notation} and \cref{appendix:lemmas}
lists all the mathematical notations and lemmas to simplify the expressions used throughout this note.

\section{Monge-Kantorovich Problem and Its Entropic Regularization}\label{sec:review}
\sectionmark{Monge-Kantorovich}

The Monge--Kantorovich problem\footnote{
  This is also called Kantorovich relaxation of the original Monge problem,
  as the problem was first proposed by \citet{monge1781} and then relaxed by \citet{kantorovich1942}.
  For mathematical foundation of Optimal Transport theory, the classic references are \citet{villani2003,villani2008};
  for its application in Economics, see \citet{galichon2016};
  for an introduction for applied mathematicians, see \citet{santambrogio2015a};
  for the computational aspects (algorithms and their properties) and its applications in data science and machine learning,
  see \citet{peyre2019}.
} states that:
given two probability vectors\footnote{
  Sometimes also called ``histograms'' in the Computational Optimal Transport community, for example in \citet{peyre2019}.
} $\mathbf{a} \in \Sigma_M$ and $\mathbf{b} \in \Sigma_N$,
how to find a \textit{coupling matrix} $\mathbf{P} \in \mathbb{R}_+^{M \times N}$
where each $\mathbf{P}_{ij}$ describes the flow of masses from bin $i$ to bin $j$,
such that the total cost of moving masses are optimal (minimal).
The cost of moving one unit of mass from bin $i$ to bin $j$ is $\mathbf{C}_{ij}$ and $\mathbf{C} \in \mathbb{R}_+^{M\times N}$.

\begin{problem}[Monge--Kantorovich]
Let $\mathbf{a} \in \Sigma_M$ and $\mathbf{b} \in \Sigma_N$.
The set of all coupling matrices is:
\begin{equation}\label{eqn-def:coupling-matrix}
  \begin{aligned}
    \mathbf{U(a,b)} \equiv \left\{
    \mathbf{P} \in \mathbb{R}_+^{M \times N}:
    \mathbf{P} \cdot \mathbb{1}_N = \mathbf{a}
    \text{ and }
    \mathbf{P}^\top \cdot \mathbb{1}_M = \mathbf{b}
    \right\}.
  \end{aligned}
\end{equation}

Given a cost matrix $\mathbf{C} \in \mathbb{R}^{M \times N}$,
the Kantorovich optimal transport problem tries to find the solution of the following
\begin{equation}\label{eqn:loss-kantorovich}
  \begin{aligned}
    \ell_{\mathbf{C}}(\mathbf{a}, \mathbf{b})
    \equiv \min_{\mathbf{P} \in \mathbf{U}(\mathbf{a},\mathbf{b})}
    \langle \mathbf{C}, \mathbf{P}\rangle
    = \sum_{i,j} \mathbf{C}_{ij} \mathbf{P}_{ij}.
  \end{aligned}
\end{equation}
\end{problem}

To solve this problem in practice, one needs to resort to linear programming\footnote{
  See \citet[Chapter 3]{peyre2019} for a historical overview and related algorithms.
}, which can be very challenging when the problem becomes large enough.
Instead of finding the exact solution to the above problem, one can actually add an entropy regularization term
to make this problem convex and hence greatly ease the computation.
This was first proposed in the seminal work of \citet{cuturi2013}
which leverages Sinkhorn--Knopp double scaling algorithm
\footnote{
  Therefore, in Optimal Transport literature,
  people usually refer to the numeric algorithm to solve the entropic regularized Kantorovich problem
  as the ``Sinkhorn algorithm.''
  Hence, its extensions are also named as ``X-khorn algorithms'',
  for example ``Greenkhorn algorithm'' \citep{altschuler2017} or ``Screenkhorn'' \citep{alaya2019}.
  Because of its simplicity, it has been re-discovered multiple times in history and renamed accordingly,
  for example,
  it can be linked to Sch\"odinger bridge problem \citep{schrodinger1931} or the RAS method \citep{bacharach1970}.
  See \citet{knight2008,leonard2013,modin2024} for historical accounts.
}\citep{sinkhorn1964,sinkhorn1967,knight2008}.

\begin{problem}[Entropic Regularized OT Problem, \citet{cuturi2013,peyre2019}]\label{thm:entropic-regularized-OT-problem}
Given the coupling matrix $\mathbf{P} \in \mathbb{R}_+^{M \times N}$, its discrete entropy is defined as
\begin{equation*}
  \begin{aligned}
    \mathbf{H}(\mathbf{P})
    \equiv - \langle \mathbf{P}, \log(\mathbf{P}) - 1\rangle
    = - \sum_{i,j} \mathbf{P}_{ij} (\log(\mathbf{P}_{ij}) - 1).
  \end{aligned}
\end{equation*}

With the cost matrix $\mathbf{C}$, the entropic regularized loss function is
\begin{equation}\label{eqn:entropic-regularized-OT-loss}
  \begin{aligned}
    \ell^\varepsilon_{\mathbf{C}}(\mathbf{a}, \mathbf{b})
    \equiv
    \langle \mathbf{P}, \mathbf{C}\rangle - \varepsilon \mathbf{H}(\mathbf{P}),
  \end{aligned}
\end{equation}

and the entropic regularized OT problem is thus
\begin{equation}\label{eqn:entropic-regularized-OT-problem}
  \begin{aligned}
    \min_{\mathbf{P} \in \mathbf{U}(\mathbf{a},\mathbf{b})}
    \ell^\varepsilon_{\mathbf{C}}(\mathbf{a}, \mathbf{b})
    \equiv
    \min_{\mathbf{P} \in \mathbf{U}(\mathbf{a},\mathbf{b})}
    \langle \mathbf{P}, \mathbf{C}\rangle - \varepsilon \mathbf{H}(\mathbf{P}).
  \end{aligned}
\end{equation}
\end{problem}

The Sinkhorn problem has a unique optimal solution, since it is $\varepsilon$-strongly convex \citep{peyre2019}.
\begin{proposition}[Uniqueness of the Entropic OT Solution]\label{thm:prop-uniqueness-ot-solution}
  The solution to \cref{eqn:entropic-regularized-OT-problem} is unique and has the form
  \begin{equation}\label{eqn:regularized-OT-solution}
    \begin{aligned}
      \mathbf{P}_{ij} = \mathbf{u}_i \mathbf{K}_{ij} \mathbf{v}_j,
    \end{aligned}
  \end{equation}

  or in matrix notation
  \begin{equation}\label{eqn:regularized-OT-solution-matrix-form}
    \begin{aligned}
      \mathbf{P} = \diag\mathbf{u} \cdot \mathbf{K} \cdot \diag\mathbf{v}
    \end{aligned}
  \end{equation}
  where $i \in \left\{1, \ldots, M\right\}$ and $j \in \left\{1, \ldots, N\right\}$ with two (unknown)
  scaling variables $\mathbf{u} \in \mathbb{R}^M$ and $\mathbf{v} \in \mathbb{R}^N$.
\end{proposition}

\begin{proof}
  This is essentially Proposition 4.3 in \citet{peyre2019}.
  Let $\mathbf{f} \in \mathbb{R}^M$ and $\mathbf{g} \in \mathbb{R}^N$ be two dual variables for the two constraints
  $\mathbf{P}\cdot \mathbb{1}_N = \mathbf{a}$ and $\mathbf{P}^\top\cdot \mathbb{1}_M = \mathbf{b}$, respectively.
  The Lagrangian of \cref{eqn:entropic-regularized-OT-problem} becomes
  \begin{equation}\label{eqn:regularized-OT-lagrangian}
    \begin{aligned}
      \mathscr{L}(\mathbf{P}, \mathbf{f}, \mathbf{g}) =
      \langle \mathbf{P}, \mathbf{C}\rangle - \varepsilon \mathbf{H}(\mathbf{P})
      -
      \langle \mathbf{f}, \mathbf{P} \cdot \mathbb{1}_m - \mathbf{a}\rangle
      -
      \langle \mathbf{g}, \mathbf{P}^\top \cdot \mathbb{1}_n - \mathbf{b}\rangle.
    \end{aligned}
  \end{equation}

  The FOC (first order condition) gives us
  \begin{equation*}
    \begin{aligned}
      \frac{
        \partial \mathscr{L}(\mathbf{P}, \mathbf{f}, \mathbf{g})
      }{
        \partial \mathbf{P}_{ij}
      } = \mathbf{C}_{ij} + \varepsilon \log(\mathbf{P}_{ij}) - \mathbf{f}_i - \mathbf{g}_j = 0
    \end{aligned}
  \end{equation*}

  Therefore we can solve for the optimal coupling matrix $\mathbf{P}$ as:
  \begin{equation*}
    \begin{aligned}
      \log(\mathbf{P}_{ij})
      =
      \frac{\mathbf{f}_i}{\varepsilon}
      - \frac{\mathbf{C}_{ij}}{\varepsilon}
      + \frac{\mathbf{g}_j}{\varepsilon}
    \end{aligned}
  \end{equation*}

  or
  \begin{equation}\label{eqn:optimal-coupling}
    \begin{aligned}
      \mathbf{P}_{ij}
      =
      \exp(\frac{\mathbf{f}_i}{\varepsilon})
      \exp(-\frac{\mathbf{C}_{ij}}{\varepsilon})
      \exp(\frac{\mathbf{g}_j}{\varepsilon})
      =
      \mathbf{u}_i \mathbf{K}_{ij} \mathbf{v}_j,
    \end{aligned}
  \end{equation}

  where $\mathbf{u}_i = \exp(\frac{\mathbf{f}_i}{\varepsilon})$,
  $\mathbf{v}_j = \exp(\frac{\mathbf{g}_j}{\varepsilon})$, and
  $\mathbf{K}_{ij} = \exp(-\frac{\mathbf{C}_{ij}}{\varepsilon})$, respectively.
\end{proof}

\begin{remark}[]
  Note that computation of $\mathbf{P}_{ij}$ are effectively in the exponential domain.
  If we choose to find
  optimal solution for $\mathbf{P}$ by updating $\mathbf{u}$ and $\mathbf{v}$,
  which would grow exponentially in $\mathbf{f}$ and $\mathbf{g}$.
  This is indeed the case for the vanilla Sinkhorn Algorithm, and thus it has numeric instability issue
  which need to be overcome by the ``log-stabilization'' technique.
  The idea is simple: just move all computation to the log-domain
  (see \cref{subsec:log-sinkhorn} for details).
\end{remark}

\section{Sinkhorn Algorithms}\label{sec:sinkhorn-algorithm}
\sectionmark{Sinkhorn}

In this section, I will discuss three different versions of Sinkhorn algorithms.

\subsection{Vanilla Sinkhorn Algorithm}\label{subsec:vanilla-sinkhorn}

Given the optimal solution to the Sinkhorn problem in \cref{eqn:regularized-OT-solution-matrix-form},
with the constraints in \cref{eqn-def:coupling-matrix}, we have
\begin{equation}
  \begin{aligned}
    \diag\mathbf{u} \cdot \mathbf{K} \cdot \diag\mathbf{v} \cdot \mathbb{1}_N  = \mathbf{a}
    \quad\text{and}\quad
    \diag\mathbf{v} \cdot \mathbf{K}^\top \cdot \diag\mathbf{u} \cdot \mathbb{1}_M  = \mathbf{b}.
  \end{aligned}
\end{equation}

Since we also have $\diag\mathbf{v} \cdot \mathbb{1}_N = \mathbf{v}$ and $\diag\mathbf{u}\cdot \mathbb{1}_M = \mathbf{u}$,
we have that
\begin{equation}\label{eqn:vanilla-sinkhorn-solution-multiply}
  \begin{aligned}
    \mathbf{u} \odot (\mathbf{K} \mathbf{v}) = \mathbf{a}
    \quad\text{and}\quad
    \mathbf{v} \odot (\mathbf{K}^\top \mathbf{u}) = \mathbf{b}.
  \end{aligned}
\end{equation}

Equivalently, we can also write the solution as
\begin{equation}\label{eqn:vanilla-sinkhorn-solution-division}
  \begin{aligned}
    \mathbf{u} = \mathbf{a} \oslash (\mathbf{K} \mathbf{v})
    \quad\text{and}\quad
    \mathbf{v} = \mathbf{b} \oslash (\mathbf{K}^\top \mathbf{u}).
  \end{aligned}
\end{equation}

In the above equations, $\odot$ and $\oslash$ refer to element-wise (Hadamard) multiplication and division, respectively.
Therefore, we can leverage \cref{eqn:vanilla-sinkhorn-solution-division} to update $\mathbf{u}$ and $\mathbf{v}$
till convergence and conclude with optimal coupling $\mathbf{P}$ by \cref{eqn:regularized-OT-solution-matrix-form},
as we also know that the solution is unique as shown in \cref{thm:prop-uniqueness-ot-solution}.
\begin{update}[Updating Equations for Vanilla Sinkhorn Algorithm]\label{update:vanilla-sinkhorn}
  At any iteration $\ell = 1, \ldots, L$ with initialized $\mathbf{v}^{(0)} = \mathbb{1}_N$, we have
  \begin{equation}
    \begin{aligned}
      \mathbf{u}^{(\ell)} \equiv \mathbf{a} \oslash (\mathbf{K}\mathbf{v}^{(\ell-1)})
      \quad\text{and}\quad
      \mathbf{v}^{(\ell)} \equiv \mathbf{b} \oslash (\mathbf{K}^\top \mathbf{u}^{(\ell)}).
    \end{aligned}
  \end{equation}
\end{update}

\begin{remark}
  In \cref{update:vanilla-sinkhorn}, the choice of $\mathbf{v}^{(0)}$ is rather arbitrary,
  as long as it is not $\mathbb{0}$.
  The most common choice is to consider the vectors of ones, i.e. $\mathbb{1}_N$.
  Changing the initialization vector will result in different solutions to $\mathbf{u}$ and $\mathbf{v}$,
  since they are only defined up to a multiplicative constant.
  But $\mathbf{u}$ and $\mathbf{v}$ will converge regardless of the initialization,
  and result in the same optimal coupling matrix $\mathbf{P}$ upon convergence.
  See \citet[p.82]{peyre2019} for discussion.
\end{remark}

After $L$ iterations, the scaling variables $\mathbf{u}^{(L)}$ and $\mathbf{v}^{(L)}$ can be used to compute the
optimal coupling matrix
$\mathbf{P}^{(L)} = \diag\mathbf{u}^{(L)} \cdot \mathbf{K} \cdot \diag \mathbf{v}^{(L)}$.
From \cref{update:vanilla-sinkhorn}, we can arrive at the vanilla Sinkhorn algorithm (\cref{algo:vanilla-sinkhorn}).

\begin{algorithm}[H]
  \caption{Vanilla Sinkhorn Algorithm}
  \begin{algorithmic}[1]\label{algo:vanilla-sinkhorn}
    \Require $\mathbf{a} \in \Sigma_M$, $\mathbf{b} \in \Sigma_N$, $\mathbf{C} \in \mathbb{R}^{M\times N}$, $\varepsilon > 0$.
    \Initialize $\mathbf{u} = \mathbb{1}_M$, $\mathbf{v} = \mathbb{1}_N$.
    \State $\mathbf{K} = \exp(-\frac{\mathbf{C}}{\varepsilon})$
    \While{Not Convergence}
    \State $\mathbf{u} = \mathbf{a} \oslash (\mathbf{K} \mathbf{v})$
    \State $\mathbf{v} = \mathbf{b} \oslash (\mathbf{K}^\top \mathbf{u})$
    \EndWhile
    \State $\mathbf{P} = \diag\mathbf{u} \cdot \mathbf{K} \cdot \diag\mathbf{v}$
    \Ensure Optimal Coupling $\mathbf{P}$
  \end{algorithmic}
\end{algorithm}

\begin{remark}[Convergence Condition]
  Since we have \cref{eqn:vanilla-sinkhorn-solution-multiply}, we can use those equations to check if the scaling variables
  $\mathbf{u}$ and $\mathbf{v}$ at the current iteration could approximate $\mathbf{a}$ and $\mathbf{b}$ well.
  That is, for some small $\rho > 0$ the algorithm terminates at iteration $L$, such that
  \begin{equation*}
    \begin{aligned}
      \lVert \mathbf{u}^{(L)} \odot \left(\mathbf{K} \mathbf{v}^{(L)}\right) - \mathbf{a}\rVert_2 \le \rho
      \quad\text{and}\quad
      \lVert \mathbf{v}^{(L)} \odot \left(\mathbf{K}^\top \mathbf{u}^{(L)}\right) - \mathbf{b}\rVert_2 \le \rho
    \end{aligned}
  \end{equation*}
  If the norms are smaller than some threshold $\rho$, we can then terminate the program.
\end{remark}

\subsection{Log-stabilized Sinkhorn}\label{subsec:log-sinkhorn}

As can be seen in \cref{eqn:optimal-coupling,eqn:vanilla-sinkhorn-solution-division},
the updates on $\mathbf{u}$ and $\mathbf{v}$ are in the exponential domain.
Therefore, if either $\mathbf{C}_{ij}$ is large or $\varepsilon$ is small,
the corresponding $\mathbf{K}_{ij}$ would be very large,
and the subsequent updates will fail because of the Inf/NaN values introduced by the numeric overflow.
This is commonly known as the numeric instability issue of the Sinkhorn algorithm \citep[Chapter 4.4]{peyre2019}.

To circumvent this instability issue, we could move the computation into the log-domain,
therefore without the need for computation of the exponential function.
To this end, instead of updating the $\mathbf{u}$ and $\mathbf{v}$ as in \cref{update:vanilla-sinkhorn},
we can consider the dual variables $\mathbf{f}$ and $\mathbf{g}$ in the Lagrangian in \cref{eqn:regularized-OT-lagrangian}.
Since $\mathbf{u} = \exp \left(\frac{\mathbf{f}}{\varepsilon}\right)$
and $\mathbf{v} = \exp \left(\frac{\mathbf{g}}\varepsilon\right)$, we have
\begin{equation}\label{eqn:fg-expression-as-uv}
  \begin{aligned}
    \mathbf{f} = \varepsilon \log(\mathbf{u})
    \quad\text{and}\quad
    \mathbf{g} = \varepsilon \log(\mathbf{v}).
  \end{aligned}
\end{equation}

Then, we can also rewrite the elements in the coupling matrix $\mathbf{P}$ as
\begin{equation*}
  \begin{aligned}
    \mathbf{P}_{ij}^{(\ell)}
    = \mathbf{u}_i^{(\ell)} \mathbf{K}_{ij} \mathbf{v}_j^{(\ell)}
    = \exp \left(
    -\frac{\mathbf{C}_{ij} - \mathbf{f}^{(\ell)}_i - \mathbf{g}^{(\ell)}_j}{\varepsilon}
    \right)
  \end{aligned}
\end{equation*}

and the matrix form of $\mathbf{P}$ expressed in $\mathbf{f}$ and $\mathbf{g}$ is
\begin{equation}\label{eqn:optimal-coupling-P-as-expR}
  \begin{aligned}
    \mathbf{P}^{(\ell)}
    = \exp \left(
    - \frac{\mathbf{C} - \mathbf{f}^{(\ell)} \cdot \mathbb{1}_N^\top -
      \mathbb{1}_M \cdot \mathbf{g}^{(\ell)\top}}{\varepsilon}
    \right)
    = \exp \left(
    -\frac{
        \mathbf{R} \left(\mathbf{f}^{(\ell)}, \mathbf{g}^{(\ell)}\right)
      }{\varepsilon}
    \right),
  \end{aligned}
\end{equation}

where we define\footnote{
  \cref{eqn:function-R-as-fg} can also be written as
  $\mathbf{R}_{i,j} = \mathbf{C}_{i,j} - \mathbf{f}_i - \mathbf{g}_j$,
  which is the same as in \citet[Remark 4.23]{peyre2019}.
}
\begin{equation}\label{eqn:function-R-as-fg}
  \begin{aligned}
    \mathbf{R} \left(\mathbf{f}^{(\ell)}, \mathbf{g}^{(\ell)}\right)
    =
    \mathbf{C} - \mathbf{f}^{(\ell)} \cdot \mathbb{1}_N^\top -
    \mathbb{1}_M \cdot \mathbf{g}^{(\ell)\top}.
  \end{aligned}
\end{equation}

Note that $\mathbf{f}$ and $\mathbf{g}$ are dual variables in the Lagrangian in \cref{eqn:regularized-OT-lagrangian}.
Then we can plug in $\mathbf{u}$ and $\mathbf{v}$ to get the following updating equations
\begin{equation}\label{eqn:log-sinkhorn-update1}
  \begin{aligned}
    \mathbf{f}^{(\ell)}
     & = \varepsilon\log \mathbf{u}^{(\ell)}
    = \varepsilon\log \mathbf{a} - \varepsilon\log (\mathbf{K} \mathbf{v}^{(\ell-1)}), \\
    \mathbf{g}^{(\ell)}
     & = \varepsilon\log \mathbf{v}^{(\ell)}
    = \varepsilon\log \mathbf{b} - \varepsilon\log (\mathbf{K}^\top \mathbf{u}^{(\ell)}).
  \end{aligned}
\end{equation}

Since $\mathbf{K}\mathbf{v}^{(\ell-1)} \in \mathbb{R}^m$, $\forall i$,
the $i$-th element is %
\begin{equation*}
  \begin{aligned}
    \left(\mathbf{K}\mathbf{v}^{(\ell-1)}\right)_i
     & = \sum_j \mathbf{K}_{ij} \mathbf{v}_j^{(\ell-1)}                                                                      \\
     & = \sum_j \exp(-\frac{\mathbf{C}_{ij}}{\varepsilon}) \exp(\frac{\mathbf{g}_j^{(\ell-1)}}{\varepsilon})                 \\
     & = \left[\sum_j \exp(- \frac{\mathbf{C}_{ij} - \mathbf{f}_i^{(\ell-1)} - \mathbf{g}_j^{(\ell-1)}}{\varepsilon})\right]
    \exp(-\frac{\mathbf{f}_i^{(\ell-1)}}{\varepsilon}),
  \end{aligned}
\end{equation*}

and
\begin{equation*}
  \begin{aligned}
    -\varepsilon\log \left(\mathbf{K} \mathbf{v}^{(\ell-1)}\right)_i
     & =
    -\varepsilon\log \sum_j \exp(- \frac{\mathbf{C}_{ij} - \mathbf{f}_i^{(\ell-1)} - \mathbf{g}_j^{(\ell-1)}}{\varepsilon})
    + \mathbf{f}_i^{(\ell-1)} \\
     & =
    -\varepsilon\log \left[
      \exp \left(
      -\frac{\mathbf{C}
        - \mathbf{f}^{(\ell-1)} \cdot \mathbb{1}_N^\top
        - \mathbb{1}_M \cdot (\mathbf{g}^{(\ell-1)})^\top}{\varepsilon}
      \right)
      \cdot \mathbb{1}_N
      \right]_i
    + \mathbf{f}_i^{(\ell-1)} \\
     & =
    -\varepsilon\log
    \left[
      \exp \left(
      -\frac{
        \mathbf{R} \left(\mathbf{f}^{(\ell-1)}, \mathbf{g}^{(\ell-1)}\right)
      }{\varepsilon}
      \right)
      \cdot \mathbb{1}_N
      \right]_i + \mathbf{f}_i^{(\ell-1)}.
  \end{aligned}
\end{equation*}

Therefore, we have
\begin{equation}\label{eqn:neg-log-kv-vec}
  \begin{aligned}
    - \varepsilon \log \left(\mathbf{K} \mathbf{v}^{(\ell-1)}\right)
    = - \varepsilon \log \left[
      \exp \left(
      -\frac{
        \mathbf{R} \left(\mathbf{f}^{(\ell-1)}, \mathbf{g}^{(\ell-1)}\right)
      }{\varepsilon}
      \right)
      \cdot \mathbb{1}_N
      \right] + \mathbf{f}^{(\ell-1)}.
  \end{aligned}
\end{equation}

And we have the updating equation for $\mathbf{f}^{(\ell)}$.
If we follow similar steps, we can also get it for $\mathbf{g}^{(\ell)}$.
Putting everything together,
\begin{equation}\label{eqn:update-fg-by-R}
  \begin{aligned}
    \mathbf{f}^{(\ell)}
     & = \mathbf{f}^{(\ell-1)} + \varepsilon\log \mathbf{a} -\varepsilon\log
    \left[
      \exp \left(
      -\frac{
        \mathbf{R} \left(\mathbf{f}^{(\ell-1)}, \mathbf{g}^{(\ell-1)}\right)
      }{\varepsilon}
      \right)
      \cdot \mathbb{1}_N
    \right],                                                                 \\
    \mathbf{g}^{(\ell)}
     & = \mathbf{g}^{(\ell-1)} + \varepsilon\log \mathbf{b} -\varepsilon\log
    \left[
      \exp \left(
      -\frac{
        \mathbf{R} \left(\mathbf{f}^{(\ell)}, \mathbf{g}^{(\ell-1)}\right)
      }{\varepsilon}
      \right)^\top
      \cdot \mathbb{1}_M
      \right].
  \end{aligned}
\end{equation}

Note that \cref{eqn:update-fg-by-R}
is essentially a ``log-sum-exp'' form, and we can apply the ``log-sum-exp'' trick\footnote{
  \url{https://en.wikipedia.org/wiki/LogSumExp}
} to make this equation numeric stable while computing it in practice.

To see this, first we need to write the \textit{soft-minimum} function \citep[Remark 4.22]{peyre2019},
let $\mathbf{z}$ be a vector of length $n$,
\begin{equation*}
  \begin{aligned}
    \min_{\varepsilon} \mathbf{z} = - \varepsilon\log \sum_i \exp \left(-\frac{\mathbf{z}_i}{\varepsilon}\right).
  \end{aligned}
\end{equation*}

Notice that when $\mathbf{z}_i$ is negative but large in absolute value,
the exponential and hence its sum will likely encounter numeric overflow,
which makes this computation not numerically stable.
To circumvent this, let use denote $\ubar{z} = \min \mathbf{z}$
and notice that the \textit{soft-minimum} is the same as
\begin{equation}\label{eqn:soft-minimum-stable}
  \begin{aligned}
    \min_\varepsilon \mathbf{z}
    = \ubar{z}
    - \varepsilon \log \sum_i \exp \left(-\frac{\mathbf{z}_i-\ubar{z}}{\varepsilon}\right).
  \end{aligned}
\end{equation}

This equivalent form guarantees that the maximum value of the vector before exponentiation
is 0 and hence remains numerically stable.
Then let us denote any matrix $\mathbf{R} \in \mathbb{R}^{M \times N}$
and $-\log \left[\exp \left(-\frac{\mathbf{R}}{\varepsilon} \right) \cdot \mathbb{1}_N\right]$
is equivalent to applying the \textit{soft-minimum} function to each row of the matrix $\mathbf{R}$.
Similarly, $-\log \left[\exp \left(-\frac{\mathbf{R}}{\varepsilon}\right)^\top \cdot \mathbb{1}_M\right]$
is equivalent to applying the \textit{soft-minimum} function to each column of the matrix $\mathbf{R}$.
Let us define the \textit{row-wise soft-minimum} and \textit{column-wise soft-minimum}
as in Equations (4.40) and (4.41) from \citet{peyre2019},
\begin{update}[Updating Equations for Log-Stabilized Sinkhron Algorithm]\label{update:update-fg-by-minrow-mincol}
  At any iteration $\ell = 1, \ldots, L$ and initialized $\mathbf{f}^{(0)}$ and $\mathbf{g}^{(0)}$, we have
  \begin{equation*}%
    \begin{aligned}
      \mathbf{f}^{(\ell)}
       & = \mathbf{f}^{(\ell-1)}
      + \varepsilon \log \mathbf{a}
      + \text{Min}^{row}_{\varepsilon} \,
      \mathbf{R} \left(\mathbf{f}^{(\ell-1)}, \mathbf{g}^{(\ell-1)}\right), \\
      \mathbf{g}^{(\ell)}
       & = \mathbf{g}^{(\ell-1)}
      + \varepsilon \log \mathbf{b}
      + \text{Min}^{col}_{\varepsilon} \,
      \mathbf{R} \left(\mathbf{f}^{(\ell)}, \mathbf{g}^{(\ell-1)}\right),
    \end{aligned}
  \end{equation*}
\end{update}

where for a matrix $\mathbf{R} \in \mathbb{R}^{M \times N}$ we have
\begin{equation*}
  \begin{aligned}
    \text{Min}^{row}_{\varepsilon} \left(\mathbf{R}\right)
     & = \left(
    \min_{\varepsilon} \left(\mathbf{R}_{i,j}\right)_j
    \right)_i \in \mathbb{R}^M, \\
    \text{Min}^{col}_{\varepsilon} \left(\mathbf{R}\right)
     & = \left(
    \min_{\varepsilon} \left(\mathbf{R}_{i,j}\right)_i
    \right)_j \in \mathbb{R}^N.
  \end{aligned}
\end{equation*}

In other words,
\begin{equation*}
  \begin{aligned}
    \text{Min}^{row}_{\varepsilon} \left(\mathbf{R}\right)
     & = \left[
      \min_{\varepsilon} \mathbf{R}_{1, \boldsymbol\cdot}\,,
      \ldots,\,
      \min_{\varepsilon} \mathbf{R}_{M, \boldsymbol\cdot}
    \right]^\top, \\
    \text{Min}^{col}_{\varepsilon} \left(\mathbf{R}\right)
     & = \left[
      \min_{\varepsilon} \mathbf{R}_{\boldsymbol\cdot, 1}\,,
      \ldots,\,
      \min_{\varepsilon} \mathbf{R}_{\boldsymbol\cdot, N}
      \right]^\top,
  \end{aligned}
\end{equation*}

where $\mathbf{R}_{i, \boldsymbol\cdot}$ and $\mathbf{R}_{\boldsymbol\cdot, j}$
denote the $i$-th row and $j$-th column of $\mathbf{R}$ respectively.

Thus, we have obtained the numerically stable version of the updating equations for $\mathbf{f}$ and $\mathbf{g}$.

\begin{remark}[]
  As mentioned in \cref{subsec:vanilla-sinkhorn}, the initialization for $\mathbf{u}$ and $\mathbf{v}$ are vectors of ones.
  This is to say that, the initializations of $\mathbf{f}$ and $\mathbf{g}$ are vectors of zeros,
  i.e. $\mathbf{f}^{(0)} = \mathbb{0}_M$ and $\mathbf{g}^{(0)} = \mathbb{0}_N$.
\end{remark}

Therefore, we can have the Log-Stabilized Sinkhorn algorithm.

\begin{algorithm}[H]
  \caption{Log-Stabilized Sinkhorn Algorithm}
  \begin{algorithmic}[1]\label{algo:log-sinkhorn}
    \Require $\mathbf{a} \in \Sigma_M$, $\mathbf{b} \in \Sigma_N$, $\mathbf{C} \in \mathbb{R}^{M\times N}$, $\varepsilon > 0$.
    \Initialize $\mathbf{f} = \mathbb{0}_M$, $\mathbf{g} = \mathbb{0}_N$.
    \While{Not Convergence}
    \State \# update $f$
    \State $\mathbf{R} = \mathbf{C} - \mathbf{f}\cdot \mathbb{1}_N^\top - \mathbb{1}_M\cdot \mathbf{g}^\top$
    \State $\mathbf{f} = \mathbf{f} + \varepsilon\log \mathbf{a}
      + \text{Min}^{row}_{\varepsilon} \, \mathbf{R}$
    \State
    \State \# update $g$
    \State $\mathbf{R} = \mathbf{C} - \mathbf{f}\cdot \mathbb{1}_N^\top - \mathbb{1}_M\cdot \mathbf{g}^\top$
    \State $\mathbf{g} = \mathbf{g} + \varepsilon\log \mathbf{b}
      + \text{Min}^{col}_{\varepsilon} \, \mathbf{R}$
    \EndWhile
    \State $\mathbf{P} = \exp \left(
      - \frac{\mathbf{C} - \mathbf{f}\cdot\mathbb{1}_N^\top - \mathbb{1}_M\cdot\mathbf{g}^\top}{\varepsilon}
      \right)$
    \Ensure $\mathbf{P}$
  \end{algorithmic}
\end{algorithm}

\begin{remark}[]
  The advantage of using log-stabilized Sinkhorn is to make sure the computation is stable for any arbitrary $\varepsilon$.
  The disadvantage, however, is that we cannot easily parallel its computation over several margins as in
  \cref{subsec:parallel-sinkhorn}.
  Moreover, the computation of matrix $\mathbf{Q}$ involves $\exp$ and $\log$ at every iteration twice
  for both updating $\mathbf{f}$ and $\mathbf{g}$.
  Thus, the log-stabilized Sinkhorn is less efficient than Vanilla and Parallel Sinkhorn algorithms.
  This is the tradeoff we need to make if we need to make sure the computation is stable,
  especially if we want to use a very small regularization constant $\varepsilon$.
  As in \citet[Chapter 4.1]{peyre2019},
  the Sinkhorn problem converges to original Kantorovich problem when $\varepsilon$ is small,
  i.e. $\ell_{\mathbf{C}}^\varepsilon(\mathbf{a},\mathbf{b}) \to \ell_{\mathbf{C}}(\mathbf{a},\mathbf{b})$
  as $\varepsilon \to 0$.
\end{remark}

\begin{remark}[Convergence Condition]\label{remark:conv-cond-parallel}
  Recall that in \cref{subsec:vanilla-sinkhorn},
  we can compare the marginals reconstructed by the scaling vectors versus the real marginals,
  and now we can do the same step albeit in the log domain.
  That is, for some $\rho > 0$ the algorithm terminates at iteration $L$, such that
  \begin{equation*}
    \begin{aligned}
      \left\lVert
      -\frac1{\varepsilon}\text{Min}^{row}_\varepsilon \, \mathbf{R} \left(\mathbf{f}^{(L)}, \mathbf{g}^{(L)}\right)
      - \log \mathbf{a}
      \right\rVert_2 & \le \rho, \\
      \left\lVert
      -\frac1{\varepsilon}\text{Min}^{col}_\varepsilon \, \mathbf{R} \left(\mathbf{f}^{(L)}, \mathbf{g}^{(L)}\right)
      - \log \mathbf{b}
      \right\rVert_2 & \le \rho.
    \end{aligned}
  \end{equation*}
\end{remark}

\subsection{Parallel Sinkhorn Algorithm}\label{subsec:parallel-sinkhorn}

Note that in the Vanilla Sinkhorn Algorithm, the computation are carried out by vector-vector element-wise computation,
i.e., $\mathbf{u}, \mathbf{a}, \text{and } \mathbf{K}\mathbf{v} \in \mathbb{R}^M$
and $\mathbf{v}, \mathbf{b}, \text{and } \mathbf{K}^\top \mathbf{u} \in \mathbb{R}^N$.
That is, if we have a single source $\mathbf{a}$ and a single target $\mathbf{b}$,
we are computing the ``one-to-one'' transportation.
However, sometimes we need to calculate one-to-many, many-to-one, or many-to-many transport
and this can be easily computed in a parallelized fashion.
All we need to do is to define the matrices $\mathbf{A}$ and $\mathbf{B}$ instead of the vectors $\mathbf{a}$ and $\mathbf{b}$,
and all that remains is to carry out the matrix computation instead of the vector ones.
Needless to say that the parallel algorithm does not immune to the numeric instability issue,
as discussed in the previous subsection \cref{subsec:log-sinkhorn},
since it's an extension to the original vanilla algorithm.

\begin{definition}
  For some $S \in \mathbb{N}_+$:
  \begin{itemize}
    \item Suppose we have $\mathbf{a} \in \Sigma_M$ and $\mathbf{b}_1, \ldots, \mathbf{b}_S \in \Sigma_N$.
          Then we can define
          \begin{equation*}
            \begin{aligned}
              \mathbf{A} = \left[\mathbf{a}, \ldots, \mathbf{a}\right] \in \mathbb{R}^{M \times S}
              \quad\text{and}\quad
              \mathbf{B} = \left[\mathbf{b}_1, \ldots, \mathbf{b}_S \right] \in \mathbb{R}^{N \times S}.
            \end{aligned}
          \end{equation*}
    \item Suppose we have $\mathbf{a}_1, \ldots, \mathbf{a}_S \in \Sigma_M$ and $\mathbf{b} \in \Sigma_N$.
          Then we can define
          \begin{equation*}
            \begin{aligned}
              \mathbf{A} = \left[\mathbf{a}_1, \ldots, \mathbf{a}_S\right] \in \mathbb{R}^{M \times S}
              \quad\text{and}\quad
              \mathbf{B} = \left[\mathbf{b}, \ldots, \mathbf{b} \right] \in \mathbb{R}^{N \times S}.
            \end{aligned}
          \end{equation*}
    \item Suppose we have $\mathbf{a}_1, \ldots, \mathbf{a}_S \in \Sigma_M$ and $\mathbf{b}_1, \ldots, \mathbf{b}_S \in \Sigma_N$.
          Then we can define
          \begin{equation*}
            \begin{aligned}
              \mathbf{A} = \left[\mathbf{a}_1, \ldots, \mathbf{a}_S\right] \in \mathbb{R}^{M \times S}
              \quad\text{and}\quad
              \mathbf{B} = \left[\mathbf{b}_1, \ldots, \mathbf{b}_S \right] \in \mathbb{R}^{N \times S}.
            \end{aligned}
          \end{equation*}
  \end{itemize}
\end{definition}

\begin{remark}[]
  It should be fairly obvious that the first two scenarios are special cases of the third,
  i.e.~many-to-many transport.
  In the many-to-many case, the number of columns in $\mathbf{A}$ should be equal to that of $\mathbf{B}$;
  whereas in one-to-many or many-to-one cases, we need to duplicate the single vector $S$ times,
  as if we had a $S$-column matrix to begin with.
  In many modern linear algebra packages, this can be achieved automatically by ``broadcasting,''\footnote{
    For example,
    Eigen (C++): \url{https://eigen.tuxfamily.org/dox/group__TutorialReductionsVisitorsBroadcasting.html},
    NumPy (Python): \url{https://numpy.org/doc/stable/user/basics.broadcasting.html}.
  }
  without the need of actually doing the duplication of the column vectors.
  Or in terms of matrix computation,
  for example,
  $\mathbf{B}
    = \left[\mathbf{b}, \ldots, \mathbf{b}\right]
    = \mathbb{1}_S^\top \otimes \mathbf{b}
    = \mathbf{b}\cdot \mathbb{1}_S^\top$.
  This will come into handy in \cref{sec:wasserstein-barycenter}.
\end{remark}

Therefore, we can have the matrix version\footnote{
  Instead of computing the vector-vector computation in \cref{eqn:vanilla-sinkhorn-solution-division} $S$ times,
  we are computing a matrix-matrix product/division once.
  They are mathematically equivalent, but in practice,
  the latter would be much faster due to the implementation of linear algebra routines (e.g. BLAS).
  Hence, this simple extension is considered parallel since we process $S$ columns of data in one batch.
} of \cref{eqn:vanilla-sinkhorn-solution-division} with scaling variable
$\mathbf{U} \in \mathbb{R}^{M\times S}$ and $\mathbf{V} \in \mathbb{R}^{N\times S}$,
\begin{equation}\label{eqn:update-UV-parallel}
  \begin{aligned}
    \mathbf{U} = \mathbf{A} \oslash (\mathbf{K}\mathbf{V})
    \quad\text{and}\quad
    \mathbf{V} = \mathbf{B} \oslash (\mathbf{K}^\top \mathbf{V}).
  \end{aligned}
\end{equation}

Hence, we could have the matrix version of the updating equations of Sinkhorn algorithm.
\begin{update}[Updating Equations for Parallel Sinkhorn Algorithm]\label{update:parallel-sinkhorn}
  At any iteration $\ell = 1, \ldots, L$ and initialized $\mathbf{V}^{(0)} \in \mathbb{R}^{N\times S}$, we have
  \begin{equation}
    \begin{aligned}
      \mathbf{U}^{(\ell)} = \mathbf{A} \oslash (\mathbf{K} \mathbf{V}^{(\ell-1)})
      \quad\text{and}\quad
      \mathbf{V}^{(\ell)} = \mathbf{B} \oslash (\mathbf{K}^\top \mathbf{U}^{(\ell)}).
    \end{aligned}
  \end{equation}
\end{update}

\begin{remark}[]
  Since \cref{update:parallel-sinkhorn} only involves matrix operations,
  the Sinkhorn algorithm are intrinsically GPU friendly, as matrix operations are even more efficient on GPU than on CPU.
  This why some applications of OT in machine learning leverage GPU acceleration via neural network libraries \citep{flamary2021}.
\end{remark}

\begin{algorithm}[H]
  \caption{Parallel Sinkhorn Algorithm}
  \begin{algorithmic}[1]\label{algo:parallel-sinkhorn}
    \Require $\mathbf{A} \in \Sigma_{M\times S}$, $\mathbf{B}\in \Sigma_{N\times S}$, $\mathbf{C} \in \mathbb{R}^{M\times N}$, $\varepsilon > 0$.
    \Initialize $\mathbf{U} = \mathbb{1}_{M \times S}$, $\mathbf{V} = \mathbb{1}_{N \times S}$.
    \State $\mathbf{K} = \exp(-\frac{\mathbf{C}}{\varepsilon})$
    \While{Not Convergence}

    \State $\mathbf{U} = \mathbf{A} \oslash \left(\mathbf{K} \mathbf{V}\right)$

    \State $\mathbf{V} = \mathbf{B} \oslash \left(\mathbf{K}^\top \mathbf{U}\right)$

    \EndWhile
    \State $\forall s, \mathbf{P}_s = \diag\mathbf{U}_s \cdot \mathbf{K} \cdot \diag\mathbf{V}_s$
    \Ensure Optimal Coupling $\mathbf{P}_s, \forall s$
  \end{algorithmic}
\end{algorithm}

\begin{remark}[Convergence Condition]
  Similar to the termination condition in \cref{subsec:vanilla-sinkhorn}, we now have the matrix norms to check convergence.
  For some $\rho > 0$ the algorithm terminates at iteration $L$, such that
  \begin{equation*}
    \begin{aligned}
      \lVert
      \mathbf{U}^{(L)} \odot \left(\mathbf{K} \mathbf{V}^{(L)}\right) - \mathbf{A}
      \rVert_2 \le \rho,
      \quad\text{and}\quad
      \lVert
      \mathbf{V}^{(L)} \odot \left(\mathbf{K}^\top \mathbf{U}^{(L)}\right) - \mathbf{B}
      \rVert_2 \le \rho.
    \end{aligned}
  \end{equation*}
\end{remark}

\section{Sinkhorn Algorithms with Gradients}\label{sec:sinkhorn-gradient}
\sectionmark{Sinkhorn Gradients}

In this section, I will derive the gradients of the loss function with respect to its arguments for the algorithms
described in \cref{sec:sinkhorn-algorithm}.
This is useful when one wishes to use regularized loss in \cref{eqn:entropic-regularized-OT-loss}
as a loss function for training purpose in machine learning, for example,
Wasserstein Dictionary Learning in \crefrange{sec:wasserstein-barycenter}{sec:wdl}.
Similar to the structure in \cref{sec:sinkhorn-algorithm},
I will derive the gradients for the vanilla and log-stabilized algorithms respectively.

\subsection{Gradients for Vanilla Sinkhorn}\label{subsec:gradient-vanilla-sinkhorn}

\begin{figure}[H]
  \centering

  \tikzset{every picture/.style={line width=0.75pt}} %

  \begin{tikzpicture}[x=0.75pt,y=0.75pt,yscale=-1,xscale=1]
    \draw    (207.77,55.83) -- (207.77,83.83) ;
    \draw [shift={(207.77,85.83)}, rotate = 270] [color={rgb, 255:red, 0; green, 0; blue, 0 }  ][line width=0.75]    (10.93,-3.29) .. controls (6.95,-1.4) and (3.31,-0.3) .. (0,0) .. controls (3.31,0.3) and (6.95,1.4) .. (10.93,3.29)   ;
    \draw    (271.77,55.83) -- (271.77,83.83) ;
    \draw [shift={(271.77,85.83)}, rotate = 270] [color={rgb, 255:red, 0; green, 0; blue, 0 }  ][line width=0.75]    (10.93,-3.29) .. controls (6.95,-1.4) and (3.31,-0.3) .. (0,0) .. controls (3.31,0.3) and (6.95,1.4) .. (10.93,3.29)   ;
    \draw    (394.1,55.83) -- (394.1,83.83) ;
    \draw [shift={(394.1,85.83)}, rotate = 270] [color={rgb, 255:red, 0; green, 0; blue, 0 }  ][line width=0.75]    (10.93,-3.29) .. controls (6.95,-1.4) and (3.31,-0.3) .. (0,0) .. controls (3.31,0.3) and (6.95,1.4) .. (10.93,3.29)   ;
    \draw    (207.77,110.83) -- (207.77,138.83) ;
    \draw [shift={(207.77,140.83)}, rotate = 270] [color={rgb, 255:red, 0; green, 0; blue, 0 }  ][line width=0.75]    (10.93,-3.29) .. controls (6.95,-1.4) and (3.31,-0.3) .. (0,0) .. controls (3.31,0.3) and (6.95,1.4) .. (10.93,3.29)   ;
    \draw    (394.1,113.83) -- (394.1,141.83) ;
    \draw [shift={(394.1,143.83)}, rotate = 270] [color={rgb, 255:red, 0; green, 0; blue, 0 }  ][line width=0.75]    (10.93,-3.29) .. controls (6.95,-1.4) and (3.31,-0.3) .. (0,0) .. controls (3.31,0.3) and (6.95,1.4) .. (10.93,3.29)   ;
    \draw    (159.1,138.5) -- (187.59,113.81) ;
    \draw [shift={(189.1,112.5)}, rotate = 139.09] [color={rgb, 255:red, 0; green, 0; blue, 0 }  ][line width=0.75]    (10.93,-3.29) .. controls (6.95,-1.4) and (3.31,-0.3) .. (0,0) .. controls (3.31,0.3) and (6.95,1.4) .. (10.93,3.29)   ;
    \draw    (348.1,138.5) -- (376.59,113.81) ;
    \draw [shift={(378.1,112.5)}, rotate = 139.09] [color={rgb, 255:red, 0; green, 0; blue, 0 }  ][line width=0.75]    (10.93,-3.29) .. controls (6.95,-1.4) and (3.31,-0.3) .. (0,0) .. controls (3.31,0.3) and (6.95,1.4) .. (10.93,3.29)   ;
    \draw    (477.43,121.5) -- (507.77,121.81) ;
    \draw [shift={(509.77,121.83)}, rotate = 180.59] [color={rgb, 255:red, 0; green, 0; blue, 0 }  ][line width=0.75]    (10.93,-3.29) .. controls (6.95,-1.4) and (3.31,-0.3) .. (0,0) .. controls (3.31,0.3) and (6.95,1.4) .. (10.93,3.29)   ;
    \draw    (417.77,145.17) -- (447.34,129.12) ;
    \draw [shift={(449.1,128.17)}, rotate = 151.52] [color={rgb, 255:red, 0; green, 0; blue, 0 }  ][line width=0.75]    (10.93,-3.29) .. controls (6.95,-1.4) and (3.31,-0.3) .. (0,0) .. controls (3.31,0.3) and (6.95,1.4) .. (10.93,3.29)   ;
    \draw    (417.1,106.17) -- (447.25,118.73) ;
    \draw [shift={(449.1,119.5)}, rotate = 202.62] [color={rgb, 255:red, 0; green, 0; blue, 0 }  ][line width=0.75]    (10.93,-3.29) .. controls (6.95,-1.4) and (3.31,-0.3) .. (0,0) .. controls (3.31,0.3) and (6.95,1.4) .. (10.93,3.29)   ;
    \draw    (228.77,138.5) -- (257.26,113.81) ;
    \draw [shift={(258.77,112.5)}, rotate = 139.09] [color={rgb, 255:red, 0; green, 0; blue, 0 }  ][line width=0.75]    (10.93,-3.29) .. controls (6.95,-1.4) and (3.31,-0.3) .. (0,0) .. controls (3.31,0.3) and (6.95,1.4) .. (10.93,3.29)   ;

    \draw (203.5,35.5) node [anchor=north west][inner sep=0.75pt]   [align=left] {$\displaystyle \mathbf{a}$};
    \draw (266.5,35.5) node [anchor=north west][inner sep=0.75pt]   [align=left] {$\displaystyle \mathbf{a}$};
    \draw (389.5,35.5) node [anchor=north west][inner sep=0.75pt]   [align=left] {$\displaystyle \mathbf{a}$};
    \draw (195,89) node [anchor=north west][inner sep=0.75pt]   [align=left] {$\displaystyle \mathbf{u}^{( 1)}$};
    \draw (383.5,89) node [anchor=north west][inner sep=0.75pt]   [align=left] {$\displaystyle \mathbf{u}^{( L)}$};
    \draw (258,89) node [anchor=north west][inner sep=0.75pt]   [align=left] {$\displaystyle \mathbf{u}^{( 2)}$};
    \draw (142,143.83) node [anchor=north west][inner sep=0.75pt]   [align=left] {$\displaystyle \mathbf{v}^{( 0)}$};
    \draw (195,143.83) node [anchor=north west][inner sep=0.75pt]   [align=left] {$\displaystyle \mathbf{v}^{( 1)}$};
    \draw (322,143.83) node [anchor=north west][inner sep=0.75pt]   [align=left] {$\displaystyle \mathbf{v}^{( L-1)}$};
    \draw (384,143.83) node [anchor=north west][inner sep=0.75pt]   [align=left] {$\displaystyle \mathbf{v}^{( L)}$};
    \draw (293.33,112.5) node [anchor=north west][inner sep=0.75pt]   [align=left] {$\displaystyle \mathbf{\dotsc }$};
    \draw (455.67,112.5) node [anchor=north west][inner sep=0.75pt]   [align=left] {$\displaystyle \mathbf{P}^{*}$};
    \draw (515,112.5) node [anchor=north west][inner sep=0.75pt]   [align=left] {$\displaystyle \mathcal{L}$};

  \end{tikzpicture}

  \caption{Computational graph for the vanilla Sinkhorn algorithm.}\label{fig:sinkhorn-vanilla-comp-graph}
\end{figure}

To recover the gradient of the loss $\mathcal{L}$ w.r.t.~$\mathbf{a}$,
we consider manually deriving the gradients by analyzing the computational graph as in \cref{fig:sinkhorn-vanilla-comp-graph}
and leverage the reverse mode automatic differentiation\footnote{
  Interested readers could refer to \citet[Chapter 8]{nocedal2006} for an introduction.
  Note that reverse mode AD is also the building block for most of the modern neural network libraries,
  though in this paper I derive the derivatives by hand, and thus eschewing the need for AD libraries.
}.

As supposed to the forward mode where we accumulates the gradients while doing the main computation,
in the reverse mode, we need to first carry out the main computation in the forward pass while recording
the intermediate values of the variables, and then at the end trace back the gradient
from the final loss.
For example, in the case of \cref{algo:vanilla-sinkhorn},
we can keep track of all the $\mathbf{u}^{(\ell)}$'s and $\mathbf{v}^{(\ell)}$'s,
i.e.~$\mathbf{u}^{(1)}$, $\ldots$, $\mathbf{u}^{(L)}$, and
$\mathbf{v}^{(1)}$, $\ldots$, $\mathbf{v}^{(L)}$.

To see the gradient of loss $\mathcal{L}$ w.r.t.~$\mathbf{a}$,
we can denote the \textit{adjoint} of $\mathbf{a}$ as $\widebar{\mathbf{a}}$,
and
\begin{equation*}
  \begin{aligned}
    \widebar{\mathbf{a}}
    = \sum_{l = 1}^{L} \left(\frac{\partial \mathbf{u}^{(\ell)}}{\partial \mathbf{a}}\right)^\top
    \cdot \frac{\partial \mathcal{L}}{\partial \mathbf{u}^{(\ell)}}
    = \sum_{l = 1}^{L} \left(\frac{\partial \mathbf{u}^{(\ell)}}{\partial \mathbf{a}}\right)^\top
    \cdot \widebar{\mathbf{u}}^{(\ell)}
    = \sum_{l = 1}^{L}
    \frac{\widebar{\mathbf{u}}^{(\ell)}}{\mathbf{K} \mathbf{v}^{(\ell-1)}}.
  \end{aligned}
\end{equation*}

Therefore, we will also need the adjoints $\widebar{\mathbf{u}}^{(\ell)}$, $\widebar{\mathbf{v}}^{(\ell)}$
First, to see the adjoints $\widebar{\mathbf{u}}^{(L)}$ and $\widebar{\mathbf{v}}^{(L)}$,
we have for $\ell = L$
\begin{equation*}
  \begin{aligned}
    \frac{\partial \mathcal{L}}{\partial \mathbf{u}^{(L)}_i}
     & = \sum_j \sum_k \frac{\partial \mathcal{L}}{\partial \mathbf{P}_{kj}}
    \cdot \frac{\partial \mathbf{P}_{kj}}{\partial \mathbf{u}^{(L)}_i}
    + \sum_j \frac{\partial \mathcal{L}}{\partial \mathbf{v}^{(L)}_j}
    \cdot \frac{\partial \mathbf{v}^{(L)}_j}{\partial \mathbf{u}^{(L)}_i}                          \\
     & = \sum_j \widebar{\mathbf{P}}_{ij} \cdot \mathbf{K}_{ij} \cdot \mathbf{v}^{(L)}_j
    - \sum_j \widebar{\mathbf{v}}^{(L)}_j \cdot
    \left(\frac{\mathbf{v}^{(L)}}{\mathbf{K}^\top \mathbf{u}^{(L)}}\right)_j \cdot \mathbf{K}_{ij} \\
     & = \left[
      \left(\widebar{\mathbf{P}} \odot \mathbf{K}\right)\cdot \mathbf{v}^{(L)} -
      \mathbf{K} \cdot \frac{\widebar{\mathbf{v}}^{(L)} \odot \mathbf{v}^{(L)}}{\mathbf{K}^\top \mathbf{u}^{(L)}}
    \right]_i,                                                                                     \\
    \frac{\partial \mathcal{L}}{\partial \mathbf{v}^{(L)}_j}
     & = \sum_i \frac{\partial \mathcal{L}}{\partial \mathbf{P}_{ij}} \cdot
    \frac{\partial \mathbf{P}_{ij}}{\partial \mathbf{v}^{(L)}_j}
    = \sum_i \mathbf{u}^{(L)}_i \cdot \mathbf{K}_{ij} \cdot \widebar{\mathbf{P}}_{ij}              \\
     & = \left[\left(\widebar{\mathbf{P}} \odot \mathbf{K}\right)^\top \mathbf{u}^{(L)}\right]_j,
  \end{aligned}
\end{equation*}

and
\begin{equation*}
  \begin{aligned}
    \widebar{\mathbf{u}}^{(L)}
     & = \left(\widebar{\mathbf{P}} \odot \mathbf{K}\right)\cdot \mathbf{v}^{(L)} -
    \mathbf{K} \cdot \frac{\widebar{\mathbf{v}}^{(L)} \odot \mathbf{v}^{(L)}}{\mathbf{K}^\top \mathbf{u}^{(L)}}, \\
    \widebar{\mathbf{v}}^{(L)}
     & = \left(\widebar{\mathbf{P}} \odot \mathbf{K}\right)^\top \mathbf{u}^{(L)}.
  \end{aligned}
\end{equation*}

Note that in the above display, we will also need the adjoint of the optimal coupling matrix $\mathbf{P}$.
To see this, we have from \cref{eqn:entropic-regularized-OT-loss} that
\begin{equation*}
  \begin{aligned}
    \frac{\partial \mathcal{L}}{\partial \mathbf{P}_{ij}}
     & = \frac{\partial }{\partial \mathbf{P}_{ij}}
    \left[
      \mathbf{P}_{ij} \cdot
      \left(\mathbf{C}_{ij} + \varepsilon\log \mathbf{P}_{ij} - \varepsilon\right)
    \right]                                                                      \\
     & = \mathbf{C}_{ij} + \varepsilon\log \mathbf{P}_{ij} - \varepsilon
    + \mathbf{P}_{ij} \cdot \frac{\partial }{\partial \mathbf{P}_{ij}}
    \left(\mathbf{C}_{ij} - \varepsilon\log \mathbf{P}_{ij} - \varepsilon\right) \\
     & = \mathbf{C}_{ij} + \varepsilon\log \mathbf{P}_{ij},
  \end{aligned}
\end{equation*}

and thus the adjoint of $\mathbf{P}$ is $\widebar{\mathbf{P}} = \mathbf{C} + \varepsilon \log \mathbf{P}$.
For all other $\ell = 1, \ldots, L-1$, we have
\begin{equation*}
  \begin{aligned}
    \frac{\partial \mathcal{L}}{\partial \mathbf{u}^{(\ell)}_i}
     & = \sum_j \frac{\partial \mathcal{L}}{\partial \mathbf{v}^{(\ell)}_j}
    \cdot \frac{\partial \mathbf{v}^{(\ell)}_j}{\partial \mathbf{u}^{(\ell)}_j}
    = \sum_j \widebar{\mathbf{v}}^{(\ell)}_j
    \cdot \frac{\partial \mathbf{v}^{(\ell)}_j}{\partial \mathbf{u}^{(\ell)}_j}, \\
    \frac{\partial \mathcal{L}}{\partial \mathbf{v}^{(\ell)}_j}
     & = \sum_i \frac{\partial \mathcal{L}}{\partial \mathbf{u}^{(\ell+1)}_i}
    \cdot \frac{\partial \mathbf{u}^{(\ell+1)}_i}{\partial \mathbf{v}^{(\ell)}_j}
    = \sum_i \widebar{\mathbf{u}}^{(\ell+1)}_i
    \cdot \frac{\partial \mathbf{u}^{(\ell+1)}_i}{\partial \mathbf{v}^{(\ell)}_j}.
  \end{aligned}
\end{equation*}

Now that we can easliy calculate the following
\begin{equation*}
  \begin{aligned}
    \frac{\partial \mathbf{v}^{(\ell)}_j}{\partial \mathbf{u}^{(\ell)}_j}
     & =
    - \left(
    \frac{\mathbf{v}^{(\ell)}}{\mathbf{K}^\top \mathbf{u}^{(\ell)}}
    \right)_j \cdot \mathbf{K}_{ij} \\
    \frac{\partial \mathbf{u}^{(\ell+1)}_i}{\partial \mathbf{v}^{(\ell)}_j}
     & =
    - \left(
    \frac{\mathbf{u}^{(\ell+1)}}{\mathbf{K} \mathbf{v}^{(\ell)}}
    \right)_i \cdot \mathbf{K}_{ij},
  \end{aligned}
\end{equation*}

and we could therefore have the adjoints
$\widebar{\mathbf{u}}^{(\ell)}$ and $\widebar{\mathbf{v}}^{(\ell)}$ for all $\ell = 1, \ldots, L-1$,
\begin{equation*}
  \begin{aligned}
    \widebar{\mathbf{u}}^{(\ell)}
     & = - \mathbf{K}\cdot
    \frac{\widebar{\mathbf{v}}^{(\ell)} \odot \mathbf{v}^{(\ell)}}{\mathbf{K}^\top \mathbf{u}^{(\ell)}}, \\
    \widebar{\mathbf{v}}^{(\ell)}
     & = - \mathbf{K}^\top \cdot
    \frac{\widebar{\mathbf{u}}^{(\ell+1)} \odot \mathbf{u}^{(\ell+1)}}{\mathbf{K} \mathbf{v}^{(\ell)}}.
  \end{aligned}
\end{equation*}

\begin{algorithm}[H]
  \caption{Vanilla Sinkhorn Algorithm with Gradient}
  \begin{algorithmic}[1]\label{algo:vanillia-sinkhorn-with-gradient}
    \Require $\mathbf{a} \in \Sigma_M$, $\mathbf{b} \in \Sigma_N$, $\mathbf{C} \in \mathbb{R}^{M\times N}$, $\varepsilon > 0$.
    \Initialize $\mathbf{u} = \mathbb{1}_M$, $\mathbf{v} = \mathbb{1}_N$.
    \State $\mathbf{K} = \exp(-\frac{\mathbf{C}}{\varepsilon})$

    \State \# forward loop
    \While{$\ell = 1, \ldots, L$}
    \State $\mathbf{u} = \mathbf{a} \oslash (\mathbf{K} \mathbf{v})$         \# update u
    \State $\mathbf{v} = \mathbf{b} \oslash (\mathbf{K}^\top \mathbf{u})$    \# update v
    \EndWhile

    \State \# backward loop
    \State $\mathbf{P} = \diag \mathbf{u} \cdot \mathbf{K} \cdot \diag \mathbf{v}$
    \State $\widebar{\mathbf{P}} = \mathbf{C} + \varepsilon \log \mathbf{P}$
    \State $\widebar{\mathbf{v}} = \left(\widebar{\mathbf{P}}\odot\mathbf{K}\right)^\top \mathbf{u}$
    \State $\widebar{\mathbf{u}} = \left(\widebar{\mathbf{P}}\odot\mathbf{K}\right) \mathbf{v} -
      \mathbf{K} \cdot \frac{\widebar{\mathbf{v}} \odot \mathbf{v}}{\mathbf{K}^\top \mathbf{u}}$
    \State $\widebar{\mathbf{a}} = \widebar{\mathbf{u}} \oslash \left(\mathbf{K} \mathbf{v}^{(L-1)}\right)$
    \For{$\ell = L-1, \ldots, 1$}
    \State $\widebar{\mathbf{v}} = - \mathbf{K}^\top\cdot
      \frac{\widebar{\mathbf{u}} \odot \mathbf{u}^{(\ell+1)}}{\mathbf{K}^\top \mathbf{v}^{(\ell)}}$
    \State $\widebar{\mathbf{u}} = - \mathbf{K} \cdot
      \frac{\widebar{\mathbf{v}} \odot \mathbf{v}^{(\ell)}}{\mathbf{K}^\top \mathbf{u}^{(\ell)}}$
    \State $\widebar{\mathbf{a}} = \widebar{\mathbf{a}} + \widebar{\mathbf{u}} \oslash \left(\mathbf{K} \mathbf{v}^{(\ell-1)}\right)$
    \EndFor

    \Ensure $\mathbf{P}$, $\widebar{\mathbf{a}}$
  \end{algorithmic}
\end{algorithm}

\subsection{Gradients for Log-Stabilized Sinkhorn}\label{subsec:gradient-log-sinkhorn}

\begin{figure}[H]
  \centering

  \tikzset{every picture/.style={line width=0.75pt}} %

  \begin{tikzpicture}[x=0.75pt,y=0.75pt,yscale=-1,xscale=1]
    \draw    (207.77,55.83) -- (207.77,83.83) ;
    \draw [shift={(207.77,85.83)}, rotate = 270] [color={rgb, 255:red, 0; green, 0; blue, 0 }  ][line width=0.75]    (10.93,-3.29) .. controls (6.95,-1.4) and (3.31,-0.3) .. (0,0) .. controls (3.31,0.3) and (6.95,1.4) .. (10.93,3.29)   ;
    \draw    (271.77,55.83) -- (271.77,83.83) ;
    \draw [shift={(271.77,85.83)}, rotate = 270] [color={rgb, 255:red, 0; green, 0; blue, 0 }  ][line width=0.75]    (10.93,-3.29) .. controls (6.95,-1.4) and (3.31,-0.3) .. (0,0) .. controls (3.31,0.3) and (6.95,1.4) .. (10.93,3.29)   ;
    \draw    (394.1,55.83) -- (394.1,83.83) ;
    \draw [shift={(394.1,85.83)}, rotate = 270] [color={rgb, 255:red, 0; green, 0; blue, 0 }  ][line width=0.75]    (10.93,-3.29) .. controls (6.95,-1.4) and (3.31,-0.3) .. (0,0) .. controls (3.31,0.3) and (6.95,1.4) .. (10.93,3.29)   ;
    \draw    (207.77,110.83) -- (207.77,138.83) ;
    \draw [shift={(207.77,140.83)}, rotate = 270] [color={rgb, 255:red, 0; green, 0; blue, 0 }  ][line width=0.75]    (10.93,-3.29) .. controls (6.95,-1.4) and (3.31,-0.3) .. (0,0) .. controls (3.31,0.3) and (6.95,1.4) .. (10.93,3.29)   ;
    \draw    (394.1,113.83) -- (394.1,141.83) ;
    \draw [shift={(394.1,143.83)}, rotate = 270] [color={rgb, 255:red, 0; green, 0; blue, 0 }  ][line width=0.75]    (10.93,-3.29) .. controls (6.95,-1.4) and (3.31,-0.3) .. (0,0) .. controls (3.31,0.3) and (6.95,1.4) .. (10.93,3.29)   ;
    \draw    (159.1,138.5) -- (187.59,113.81) ;
    \draw [shift={(189.1,112.5)}, rotate = 139.09] [color={rgb, 255:red, 0; green, 0; blue, 0 }  ][line width=0.75]    (10.93,-3.29) .. controls (6.95,-1.4) and (3.31,-0.3) .. (0,0) .. controls (3.31,0.3) and (6.95,1.4) .. (10.93,3.29)   ;
    \draw    (348.1,138.5) -- (376.59,113.81) ;
    \draw [shift={(378.1,112.5)}, rotate = 139.09] [color={rgb, 255:red, 0; green, 0; blue, 0 }  ][line width=0.75]    (10.93,-3.29) .. controls (6.95,-1.4) and (3.31,-0.3) .. (0,0) .. controls (3.31,0.3) and (6.95,1.4) .. (10.93,3.29)   ;
    \draw    (477.43,121.5) -- (507.77,121.81) ;
    \draw [shift={(509.77,121.83)}, rotate = 180.59] [color={rgb, 255:red, 0; green, 0; blue, 0 }  ][line width=0.75]    (10.93,-3.29) .. controls (6.95,-1.4) and (3.31,-0.3) .. (0,0) .. controls (3.31,0.3) and (6.95,1.4) .. (10.93,3.29)   ;
    \draw    (417.77,145.17) -- (447.34,129.12) ;
    \draw [shift={(449.1,128.17)}, rotate = 151.52] [color={rgb, 255:red, 0; green, 0; blue, 0 }  ][line width=0.75]    (10.93,-3.29) .. controls (6.95,-1.4) and (3.31,-0.3) .. (0,0) .. controls (3.31,0.3) and (6.95,1.4) .. (10.93,3.29)   ;
    \draw    (417.1,106.17) -- (447.25,118.73) ;
    \draw [shift={(449.1,119.5)}, rotate = 202.62] [color={rgb, 255:red, 0; green, 0; blue, 0 }  ][line width=0.75]    (10.93,-3.29) .. controls (6.95,-1.4) and (3.31,-0.3) .. (0,0) .. controls (3.31,0.3) and (6.95,1.4) .. (10.93,3.29)   ;
    \draw    (228.77,138.5) -- (257.26,113.81) ;
    \draw [shift={(258.77,112.5)}, rotate = 139.09] [color={rgb, 255:red, 0; green, 0; blue, 0 }  ][line width=0.75]    (10.93,-3.29) .. controls (6.95,-1.4) and (3.31,-0.3) .. (0,0) .. controls (3.31,0.3) and (6.95,1.4) .. (10.93,3.29)   ;
    \draw    (163.4,158.15) -- (187.5,158.24) ;
    \draw [shift={(189.5,158.25)}, rotate = 180.2] [color={rgb, 255:red, 0; green, 0; blue, 0 }  ][line width=0.75]    (10.93,-3.29) .. controls (6.95,-1.4) and (3.31,-0.3) .. (0,0) .. controls (3.31,0.3) and (6.95,1.4) .. (10.93,3.29)   ;
    \draw    (349.4,158.15) -- (373.5,158.24) ;
    \draw [shift={(375.5,158.25)}, rotate = 180.2] [color={rgb, 255:red, 0; green, 0; blue, 0 }  ][line width=0.75]    (10.93,-3.29) .. controls (6.95,-1.4) and (3.31,-0.3) .. (0,0) .. controls (3.31,0.3) and (6.95,1.4) .. (10.93,3.29)   ;
    \draw    (222.6,100.41) -- (246.7,100.49) ;
    \draw [shift={(248.7,100.5)}, rotate = 180.2] [color={rgb, 255:red, 0; green, 0; blue, 0 }  ][line width=0.75]    (10.93,-3.29) .. controls (6.95,-1.4) and (3.31,-0.3) .. (0,0) .. controls (3.31,0.3) and (6.95,1.4) .. (10.93,3.29)   ;
    \draw    (351.8,100.41) -- (375.9,100.49) ;
    \draw [shift={(377.9,100.5)}, rotate = 180.2] [color={rgb, 255:red, 0; green, 0; blue, 0 }  ][line width=0.75]    (10.93,-3.29) .. controls (6.95,-1.4) and (3.31,-0.3) .. (0,0) .. controls (3.31,0.3) and (6.95,1.4) .. (10.93,3.29)   ;
    \draw    (282.6,100.41) -- (306.7,100.49) ;
    \draw [shift={(308.7,100.5)}, rotate = 180.2] [color={rgb, 255:red, 0; green, 0; blue, 0 }  ][line width=0.75]    (10.93,-3.29) .. controls (6.95,-1.4) and (3.31,-0.3) .. (0,0) .. controls (3.31,0.3) and (6.95,1.4) .. (10.93,3.29)   ;
    \draw    (225,158.15) -- (249.1,158.24) ;
    \draw [shift={(251.1,158.25)}, rotate = 180.2] [color={rgb, 255:red, 0; green, 0; blue, 0 }  ][line width=0.75]    (10.93,-3.29) .. controls (6.95,-1.4) and (3.31,-0.3) .. (0,0) .. controls (3.31,0.3) and (6.95,1.4) .. (10.93,3.29)   ;
    \draw    (290.6,158.15) -- (314.7,158.24) ;
    \draw [shift={(316.7,158.25)}, rotate = 180.2] [color={rgb, 255:red, 0; green, 0; blue, 0 }  ][line width=0.75]    (10.93,-3.29) .. controls (6.95,-1.4) and (3.31,-0.3) .. (0,0) .. controls (3.31,0.3) and (6.95,1.4) .. (10.93,3.29)   ;

    \draw (203.5,35.5) node [anchor=north west][inner sep=0.75pt]   [align=left] {$\displaystyle \mathbf{a}$};
    \draw (266.5,35.5) node [anchor=north west][inner sep=0.75pt]   [align=left] {$\displaystyle \mathbf{a}$};
    \draw (389.5,35.5) node [anchor=north west][inner sep=0.75pt]   [align=left] {$\displaystyle \mathbf{a}$};
    \draw (195,89) node [anchor=north west][inner sep=0.75pt]   [align=left] {$\displaystyle \mathbf{f}^{( 1)}$};
    \draw (383.5,89) node [anchor=north west][inner sep=0.75pt]   [align=left] {$\displaystyle \mathbf{f}^{( L)}$};
    \draw (258,89) node [anchor=north west][inner sep=0.75pt]   [align=left] {$\displaystyle \mathbf{f}^{( 2)}$};
    \draw (142,143.83) node [anchor=north west][inner sep=0.75pt]   [align=left] {$\displaystyle \mathbf{g}^{( 0)}$};
    \draw (194.33,143.83) node [anchor=north west][inner sep=0.75pt]   [align=left] {$\displaystyle \mathbf{g}^{( 1)}$};
    \draw (322,143.83) node [anchor=north west][inner sep=0.75pt]   [align=left] {$\displaystyle \mathbf{g}^{( L-1)}$};
    \draw (384,143.83) node [anchor=north west][inner sep=0.75pt]   [align=left] {$\displaystyle \mathbf{g}^{( L)}$};
    \draw (291.33,125.5) node [anchor=north west][inner sep=0.75pt]   [align=left] {$\displaystyle \mathbf{\dotsc }$};
    \draw (455.67,112.5) node [anchor=north west][inner sep=0.75pt]   [align=left] {$\displaystyle \mathbf{P}^{*}$};
    \draw (515,112.5) node [anchor=north west][inner sep=0.75pt]   [align=left] {$\displaystyle \mathcal{L}$};
    \draw (318.73,97.9) node [anchor=north west][inner sep=0.75pt]   [align=left] {$\displaystyle \mathbf{\dotsc }$};
    \draw (258.13,155.7) node [anchor=north west][inner sep=0.75pt]   [align=left] {$\displaystyle \mathbf{\dotsc }$};

  \end{tikzpicture}

  \caption{Computational graph for the log-stabilized Sinkhorn algorithm.}\label{fig:sinkhorn-log-comp-graph}
\end{figure}

Similar to \cref{subsec:gradient-vanilla-sinkhorn},
we have the computational graph for the log-stabilized Sinkhorn algorithm in \cref{fig:sinkhorn-log-comp-graph}.
Again, we need to separate the cases for $\ell = L$ and $\ell = 1, \ldots, L-1$.
First, for $\ell = L$, we have
\begin{equation*}
  \begin{aligned}
    \frac{\partial \mathcal{L}}{\partial \mathbf{f}^{(L)}_i}
     & = \sum_j \sum_k \frac{\partial \mathcal{L}}{\partial \mathbf{P}_{kj}} \cdot
    \frac{\partial \mathbf{P}_{kj}}{\partial \mathbf{f}^{(L)}_i}
    + \sum_j \frac{\partial \mathcal{L}}{\partial \mathbf{g}^{(L)}_j} \cdot
    \frac{\partial \mathbf{g}^{(L)}_j}{\partial \mathbf{f}^{(L)}_i},               \\
    \frac{\partial \mathcal{L}}{\partial \mathbf{g}^{(L)}_j}
     & =
    \sum_i \frac{\partial \mathcal{L}}{\partial \mathbf{P}_{ij}} \cdot
    \frac{\partial \mathbf{P}_{ij}}{\partial \mathbf{g}^{(L)}_j}
    = \sum_i \widebar{\mathbf{P}}_{ij} \cdot
    \frac{\partial \mathbf{P}_{ij}}{\partial \mathbf{g}^{(L)}_j}
    = \frac1\varepsilon \sum_i \widebar{\mathbf{P}}_{ij} \cdot \mathbf{P}_{ij}
    = \frac1\varepsilon \left[\left(\widebar{\mathbf{P}} \odot \mathbf{P}\right)^\top \cdot \mathbb{1}_M\right]_i
  \end{aligned}
\end{equation*}

Now that $\frac{\partial \mathbf{P}_{kj}}{\partial \mathbf{f}^{(L)}_i}
  = \frac1\varepsilon \mathbf{P}_{kj} \cdot \frac{\partial \mathbf{f}^{(L)}_k}{\partial \mathbf{f}^{(L)}_i}$,
and thus
\begin{equation*}
  \begin{aligned}
    \sum_j \sum_k \frac{\partial \mathcal{L}}{\partial \mathbf{P}_{kj}} \cdot
    \frac{\partial \mathbf{P}_{kj}}{\partial \mathbf{f}^{(L)}_i}
    = \sum_j \widebar{\mathbf{P}}_{ij} \cdot \frac1\varepsilon \mathbf{P}_{ij}
    = \frac1\varepsilon \left[\left(\widebar{\mathbf{P}} \odot \mathbf{P}\right) \cdot \mathbb{1}_N\right]_i.
  \end{aligned}
\end{equation*}

Moreover, we have
\begin{equation*}
  \begin{aligned}
    \frac{\partial \mathbf{g}^{(L)}_j}{\partial \mathbf{f}^{(L)}_i}
     & = \frac{\partial }{\partial \mathbf{f}^{(L)}_i}
    \left[
      \mathbf{g}^{(L-1)}_j + \varepsilon \log \mathbf{b}_j +
      \left(
      \text{Min}_\varepsilon^{col} \mathbf{R} \left(\mathbf{f}^{(L)}, \mathbf{g}^{(L-1)}\right)
      \right)_j
    \right]                                                                                            \\
     & = \frac{\partial }{\partial \mathbf{f}^{(L)}_i}
    \min_\varepsilon \mathbf{R} \left(\mathbf{f}^{(L)}, \mathbf{g}^{(L-1)}\right)_{\boldsymbol\cdot,j} \\
     & = \sum_k
    \frac{
      \partial \min_\varepsilon \mathbf{R} \left(\mathbf{f}^{(L)}, \mathbf{g}^{(L-1)}\right)_{\boldsymbol\cdot,j}
    }{
      \partial \mathbf{R} \left(\mathbf{f}^{(L)}, \mathbf{g}^{(L-1)}\right)_{k,j}
    }
    \cdot
    \frac{
      \partial \mathbf{R} \left(\mathbf{f}^{(L)}, \mathbf{g}^{(L-1)}\right)_{k,j}
    }{
      \partial \mathbf{f}^{(L)}_i
    }
     & = -
    \frac{
      \partial \min_\varepsilon \mathbf{R} \left(\mathbf{f}^{(L)}, \mathbf{g}^{(L-1)}\right)_{\boldsymbol\cdot,j}
    }{
      \partial \mathbf{R} \left(\mathbf{f}^{(L)}, \mathbf{g}^{(L-1)}\right)_{i,j}
    }
  \end{aligned}
\end{equation*}

and by \cref{lemma:grad-soft-minimum} we have
\begin{equation*}
  \begin{aligned}
    \frac{
      \partial \mathbf{R} \left(\mathbf{f}^{(L)}, \mathbf{g}^{(L-1)}\right)_{\boldsymbol\cdot,j}
    }{
      \partial \mathbf{R} \left(\mathbf{f}^{(L)}, \mathbf{g}^{(L-1)}\right)_{i,j}
    }
    =
    \frac{
      \exp \left(-\frac{
        \mathbf{R} \left(\mathbf{f}^{(L)}, \mathbf{g}^{(L-1)}\right)_{i,j} -
        \min \mathbf{R} \left(\mathbf{f}^{(L)}, \mathbf{g}^{(L-1)}\right)_{\boldsymbol\cdot,j}
      }{\varepsilon}\right)
    }{
      \sum_{i'}
      \exp \left(-\frac{
        \mathbf{R} \left(\mathbf{f}^{(L)}, \mathbf{g}^{(L-1)}\right)_{i',j} -
        \min \mathbf{R} \left(\mathbf{f}^{(L)}, \mathbf{g}^{(L-1)}\right)_{\boldsymbol\cdot,j}
      }{\varepsilon}\right)
    }
    = \mathbf{W} \left(\mathbf{f}^{(L)}, \mathbf{g}^{(L-1)}\right)_{i,j},
  \end{aligned}
\end{equation*}

where %
\begin{equation*}
  \begin{aligned}
    \mathbf{W} \left(\mathbf{f}^{(L)}, \mathbf{g}^{(L-1)}\right)
    = \begin{bmatrix*}
        \ldots, &
        \nabla_{
          \mathbf{R} \left(\mathbf{f}^{(L)}, \mathbf{g}^{(L-1)}\right)_{\boldsymbol\cdot,j}
        } \min_\varepsilon \mathbf{R} \left(\mathbf{f}^{(L)}, \mathbf{g}^{(L-1)}\right)_{\boldsymbol\cdot,j},&
        \ldots
      \end{bmatrix*} \in \mathbb{R}^{M \times N}.
  \end{aligned}
\end{equation*}

Thus,
\begin{equation*}
  \begin{aligned}
    \frac{\partial \mathcal{L}}{\partial \mathbf{f}^{(L)}_i}
    = \frac1\varepsilon \left[\left(\widebar{\mathbf{P}} \odot \mathbf{P}\right) \cdot \mathbb{1}_N\right]_i
    - \left[\mathbf{W} \left(\mathbf{f}^{(L)}, \mathbf{g}^{(L-1)}\right) \cdot \widebar{\mathbf{g}}^{(L)}\right]_i,
  \end{aligned}
\end{equation*}

and
\begin{equation*}
  \begin{aligned}
    \widebar{\mathbf{f}}^{(L)}
     & = \frac{\partial \mathcal{L}}{\partial \mathbf{f}^{(L)}}
    = \frac1\varepsilon \left(\widebar{\mathbf{P}} \odot \mathbf{P}\right) \cdot \mathbb{1}_N
    - \mathbf{W} \left(\mathbf{f}^{(L)}, \mathbf{g}^{(L-1)}\right) \cdot \widebar{\mathbf{g}}^{(L)}, \\
    \widebar{\mathbf{g}}^{(L)}
     & = \frac{\partial \mathcal{L}}{\partial \mathbf{g}^{(L)}}
    = \frac1\varepsilon \left(\widebar{\mathbf{P}} \odot \mathbf{P}\right)^\top \cdot \mathbb{1}_M.
  \end{aligned}
\end{equation*}

Same as in \cref{subsec:gradient-vanilla-sinkhorn},
we can have the adjoint of $\mathbf{P}$ as $\widebar{\mathbf{P}}$.
Then for all $\ell = 1, \ldots, L-1$, we have
\begin{equation*}
  \begin{aligned}
    \frac{\partial \mathcal{L}}{\partial \mathbf{f}^{(\ell)}_i}
     & =
    \frac{\partial \mathcal{L}}{\partial \mathbf{f}^{(\ell+1)}_i}\cdot
    \frac{\partial \mathbf{f}^{(\ell+1)}_i}{\partial \mathbf{f}^{(\ell)}_i}
    + \sum_j \frac{\partial \mathcal{L}}{\partial \mathbf{g}^{(\ell)}_j} \cdot
    \frac{\partial \mathbf{g}^{(\ell)}_j}{\partial \mathbf{f}^{(\ell)}_i}
    = \widebar{\mathbf{f}}^{(\ell+1)} \cdot
    \frac{\partial \mathbf{f}^{(\ell+1)}_i}{\partial \mathbf{f}^{(\ell)}_i}
    + \sum_j \widebar{\mathbf{g}}^{(\ell)}_j \cdot
    \frac{\partial \mathbf{g}^{(\ell)}_j}{\partial \mathbf{f}^{(\ell)}_i}, \\
    \frac{\partial \mathcal{L}}{\partial \mathbf{g}^{(\ell)}_j}
     & =
    \frac{\partial \mathcal{L}}{\partial \mathbf{g}^{(\ell+1)}_j}\cdot
    \frac{\partial \mathbf{g}^{(\ell+1)}_j}{\partial \mathbf{g}^{(\ell)}_j}
    + \sum_i \frac{\partial \mathcal{L}}{\partial \mathbf{f}^{(\ell+1)}_i}\cdot
    \frac{\partial \mathbf{f}^{(\ell+1)}_i}{\partial \mathbf{g}^{(\ell)}_j}
    = \widebar{\mathbf{g}}^{(\ell)}_j \cdot
    \frac{\partial \mathbf{g}^{(\ell+1)}_j}{\partial \mathbf{g}^{(\ell)}_j}
    + \sum_i \widebar{\mathbf{f}}^{(\ell+1)}_i \cdot
    \frac{\partial \mathbf{f}^{(\ell+1)}_i}{\partial \mathbf{g}^{(\ell)}_j}.
  \end{aligned}
\end{equation*}

Since we also have
\begin{equation*}
  \begin{aligned}
    \frac{\partial \mathbf{f}^{(\ell+1)}_i}{\partial \mathbf{f}^{(\ell)}_i}
     & = \frac{\partial }{\partial \mathbf{f}^{(\ell)}_i}
    \left[
      \mathbf{f}^{(\ell)}_i + \varepsilon\log \mathbf{a}_i +
      \min_\varepsilon \mathbf{R} \left(\mathbf{f}^{(\ell)}, \mathbf{g}^{(\ell)}\right)_{i,\boldsymbol\cdot}
    \right]                                               \\
     & = 1 + \frac{
      \partial \min_\varepsilon \mathbf{R} \left(\mathbf{f}^{(\ell)}, \mathbf{g}^{(\ell)}\right)_{i,\boldsymbol\cdot}
    }{\partial \mathbf{f}^{(\ell)}_i}                     \\
     & = 1 + \sum_j
    \frac{
      \partial \min_\varepsilon \mathbf{R} \left(\mathbf{f}^{(\ell)}, \mathbf{g}^{(\ell)}\right)_{i,\boldsymbol\cdot}
    }{
      \partial \mathbf{R} \left(\mathbf{f}^{(\ell)}, \mathbf{g}^{(\ell)}\right)_{i,j}
    } \cdot
    \frac{
      \partial \mathbf{R} \left(\mathbf{f}^{(\ell)}, \mathbf{g}^{(\ell)}\right)_{i,j}
    }{
      \partial \mathbf{f}^{(\ell)}_i
    }                                                     \\
     & = 1 -
    \sum_j \frac{
      \exp \left(
      - \frac{
        \mathbf{R} \left(\mathbf{f}^{(\ell)}, \mathbf{g}^{(\ell)}\right)_{i,j}
        - \min \mathbf{R} \left(\mathbf{f}^{(\ell)}, \mathbf{g}^{(\ell)}\right)_{i,\boldsymbol\cdot}
      }{\varepsilon}
      \right)
    }{
      \sum_{j'} \exp \left(
      - \frac{
        \mathbf{R} \left(\mathbf{f}^{(\ell)}, \mathbf{g}^{(\ell)}\right)_{i,j'}
        - \min \mathbf{R} \left(\mathbf{f}^{(\ell)}, \mathbf{g}^{(\ell)}\right)_{i,\boldsymbol\cdot}
      }{\varepsilon}
      \right)
    }                                                     \\
     & = 0,
  \end{aligned}
\end{equation*}

and similarly,
\begin{equation*}
  \begin{aligned}
    \frac{\partial \mathbf{g}^{(\ell+1)}_j}{\partial \mathbf{g}^{(\ell)}_j} = 0.
  \end{aligned}
\end{equation*}

Also similar to the derivation for the $\ell = L$ case, we have
\begin{equation*}
  \begin{aligned}
    \frac{\partial \mathbf{g}^{(\ell)}_j}{\partial \mathbf{f}^{(\ell)}_i}
     & = - \mathbf{W} \left(\mathbf{f}^{(\ell)}, \mathbf{g}^{(\ell-1)}\right)_{i,j}, \\
    \frac{\partial \mathbf{f}^{(\ell+1)}_i}{\partial \mathbf{g}^{(\ell)}_j}
     & = - \mathbf{X} \left(\mathbf{f}^{(\ell)}, \mathbf{g}^{(\ell)}\right)_{j,i},
  \end{aligned}
\end{equation*}

where
\begin{equation}\label{eqn:W-X-log-sinkhorn}
  \begin{aligned}
    \mathbf{W} \left(\mathbf{f}^{(\ell)}, \mathbf{g}^{(\ell-1)}\right)
     & = \begin{bmatrix*}
           \ldots, &
           \nabla_{
             \mathbf{R} \left(\mathbf{f}^{(\ell)}, \mathbf{g}^{(\ell-1)}\right)_{\boldsymbol\cdot,j}
           } \min_\varepsilon \mathbf{R} \left(\mathbf{f}^{(\ell)}, \mathbf{g}^{(\ell-1)}\right)_{\boldsymbol\cdot,j},&
           \ldots
         \end{bmatrix*}, \\
    \mathbf{X} \left(\mathbf{f}^{(\ell)}, \mathbf{g}^{(\ell)}\right)
     & = \begin{bmatrix*}
           \ldots, &
           \nabla_{
             \mathbf{R} \left(\mathbf{f}^{(\ell)}, \mathbf{g}^{(\ell)}\right)_{i,\boldsymbol\cdot}
           } \min_\varepsilon \mathbf{R} \left(\mathbf{f}^{(\ell)}, \mathbf{g}^{(\ell)}\right)_{i,\boldsymbol\cdot},&
           \ldots
         \end{bmatrix*}.
  \end{aligned}
\end{equation}

Therefore, for all the $\ell = 1, \ldots, L-1$, we have the adjoints
\begin{equation*}
  \begin{aligned}
    \widebar{\mathbf{f}}^{(\ell)}
     & = - \mathbf{W} \left(\mathbf{f}^{(\ell)}, \mathbf{g}^{(\ell-1)}\right) \cdot \widebar{\mathbf{g}}^{(\ell)}, \\
    \widebar{\mathbf{g}}^{(\ell)}
     & = - \mathbf{X} \left(\mathbf{f}^{(\ell)}, \mathbf{g}^{(\ell)}\right) \cdot \widebar{\mathbf{f}}^{(\ell+1)}.
  \end{aligned}
\end{equation*}

Finally to have the gradient of the loss w.r.t.~$\mathbf{a}$, i.e.~the adjoint of $\mathbf{a}$, we have
\begin{equation*}
  \begin{aligned}
    \widebar{\mathbf{a}}
    = \sum_{\ell = 1}^{L} \left(\frac{\partial \mathbf{f}^{(\ell)}}{\partial \mathbf{a}}\right)^\top
    \cdot \frac{\partial \mathcal{L}}{\partial \mathbf{f}^{(\ell)}}
    = \sum_{\ell = 1}^{L} \left(\frac{\partial \mathbf{f}^{(\ell)}}{\partial \mathbf{a}}\right)^\top
    \cdot \widebar{\mathbf{f}}^{(\ell)}
    = \varepsilon \sum_{\ell = 1}^{L} \frac{\widebar{\mathbf{f}}^{(\ell)}}{\mathbf{a}}.
  \end{aligned}
\end{equation*}

\begin{algorithm}[H]
  \caption{Log-Stabilized Sinkhorn with Gradient}
  \begin{algorithmic}[1]\label{algo:log-sinkhorn-with-gradient}
    \Require $\mathbf{a} \in \Sigma_M$, $\mathbf{b} \in \Sigma_N$, $\mathbf{C} \in \mathbb{R}^{M\times N}$, $\varepsilon > 0$.
    \Initialize $\mathbf{f} = \mathbb{0}_M$, $\mathbf{g} = \mathbb{0}_N$.

    \State \# forward loop
    \While{$\ell = 1, \ldots, L$}
    \State $\mathbf{R} = \mathbf{C} - \mathbf{f}\cdot \mathbb{1}_N^\top - \mathbb{1}_M\cdot \mathbf{g}^\top$
    \State $\mathbf{f} =  \mathbf{f} + \varepsilon \log \mathbf{a} + \text{Min}^{row}_\varepsilon \mathbf{R}$
    \# update f
    \State $\mathbf{R} = \mathbf{C} - \mathbf{f}\cdot \mathbb{1}_N^\top - \mathbb{1}_M\cdot \mathbf{g}^\top$
    \State $\mathbf{g} = \mathbf{g} + \varepsilon \log \mathbf{b} + \text{Min}^{col}_\varepsilon \mathbf{R}$
    \# update g
    \EndWhile

    \State \# backward loop
    \State $\mathbf{P} = \exp \left(
      - \frac{\mathbf{C} - \mathbf{f}\cdot\mathbb{1}_N^\top - \mathbb{1}_M\cdot\mathbf{g}^\top}{\varepsilon}
      \right)$
    \State $\widebar{\mathbf{P}} = \mathbf{C} + \varepsilon \log \mathbf{P}$
    \State $\mathbf{W}$ from \cref{eqn:W-X-log-sinkhorn}
    \State $\widebar{\mathbf{g}}
      = \frac1\varepsilon \left(\widebar{\mathbf{P}} \odot \mathbf{P}\right)^\top \cdot \mathbb{1}_M$
    \State $\widebar{\mathbf{f}}
      = \frac1\varepsilon \left(\widebar{\mathbf{P}} \odot \mathbf{P}\right) \cdot \mathbb{1}_N
      - \mathbf{W} \cdot \widebar{\mathbf{g}}$
    \If{$\ell = L-1,\ldots,1$}
    \State $\mathbf{W}$ and $\mathbf{X}$ from \cref{eqn:W-X-log-sinkhorn}
    \State $\widebar{\mathbf{g}} = - \mathbf{X} \cdot \widebar{\mathbf{f}}$
    \State $\widebar{\mathbf{f}} = - \mathbf{W} \cdot \widebar{\mathbf{g}}$
    \State $\widebar{\mathbf{a}} = \widebar{\mathbf{a}} + \varepsilon\, \widebar{\mathbf{f}} \oslash \mathbf{a}$
    \EndIf

    \Ensure $\mathbf{P}$, $\widebar{\mathbf{a}}$.
  \end{algorithmic}
\end{algorithm}

\section{Wasserstein Barycenter}\label{sec:wasserstein-barycenter}
\sectionmark{Barycenter}

As is often in the case in practice, one wants to compute the ``mean'' or ``center'' given a sequence of points.
This is useful when we want to find a unknown center as an estimate to the group of points (usually data).
Suppose we have $\left(\mathbf{x}_s\right)_{s = 1}^S \in \mathcal{X}^S$ in a metric space $\left(\mathcal{X}, d\right)$,
where $d(\cdot)$ is a metric defined on $\mathcal{X}$,
then we have the following problem
\begin{equation*}
  \begin{aligned}
    \min_{\mathbf{x} \in \mathcal{X}} \sum_s^S \mathbf{w}_s d(x, x_s)^p,
  \end{aligned}
\end{equation*}

for a weight vector $\mathbf{w} \in \Sigma_S$ and some power $p$ but is often set to 2.
This is often called ``Fr\'echet'' or ``Karcher'' mean\footnote{
  Or Riemannian center of mass.
  See \citet{karcher2014} for a historical account on the naming.
}.

\subsection{Parallel Algorithm}\label{subsec:barycenter-parallel}
We can then state the Wasserstein Barycenter problem, similar to the idea of ``Fr\'echet mean''
as a weighted sum of the individual distances/losses.
Consider the loss function in \cref{eqn:loss-kantorovich},
a weight vector $\mathbf{w} \in \Sigma_S$,
and a sequence of discrete densities $\left\{\mathbf{a}_s\right\}_{s = 1}^S$,
where $a_s \in \Sigma_{M_s}$ and $M_s$ is the length of $a_s$,
a Wasserstein barycenter is computed by the following minimization problem
\begin{equation*}
  \begin{aligned}
    \argmin_{\mathbf{b} \in \Sigma_N} \sum_{s=1}^S \mathbf{w}_s \ell_{\mathbf{C}_s} \left(\mathbf{a}_s, \mathbf{b}\right),
  \end{aligned}
\end{equation*}

where the cost matrix $\mathbf{C}_s \in \mathbb{R}^{M_s \times N}$ specifies the cost for transporting
between the probability vectors $\mathbf{a}_s$ and $\mathbf{b}$.
Of course, this problem can be approximated by the entropic regularized version
by simply swapping the loss function by \cref{eqn:entropic-regularized-OT-loss}.
Thus, we have the entropic regularized Wasserstein barycenter problem
\begin{equation}\label{eqn:barycenter-problem-approximation}
  \begin{aligned}
    \argmin_{\mathbf{b} \in \Sigma_N} \sum_{s=1}^S \mathbf{w}_s \ell_{\mathbf{C}_s}^\varepsilon \left(\mathbf{a}_s, \mathbf{b}\right).
  \end{aligned}
\end{equation}

$\mathbf{b}$ is therefore the unknown computed barycenter as the ``mean'' of
the atoms $\mathbf{a}_1, \ldots, \mathbf{a}_S$.
By rewriting the problem in \cref{eqn:barycenter-problem-approximation}
as (weighted) Kullback-Leibler (KL) projection problem \citep{benamou2015,schmitz2018},
\begin{equation*}
  \begin{aligned}
    \min_{\left(\mathbf{P}_s\right)_s}
    \left\{
    \sum_s w_s \, \varepsilon \, \mathbf{KL} \left(\mathbf{P}_s \big| \mathbf{K}_s\right) :
    \mathbf{P}_s \cdot \mathbb{1}_N = \mathbf{a}_s,
    \mathbf{P}_1^\top \cdot \mathbb{1}_{M_1} = \ldots = \mathbf{P}_S^\top \cdot \mathbb{1}_{M_S} = \mathbf{b}
    \right\},
  \end{aligned}
\end{equation*}

we will arrive at the Sinkhorn-like updates:
\begin{equation}\label{eqn:barycenter-sinkhorn-like-vector-update}
  \begin{aligned}
    \mathbf{u}_s^{(\ell)}
     & = \frac{\mathbf{a}_s}{\mathbf{K}_s \mathbf{v}_s^{(\ell-1)}},                                      \\
    \mathbf{b}^{(\ell)}
     & = \boldsymbol{\Pi}_{s = 1}^S \left(\mathbf{K}_s^\top \mathbf{u}_s^{(\ell)}\right)^{\mathbf{w}_s}, \\
    \mathbf{v}_s^{(\ell)}
     & = \frac{\mathbf{b}^{(\ell)}}{\mathbf{K}_s^\top \mathbf{u}_s^{(\ell)}},
  \end{aligned}
\end{equation}

where $\mathbf{K}_s = \exp \left(-\frac{\mathbf{C}_s}{\varepsilon}\right)$,
and $\mathbf{u}_s$, $\mathbf{v}_s$ are scaling vectors for $s = 1, \ldots, S$.
Also, $\frac{\boldsymbol\cdot}{\boldsymbol\cdot}$ and $\left(\boldsymbol\cdot\right)^{(\boldsymbol\cdot)}$
are to be understood as element-wise operations.

In practical applications, however,
we usually have $M_1 = \ldots = M_S = M$
and $\mathbf{K} = \mathbf{K}_1 = \cdots = \mathbf{K}_S$.
We can thus rewrite \cref{eqn:barycenter-sinkhorn-like-vector-update}
similar to \cref{subsec:parallel-sinkhorn} as the following updating equations.
In some applications, for example Wasserstein Dictionary Learning~\citep{schmitz2018},
even $M = N$,
since both ``source'' densities and ``target'' densities are token distributions of same length.

\begin{update}[Updating Equations for Parallel Wasserstein Barycenter Algorithm]\label{update:parallel-barycenter}
  \begin{equation}\label{eqn:barycenter-sinkhorn-like-matrix-update}
    \begin{aligned}
      \mathbf{U}^{(\ell)}
       & = \frac{\mathbf{A}}{\mathbf{K} \mathbf{V}^{(\ell-1)}},                                     \\
      \mathbf{b}^{(\ell)}
       & = \boldsymbol{\Pi}_{row}
      \left(\mathbf{K}^\top \mathbf{U}^{(\ell)}\right)^{\mathbb{1}_N \cdot \mathbf{w}^\top},        \\
      \mathbf{V}^{(\ell)}
       & = \frac{\mathbf{b}^{(\ell)} \cdot \mathbb{1}_S^\top}{\mathbf{K}^\top \mathbf{U}^{(\ell)}},
    \end{aligned}
  \end{equation}
  where
  $\mathbf{A} = \left[\mathbf{a}, \ldots, \mathbf{a}\right] \in \mathbb{R}^{M \times S}$
  and the product operator $\boldsymbol\Pi_{row} \left(\cdot\right)$ refers to the row-wise product on the matrix.
\end{update}

The row-wise product will reduce the matrix into a length $N$ vector,
and the power operation in the second line refers to the element-wise power operation.
The matrix $\mathbb{1}_N \cdot \mathbf{w}^\top$ simply populates the length $S$ vector
$\mathbf{w}$ into an $N\times S$ matrix,
which can also be seen as applying each element of the $\mathbf{w}$ vector to the base matrix column-wise.

\begin{algorithm}[H]
  \caption{Parallel Wasserstein Barycenter Algorithm}
  \begin{algorithmic}[1]\label{algo:parallel-barycenter}
    \Require $\mathbf{A} \in \Sigma_{M \times S}$, $\mathbf{C} \in \mathbb{R}^{M \times N}$,
    $\mathbf{w} \in \Sigma_S$, $\varepsilon > 0$.
    \Initialize $\mathbf{U} = \mathbb{1}_{M \times S}$, $\mathbf{V} = \mathbb{1}_{N \times S}$,
    $\mathbf{b} = \mathbb{0}_N$.
    \State $\mathbf{K} = \exp(-\frac{\mathbf{C}}{\varepsilon})$
    \While{Not Convergence}
    \State $\mathbf{U} = \mathbf{A} \oslash (\mathbf{K} \mathbf{V})$
    \State $\mathbf{b} =
      \boldsymbol\Pi_{row}
      \left(\mathbf{K}^\top \mathbf{U}\right)^{\mathbb{1}_N \cdot \mathbf{w}^\top}
    $
    \State $\mathbf{V} = \left(\mathbf{b} \cdot \mathbb{1}_S^\top\right) \oslash (\mathbf{K}^\top \mathbf{U})$
    \EndWhile
    \Ensure $\mathbf{b}$
  \end{algorithmic}
\end{algorithm}

\begin{remark}[]
  In Line 4 of \cref{algo:parallel-barycenter},
  the notation $\boldsymbol\Pi_s$ refers to row-wise product reduction
  as in \cref{eqn:barycenter-sinkhorn-like-matrix-update}.
  In Line 5, the numerator $\mathbf{b}$ can be broadcast to be the same numerator of all columns of the denominator matrix
  as discussed in \cref{subsec:parallel-sinkhorn}.
\end{remark}

\begin{remark}[Convergence Condition]
  Similar to the Convergence Condition in \cref{remark:conv-cond-parallel},
  we can check the following condition for convergence. For some $\rho > 0$,
  \begin{equation*}
    \begin{aligned}
      \lVert
      \mathbf{U}^{(\ell)} \odot \left(\mathbf{K} \mathbf{V}^{(\ell)}\right) - \mathbf{A}
      \rVert_2 \le \rho.
    \end{aligned}
  \end{equation*}
\end{remark}

\subsection{Log-stabilized Algorithm}\label{subsec:log-barycenter}

Similar to the discussion of \cref{subsec:log-sinkhorn},
numerical instability is still an issue for \cref{algo:parallel-barycenter}
due to the computation of the exponential function and its potential numeric overflow.
Therefore, we would also need to stabilize the computation in the log-domain.
In order to do so, we will also need to state the algorithms (updating equations)
in terms of the dual variables $\mathbf{f}$ and $\mathbf{g}$.

Given \cref{update:update-fg-by-minrow-mincol}
we have that for $s = 1, \ldots, S$,

\begin{equation}\label{eqn:update-fs-gs-by-minrow-mincol}
  \begin{aligned}
    \mathbf{f}_s^{(\ell)}
     & = \mathbf{f}_s^{(\ell-1)} + \varepsilon \log \mathbf{a}
    + \text{Min}^{row}_\varepsilon \, \mathbf{R} \left(\mathbf{f}_s^{(\ell-1)}, \mathbf{g}_s^{(\ell-1)}\right), \\
    \mathbf{g}_s^{(\ell)}
     & = \mathbf{g}_s^{(\ell-1)} + \varepsilon \log \mathbf{b}^{(\ell)}
    + \text{Min}^{col}_\varepsilon \, \mathbf{R} \left(\mathbf{f}_s^{(\ell)}, \mathbf{g}_s^{(\ell-1)}\right).
  \end{aligned}
\end{equation}

Let us define

\begin{equation*}
  \begin{aligned}
    \mathbf{F}^{(\ell)}
     & \equiv \left[
    \mathbf{f}_1^{(\ell)}, \ldots, \mathbf{f}_S^{(\ell)}
    \right]
    =
    \left[
    \varepsilon \log \mathbf{u}_1^{(\ell)}, \ldots, \varepsilon \log \mathbf{u}_S^{(\ell)}
    \right]
    = \varepsilon \log \mathbf{U}^{(\ell)} \in \mathbb{R}^{M \times S}, \\
    \mathbf{G}^{(\ell)}
     & \equiv \left[
      \mathbf{g}^{(\ell)}_1, \ldots, \mathbf{g}^{(\ell)}_S
      \right]
    = \left[
      \varepsilon \log \mathbf{v}^{(\ell)}_1, \ldots, \varepsilon \log \mathbf{v}^{(\ell)}_S
      \right]
    = \varepsilon \log \mathbf{V}^{(\ell)} \in \mathbb{R}^{N \times S},
  \end{aligned}
\end{equation*}

and from \cref{eqn:update-UV-parallel} we will also need
$\log \left(\mathbf{K} \mathbf{V}^{(\ell-1)}\right)$ and
$\log \left(\mathbf{K}^\top \mathbf{U}^{(\ell-1)}\right)$.

From \cref{subsec:log-sinkhorn} we have

\begin{equation*}
  \begin{aligned}
    \log \left(\mathbf{K} \mathbf{v}_s^{(\ell-1)}\right)
     & = \log \left[
    \exp \left(-\frac{
    \mathbf{R} \left(\mathbf{f}_s^{(\ell-1)}, \mathbf{g}_s^{(\ell-1)}\right)
    }{\varepsilon}\right) \cdot \mathbb{1}_N
    \right]
    - \frac1\varepsilon \mathbf{f}_s^{(\ell-1)}, \\
    \log \left(\mathbf{K}^\top \mathbf{u}_s^{(\ell-1)}\right)
     & = \log \left[
    \exp \left(-\frac{
    \mathbf{R} \left(\mathbf{f}_s^{(\ell)}, \mathbf{g}_s^{(\ell-1)}\right)
    }{\varepsilon}\right)^\top \cdot \mathbb{1}_M
    \right]
    - \frac1\varepsilon \mathbf{g}_s^{(\ell-1)}.
  \end{aligned}
\end{equation*}

Thus,

\begin{equation}\label{eqn:updating-equations-FG-by-minrow-mincol}
  \begin{aligned}
    \mathbf{F}^{(\ell)}
     & = \varepsilon \log \mathbf{U}^{(\ell)}                                                                       \\
     & = \varepsilon \log \mathbf{A} - \varepsilon \log \left(\mathbf{K} \mathbf{V}^{(\ell-1)}\right)               \\
     & = \mathbf{F}^{(\ell-1)} + \varepsilon \log \mathbf{A}
    + \left[
    \ldots,
    - \varepsilon \log \left[
    \exp \left(-\frac{
    \mathbf{R} \left(\mathbf{f}_s^{(\ell-1)}, \mathbf{g}_s^{(\ell-1)}\right)
    }{\varepsilon}\right) \cdot \mathbb{1}_N
    \right],
    \ldots
    \right]                                                                                                         \\
     & = \mathbf{F}^{(\ell-1)} + \varepsilon \log \mathbf{A}
    + \left[
    \ldots,
    \text{Min}^{row}_\varepsilon \, \mathbf{R} \left(\mathbf{f}_s^{(\ell-1)}, \mathbf{g}_s^{(\ell-1)}\right),
    \ldots
    \right]                                                                                                         \\
     & = \mathbf{F}^{(\ell-1)} + \varepsilon \log \mathbf{A}
    + \text{Min}^{row}_\varepsilon \, \mathbf{R} \left(\mathbf{F}^{(\ell-1)}, \mathbf{G}^{(\ell-1)}\right),         \\
    \mathbf{G}^{(\ell)}
     & = \varepsilon \log \mathbf{V}^{(\ell)}                                                                       \\
     & = \varepsilon \log \mathbf{B}^{(\ell)} - \varepsilon \log \left(\mathbf{K}^\top \mathbf{U}^{(\ell-1)}\right) \\
     & = \mathbf{G}^{(\ell-1)} + \varepsilon \log \mathbf{b}^{(\ell)} \cdot \mathbb{1}_S^\top
    + \left[
    \ldots,
    - \varepsilon \log \left[
    \exp \left(-\frac{
    \mathbf{R} \left(\mathbf{f}_s^{(\ell)}, \mathbf{g}_s^{(\ell-1)}\right)
    }{\varepsilon}\right)^\top \cdot \mathbb{1}_M
    \right],
    \ldots
    \right]                                                                                                         \\
     & = \mathbf{G}^{(\ell-1)} + \varepsilon \log \mathbf{b}^{(\ell)} \cdot \mathbb{1}_S^\top
    + \left[
    \ldots,
    \text{Min}^{col}_{\varepsilon} \, \mathbf{R} \left(\mathbf{f}_s^{(\ell)}, \mathbf{g}_s^{(\ell-1)}\right),
    \ldots
    \right]                                                                                                         \\
     & = \mathbf{G}^{(\ell-1)} + \varepsilon \log \mathbf{b}^{(\ell)} \cdot \mathbb{1}_S^\top
    + \text{Min}^{col}_\varepsilon \, \mathbf{R} \left(\mathbf{F}^{(\ell)}, \mathbf{G}^{(\ell-1)}\right),
  \end{aligned}
\end{equation}

where we slightly abuse the notation that

\begin{dgroup*}
  \begin{dmath}\label{eqn:minrow-R-FG}
    \text{Min}^{row}_\varepsilon \, \mathbf{R} \left(\mathbf{F}^{(\ell-1)}, \mathbf{G}^{(\ell-1)}\right)
    = \left[
    \text{Min}^{row}_\varepsilon \, \mathbf{R} \left(\mathbf{f}_1^{(\ell-1)}, \mathbf{g}_1^{(\ell-1)}\right),
    \ldots,
    \text{Min}^{row}_\varepsilon \, \mathbf{R} \left(\mathbf{f}_s^{(\ell-1)}, \mathbf{g}_s^{(\ell-1)}\right),
    \ldots,
    \text{Min}^{row}_\varepsilon \, \mathbf{R} \left(\mathbf{f}_S^{(\ell-1)}, \mathbf{g}_S^{(\ell-1)}\right)
    \right],
  \end{dmath}
  \begin{dmath}\label{eqn:mincol-R-FG}
    \text{Min}^{col}_\varepsilon \, \mathbf{R} \left(\mathbf{F}^{(\ell)}, \mathbf{G}^{(\ell-1)}\right)
    = \left[
    \text{Min}^{col}_{\varepsilon} \, \mathbf{R} \left(\mathbf{f}_1^{(\ell)}, \mathbf{g}_1^{(\ell-1)}\right),
    \ldots,
    \text{Min}^{col}_{\varepsilon} \, \mathbf{R} \left(\mathbf{f}_s^{(\ell)}, \mathbf{g}_s^{(\ell-1)}\right),
    \ldots,
    \text{Min}^{col}_{\varepsilon} \, \mathbf{R} \left(\mathbf{f}_S^{(\ell)}, \mathbf{g}_S^{(\ell-1)}\right)
    \right].
  \end{dmath}
\end{dgroup*}

Since

\begin{equation*}
  \begin{aligned}
    \log \left(\mathbf{K}^\top \mathbf{U}^{(\ell)}\right)
     & = -\frac1\varepsilon \mathbf{G}^{(\ell-1)}
    + \left[
    \ldots,
    \log \left[
    \exp \left(-\frac{
    \mathbf{R} \left(\mathbf{f}_s^{(\ell)}, \mathbf{g}_s^{(\ell-1)}\right)
    }{\varepsilon}\right)^\top \cdot \mathbb{1}_M
    \right],
    \ldots
    \right]                                       \\
     & = -\frac1\varepsilon \mathbf{G}^{(\ell-1)}
    - \frac1\varepsilon \left[
    \ldots,
    \text{Min}^{col}_\varepsilon \, \mathbf{R} \left(\mathbf{f}_s^{(\ell)}, \mathbf{g}_s^{(\ell-1)}\right),
    \ldots
    \right]                                       \\
     & = -\frac1\varepsilon \mathbf{G}^{(\ell-1)}
    - \frac1\varepsilon
    \text{Min}^{col}_\varepsilon \, \mathbf{R} \left(\mathbf{F}^{(\ell)}, \mathbf{G}^{(\ell-1)}\right),
  \end{aligned}
\end{equation*}

then

\begin{equation}\label{eqn:updating-eqn-for-logb-by-mincol}
  \begin{aligned}
    \log \mathbf{b}^{(\ell)}
     & = \log\boldsymbol\Pi_{s} \left(\mathbf{K}^\top \mathbf{u}_s^{(\ell)}\right)^{\mathbf{w}_s}   \\
     & = \boldsymbol\Sigma_s \log \left(\mathbf{K}^\top \mathbf{u}^{(\ell)}_s\right)^{\mathbf{w}_s} \\
     & = \boldsymbol\Sigma_s \mathbf{w}_s \log \left(\mathbf{K}^\top \mathbf{u}^{(\ell)}_s\right)   \\
     & = \log \left(\mathbf{K}^\top \mathbf{U}^{(\ell)}\right) \cdot \mathbf{w}                     \\
     & = -\frac1\varepsilon \mathbf{G}^{(\ell-1)} \cdot \mathbf{w}
    - \frac1\varepsilon
    \text{Min}^{col}_\varepsilon \, \mathbf{R} \left(\mathbf{F}^{(\ell)}, \mathbf{G}^{(\ell-1)}\right)
    \cdot \mathbf{w}.
  \end{aligned}
\end{equation}

Putting everything together, we have the updating equations

\begin{update}[Updating Equations for the Log-Stabilized Wasserstein Barycenter Algorithm]
  \label{update:log-barycenter}
  \begin{equation*}
    \begin{aligned}
      \mathbf{F}^{(\ell)}
       & = \mathbf{F}^{(\ell-1)} + \varepsilon\log \mathbf{A}
      + \text{Min}_\varepsilon^{row} \, \mathbf{R} \left(\mathbf{F}^{(\ell-1)}, \mathbf{G}^{(\ell-1)}\right), \\
      \log \mathbf{b}^{(\ell)}
       & = -\frac1\varepsilon \mathbf{G}^{(\ell-1)} \cdot \mathbf{w}
      - \frac1\varepsilon
      \text{Min}_\varepsilon^{col} \,\mathbf{R} \left(\mathbf{F}^{(\ell)}, \mathbf{G}^{(\ell-1)}\right)
      \cdot \mathbf{w},                                                                                       \\
      \mathbf{G}^{(\ell)}
       & = \mathbf{G}^{(\ell-1)} + \varepsilon \log \mathbf{b}^{(\ell)} \cdot \mathbb{1}_S^\top
      + \text{Min}_\varepsilon^{col} \,\mathbf{R} \left(\mathbf{F}^{(\ell)}, \mathbf{G}^{(\ell-1)}\right).
    \end{aligned}
  \end{equation*}
\end{update}

\begin{algorithm}[H]
  \caption{Log-Stabilized Wasserstein Barycenter Algorithm}
  \begin{algorithmic}[1]\label{algo:log-barycenter}
    \Require $\mathbf{A} \in \Sigma_{M \times S}$, $\mathbf{C} \in \mathbb{R}^{M \times N}$, $\mathbf{w} \in \Sigma_S$,
    $\varepsilon > 0$.
    \Initialize $\mathbf{F} = \mathbb{0}_{M \times S}$, $\mathbf{G}_{N \times S} = \mathbb{1}_{N \times S}$,
    $\mathbf{logb} = \mathbb{0}_N$.
    \While{Not Convergence}
    \State $\mathbf{R}^{row}_{min}$ from \cref{eqn:minrow-R-FG}
    \State $\mathbf{F} = \mathbf{F} + \varepsilon \log \mathbf{A} + \mathbf{R}^{row}_{min}$ \# update F

    \State $\mathbf{R}^{col}_{min}$ from \cref{eqn:mincol-R-FG}
    \State $\mathbf{logb} = - \frac1\varepsilon \mathbf{G} \cdot \mathbf{w}
      - \frac1\varepsilon \mathbf{R}^{col}_{min} \cdot \mathbf{w}$ \# update logb

    \State $\mathbf{G} = \mathbf{G} + \varepsilon \mathbf{logb} \cdot \mathbb{1}_S^\top
      + \mathbf{R}^{col}_{min}$ \# update G

    \EndWhile
    \Ensure $\mathbf{b} = \exp \mathbf{logb}$
  \end{algorithmic}
\end{algorithm}

\section{Wasserstein Barycenter with Gradients}\label{sec:barycenter-gradient}
\sectionmark{Barycenter Gradients}

In \cref{sec:sinkhorn-algorithm}, apart from getting the optimal transport matrix $\mathbf{P}$,
we can also have the \textit{Sinkhorn loss} given by \cref{eqn:loss-kantorovich}.
Therefore, we could leverage the backward AD technique as in \cref{sec:sinkhorn-gradient}
to obtain the gradient of the final loss w.r.t.~the input $\mathbf{a}$.

However, in \cref{sec:wasserstein-barycenter} we have the ``predicted barycenter'' as output vector,
and the backward AD is most efficient when the output is scalar output.
Therefore, let us consider the quadratic loss between the barycenter and some other ``data vector''
as the final scalar loss, and we would be able to derive the gradients w.r.t.~$\mathbf{A}$ and $\mathbf{w}$
similar to \cref{sec:sinkhorn-gradient}.
Denote the ``data vector'' as $\mathbf{b}$ and the computed barycenter as $\mathbf{b}^{(L)}$ after $L$ iterations,
and note that $\mathbf{b}$ is external to the Wasserstein barycenter algorithm and is thus considered as constant.
Even though $\mathbf{b}$ is not needed for the computation of the barycenter algorithm as in
\cref{algo:parallel-barycenter,algo:log-barycenter}, the introduction of $\mathbf{b}$ and the gradients
will be useful for the Wasserstein Dictionary Learning algorithm in \cref{sec:wdl}.

\subsection{Gradients for Parallel Algorithm}\label{subsec:gradient-parallel-barycenter}

\begin{figure}[H]
  \centering
  \tikzset{every picture/.style={line width=0.75pt}} %

  \begin{tikzpicture}[x=0.75pt,y=0.75pt,yscale=-1,xscale=1]
    \draw    (215.11,62.06) -- (215.31,81) ;
    \draw [shift={(215.33,83)}, rotate = 269.39] [color={rgb, 255:red, 0; green, 0; blue, 0 }  ][line width=0.75]    (10.93,-3.29) .. controls (6.95,-1.4) and (3.31,-0.3) .. (0,0) .. controls (3.31,0.3) and (6.95,1.4) .. (10.93,3.29)   ;
    \draw    (409.78,60.72) -- (409.98,79.67) ;
    \draw [shift={(410,81.67)}, rotate = 269.39] [color={rgb, 255:red, 0; green, 0; blue, 0 }  ][line width=0.75]    (10.93,-3.29) .. controls (6.95,-1.4) and (3.31,-0.3) .. (0,0) .. controls (3.31,0.3) and (6.95,1.4) .. (10.93,3.29)   ;
    \draw    (215.78,108.72) -- (215.99,155) ;
    \draw [shift={(216,157)}, rotate = 269.74] [color={rgb, 255:red, 0; green, 0; blue, 0 }  ][line width=0.75]    (10.93,-3.29) .. controls (6.95,-1.4) and (3.31,-0.3) .. (0,0) .. controls (3.31,0.3) and (6.95,1.4) .. (10.93,3.29)   ;
    \draw    (180.44,154.06) -- (205.7,108.09) ;
    \draw [shift={(206.67,106.33)}, rotate = 118.79] [color={rgb, 255:red, 0; green, 0; blue, 0 }  ][line width=0.75]    (10.93,-3.29) .. controls (6.95,-1.4) and (3.31,-0.3) .. (0,0) .. controls (3.31,0.3) and (6.95,1.4) .. (10.93,3.29)   ;
    \draw    (375.78,150.72) -- (401.04,104.75) ;
    \draw [shift={(402,103)}, rotate = 118.79] [color={rgb, 255:red, 0; green, 0; blue, 0 }  ][line width=0.75]    (10.93,-3.29) .. controls (6.95,-1.4) and (3.31,-0.3) .. (0,0) .. controls (3.31,0.3) and (6.95,1.4) .. (10.93,3.29)   ;
    \draw    (227.11,102.72) -- (237.38,119.03) ;
    \draw [shift={(238.44,120.72)}, rotate = 237.8] [color={rgb, 255:red, 0; green, 0; blue, 0 }  ][line width=0.75]    (10.93,-3.29) .. controls (6.95,-1.4) and (3.31,-0.3) .. (0,0) .. controls (3.31,0.3) and (6.95,1.4) .. (10.93,3.29)   ;
    \draw    (421.78,102.72) -- (432.05,119.03) ;
    \draw [shift={(433.11,120.72)}, rotate = 237.8] [color={rgb, 255:red, 0; green, 0; blue, 0 }  ][line width=0.75]    (10.93,-3.29) .. controls (6.95,-1.4) and (3.31,-0.3) .. (0,0) .. controls (3.31,0.3) and (6.95,1.4) .. (10.93,3.29)   ;
    \draw    (457.78,140.06) -- (482.6,150.61) ;
    \draw [shift={(484.44,151.39)}, rotate = 203.03] [color={rgb, 255:red, 0; green, 0; blue, 0 }  ][line width=0.75]    (10.93,-3.29) .. controls (6.95,-1.4) and (3.31,-0.3) .. (0,0) .. controls (3.31,0.3) and (6.95,1.4) .. (10.93,3.29)   ;
    \draw    (497.44,106.06) -- (497.44,138.66) ;
    \draw [shift={(497.44,140.66)}, rotate = 270] [color={rgb, 255:red, 0; green, 0; blue, 0 }  ][line width=0.75]    (10.93,-3.29) .. controls (6.95,-1.4) and (3.31,-0.3) .. (0,0) .. controls (3.31,0.3) and (6.95,1.4) .. (10.93,3.29)   ;
    \draw    (239.78,140.99) -- (227.9,156.09) ;
    \draw [shift={(226.67,157.67)}, rotate = 308.17] [color={rgb, 255:red, 0; green, 0; blue, 0 }  ][line width=0.75]    (10.93,-3.29) .. controls (6.95,-1.4) and (3.31,-0.3) .. (0,0) .. controls (3.31,0.3) and (6.95,1.4) .. (10.93,3.29)   ;
    \draw    (247.44,69.99) -- (247.44,115.99) ;
    \draw [shift={(247.44,117.99)}, rotate = 270] [color={rgb, 255:red, 0; green, 0; blue, 0 }  ][line width=0.75]    (10.93,-3.29) .. controls (6.95,-1.4) and (3.31,-0.3) .. (0,0) .. controls (3.31,0.3) and (6.95,1.4) .. (10.93,3.29)   ;
    \draw    (441.44,69.99) -- (441.44,115.99) ;
    \draw [shift={(441.44,117.99)}, rotate = 270] [color={rgb, 255:red, 0; green, 0; blue, 0 }  ][line width=0.75]    (10.93,-3.29) .. controls (6.95,-1.4) and (3.31,-0.3) .. (0,0) .. controls (3.31,0.3) and (6.95,1.4) .. (10.93,3.29)   ;

    \draw (208,40) node [anchor=north west][inner sep=0.75pt]   [align=left] {$\displaystyle \mathbf{A}$};
    \draw (402.67,40) node [anchor=north west][inner sep=0.75pt]   [align=left] {$\displaystyle \mathbf{A}$};
    \draw (207.67,82.33) node [anchor=north west][inner sep=0.75pt]   [align=left] {$\displaystyle \mathbf{U}^{( 1)}$};
    \draw (402,80.67) node [anchor=north west][inner sep=0.75pt]   [align=left] {$\displaystyle \mathbf{U}^{( L)}$};
    \draw (165.33,160) node [anchor=north west][inner sep=0.75pt]   [align=left] {$\displaystyle \mathbf{V}^{( 0)}$};
    \draw (240,121.33) node [anchor=north west][inner sep=0.75pt]   [align=left] {$\displaystyle \mathbf{b}^{( 1)}$};
    \draw (435,121.33) node [anchor=north west][inner sep=0.75pt]   [align=left] {$\displaystyle \mathbf{b}^{( L)}$};
    \draw (492.67,84.67) node [anchor=north west][inner sep=0.75pt]   [align=left] {$\displaystyle \mathbf{b}$};
    \draw (365.67,159.67) node [anchor=north west][inner sep=0.75pt]   [align=left] {$\displaystyle \mathbf{V}^{( L-1)}$};
    \draw (490.33,143.67) node [anchor=north west][inner sep=0.75pt]   [align=left] {$\displaystyle \mathcal{L}$};
    \draw (208.33,159.67) node [anchor=north west][inner sep=0.75pt]   [align=left] {$\displaystyle \mathbf{V}^{( 1)}$};
    \draw (312,115.33) node [anchor=north west][inner sep=0.75pt]   [align=left] {$\displaystyle \dotsc \dotsc $};
    \draw (239.67,49.33) node [anchor=north west][inner sep=0.75pt]   [align=left] {$\displaystyle \mathbf{w}$};
    \draw (434,49.33) node [anchor=north west][inner sep=0.75pt]   [align=left] {$\displaystyle \mathbf{w}$};

  \end{tikzpicture}

  \caption{Computational graph for the parallel barycenter algorithm.}\label{fig:parallel-barycenter-comp-graph}
\end{figure}

As in \cref{sec:sinkhorn-gradient}, we will need the adjoints of $\mathbf{A}$ and $\mathbf{w}$,
i.e.~$\widebar{\mathbf{A}}$ and $\widebar{\mathbf{w}}$ as the desired output for the gradients.
To see this, we have
\begin{equation*}
  \begin{aligned}
    \widebar{\mathbf{A}} =
    \begin{bmatrix*}
      \nabla_{\mathbf{a}_1} \mathcal{L}, &
      \ldots,&
      \nabla_{\mathbf{a}_S} \mathcal{L}
    \end{bmatrix*},
  \end{aligned}
\end{equation*}

where each $\nabla_{\mathbf{a}_s} \mathcal{L}$ is
\begin{equation*}
  \begin{aligned}
    \nabla_{\mathbf{a}_s} \mathcal{L}
    = \left(\frac{\partial \mathcal{L}}{\partial \mathbf{a}_s}\right)^\top
    = \sum_{\ell = 1}^{L}
    \left(
    \frac{\partial \mathcal{L}}{\partial \mathbf{u}^{(\ell)}_s} \cdot
    \frac{\partial \mathbf{u}^{(\ell)}_s}{\partial \mathbf{a}_s}
    \right)^\top
    = \sum_{\ell = 1}^{L}
    \left(\frac{\partial \mathbf{u}^{(\ell)}_s}{\partial \mathbf{a}_s}\right)^\top \cdot
    \widebar{\mathbf{u}}^{(\ell)}_s
    = \sum_{\ell = 1}^{L} \frac{\widebar{\mathbf{u}}^{(\ell)}_s}{\mathbf{K} \mathbf{v}^{(\ell-1)}_s},
  \end{aligned}
\end{equation*}

and hence
\begin{equation*}
  \begin{aligned}
    \widebar{\mathbf{A}}
    = \begin{bmatrix*}
        \sum_{\ell = 1}^{L} \frac{\widebar{\mathbf{u}}^{(\ell)}_1}{\mathbf{K} \mathbf{v}^{(\ell-1)}_1},&
        \ldots,&
        \sum_{\ell = 1}^{L} \frac{\widebar{\mathbf{u}}^{(\ell)}_S}{\mathbf{K} \mathbf{v}^{(\ell-1)}_S}
      \end{bmatrix*}
    = \sum_{\ell = 1}^{L}
    \begin{bmatrix*}
      \frac{\widebar{\mathbf{u}}^{(\ell)}_1}{\mathbf{K} \mathbf{v}^{(\ell-1)}_1},&
      \ldots,&
      \frac{\widebar{\mathbf{u}}^{(\ell)}_S}{\mathbf{K} \mathbf{v}^{(\ell-1)}_S}
    \end{bmatrix*}
    = \sum_{\ell = 1}^{L} \frac{\widebar{\mathbf{U}}^{(\ell)}}{\mathbf{K} \mathbf{V}^{(\ell-1)}}.
  \end{aligned}
\end{equation*}

Then for the adjoint of $\mathbf{w}$, we have
\begin{equation*}
  \begin{aligned}
    \frac{\partial \mathcal{L}}{\partial \mathbf{w}_s}
    = \sum_{\ell = 1}^{L} \sum_j \frac{\partial \mathcal{L}}{\partial \mathbf{b}^{(\ell)}_j}\cdot
    \frac{\partial \mathbf{b}^{(\ell)}_j}{\partial \mathbf{w}_s}
    = \sum_{\ell = 1}^{L} \sum_j \widebar{\mathbf{b}}^{(\ell)}_j \cdot
    \frac{\partial \mathbf{b}^{(\ell)}_j}{\partial \mathbf{w}_s}.
  \end{aligned}
\end{equation*}

Now that from \cref{eqn:barycenter-sinkhorn-like-vector-update} we also have
\begin{equation*}
  \begin{aligned}
    \mathbf{b}^{(\ell)}_j
     & = \prod_t \left(\sum_{i'} \mathbf{K}_{i',j} \mathbf{u}^{(\ell)}_{i',j}\right)^{\mathbf{w}_t}, \\
    \log \mathbf{b}^{(\ell)}_j
     & = \sum_t \mathbf{w}_t \log \sum_{i'} \mathbf{K}_{i',j} \mathbf{u}^{(\ell)}_{i',t},
  \end{aligned}
\end{equation*}

and
\begin{equation*}
  \begin{aligned}
    \frac{\partial \mathbf{b}^{(\ell)}_j}{\partial \mathbf{w}_s}
     & = \frac{\exp\log \mathbf{b}^{(\ell)}_j}{\partial \mathbf{w}_s}                            \\
     & = \mathbf{b}^{(\ell)}_j \cdot
    \frac{\partial \log \mathbf{b}^{(\ell)}_j}{\partial \mathbf{w}_s}                            \\
     & = \mathbf{b}^{(\ell)}_j \cdot
    \frac{\partial }{\partial \mathbf{w}_s}
    \left(\sum_t \mathbf{w}_t \log \sum_{i'} \mathbf{K}_{i',j} \mathbf{u}^{(\ell)}_{i',t}\right) \\
     & = \mathbf{b}^{(\ell)}_j \cdot
    \log \left(\mathbf{K}^\top \mathbf{u}^{(\ell)}_s\right)_j.                                   \\
  \end{aligned}
\end{equation*}

Therefore,
\begin{equation*}
  \begin{aligned}
    \frac{\partial \mathcal{L}}{\partial \mathbf{w}_s}
     & = \sum_{\ell = 1}^{L} \sum_j \widebar{\mathbf{b}}^{(\ell)}_j \cdot
    \mathbf{b}^{(\ell)}_j \cdot \log \left(\mathbf{K}^\top \mathbf{U}^{(\ell)}_s\right)_{js} \\
     & = \sum_{\ell = 1}^{L} \left[
      \log \left(\mathbf{K}^\top \mathbf{U}^{(\ell)}\right)^\top \cdot
      \left(\widebar{\mathbf{b}}^{(\ell)} \odot \mathbf{b}^{(\ell)}\right)
      \right]_s,
  \end{aligned}
\end{equation*}

and the adjoint of $\mathbf{w}$ is
\begin{equation*}
  \begin{aligned}
    \widebar{\mathbf{w}}
    = \sum_{\ell = 1}^{L}
    \log \left(\mathbf{K}^\top \mathbf{U}^{(\ell)}\right)^\top \cdot
    \left(\widebar{\mathbf{b}}^{(\ell)} \odot \mathbf{b}^{(\ell)}\right).
  \end{aligned}
\end{equation*}

Here, all the remains to be derived are the adjoints for the intermediate variables,
$\widebar{\mathbf{U}}^{(\ell)}$, $\widebar{\mathbf{b}}^{(\ell)}$, and $\widebar{\mathbf{V}}^{(\ell)}$.
Once again, as in \cref{sec:sinkhorn-gradient}, we need to separate the cases for $\ell = L$ and $\ell = 1, \ldots, L-1$.

For $\ell=L$, we have the same $\widebar{\mathbf{b}}^{(L)} = 2 \left(\mathbf{b}^{(L)}-\mathbf{b}\right)$
from the quadratic loss\footnote{Also in Table 1 of \citet{schmitz2018} with other types of loss functions.},
and the adjoint of $\mathbf{U}^{(L)}$ is given by
\begin{equation*}
  \begin{aligned}
    \frac{\partial \mathcal{L}}{\partial \mathbf{u}^{(L)}_{is}}
    = \sum_j \frac{\partial \mathcal{L}}{\partial \mathbf{b}^{(L)}_j} \cdot
    \frac{\partial \mathbf{b}^{(L)}_j}{\partial \mathbf{u}^{(L)}_{is}}
    = \sum_j \widebar{\mathbf{b}}^{(L)} \cdot
    \frac{\partial \mathbf{b}^{(L)}_j}{\partial \mathbf{u}^{(L)}_{is}}.
  \end{aligned}
\end{equation*}

Then we have
\begin{equation*}
  \begin{aligned}
    \frac{\partial \mathbf{b}^{(L)}_j}{\partial \mathbf{u}^{(L)}_{is}}
     & = \mathbf{b}^{(L)}_j \cdot
    \frac{\partial }{\partial \mathbf{u}^{(L)}}_{is}
    \left(\sum_t \mathbf{w}_t \log \sum_{i'} \mathbf{K}_{i',j} \mathbf{u}^{(L)}_{i',t}\right) \\
     & =
    \mathbf{b}^{(L)}_j \cdot
    \mathbf{w}_s \cdot
    \frac1{\left(\mathbf{K}^\top \mathbf{u}^{(L)}_s\right)_j} \cdot \mathbf{K}_{ij}           \\
     & =
    \mathbf{w}_s \cdot
    \mathbf{K}_{ij} \cdot
    \left(\frac{\mathbf{b}^{(L)}}{\mathbf{K}^\top \mathbf{u}^{(L)}_s}\right)_j                \\
     & = \mathbf{w}_s \cdot
    \mathbf{K}_{ij} \cdot
    \mathbf{v}^{(L)}_{js},
  \end{aligned}
\end{equation*}

and then
\begin{equation*}
  \begin{aligned}
    \frac{\partial \mathcal{L}}{\partial \mathbf{u}^{(L)}_{is}}
     & = \sum_j \widebar{\mathbf{b}}^{(L)} \cdot
    \mathbf{w}_s \cdot
    \mathbf{K}_{ij} \cdot
    \mathbf{v}^{(L)}_{js}                        \\
     & = \sum_j
    \mathbf{K}_{ij} \cdot
    \left(\widebar{\mathbf{b}}^{(L)}_j \cdot \mathbf{w}_s\right) \cdot
    \mathbf{v}^{(L)}_{js}                        \\
     & = \sum_j
    \mathbf{K}_{ij}
    \left[\left(\widebar{\mathbf{b}}^{(L)} \cdot \mathbf{w}^\top\right) \odot \mathbf{V}^{(L)}\right]_{js}.
  \end{aligned}
\end{equation*}

Thus, the adjoint of $\mathbf{U}^{(L)}$ is
\begin{equation*}
  \begin{aligned}
    \widebar{\mathbf{U}}^{(L)} =
    \mathbf{K} \left[\left(\widebar{\mathbf{b}}^{(L)} \cdot \mathbf{w}^\top\right) \odot \mathbf{V}^{(L)}\right].
  \end{aligned}
\end{equation*}

Note that in the barycenter algorithm, $\mathbf{V}^{(L)}$ is not connected to the final loss,
so its adjoint is not needed for the gradient derivation.
We can then move on to the case of $\ell = 1, \ldots, L-1$, and we start from the adjoint of $\mathbf{V}^{(\ell)}$.
To see this, we have
\begin{equation*}
  \begin{aligned}
    \frac{\partial \mathcal{L}}{\partial \mathbf{v}^{(\ell)}_{js}}
    = \sum_i \frac{\partial \mathcal{L}}{\partial \mathbf{u}^{(\ell+1)}_{is}} \cdot
    \frac{\partial \mathbf{u}^{(\ell+1)}_{is}}{\partial \mathbf{v}^{(\ell)}_{js}}
    = \sum_i \widebar{\mathbf{U}}^{(\ell+1)}_{is}
    \frac{\partial \mathbf{u}^{(\ell+1)}_{is}}{\partial \mathbf{v}^{(\ell)}_{js}}.
  \end{aligned}
\end{equation*}

From \cref{eqn:barycenter-sinkhorn-like-vector-update} we have
\begin{equation*}
  \begin{aligned}
    \frac{\partial \mathbf{u}^{(\ell+1)}_{is}}{\partial \mathbf{v}^{(\ell)}_{js}}
     & = \frac{\partial }{\partial \mathbf{v}^{(\ell)}_{js}}
    \frac{\mathbf{a}_{is}}{\sum_{j'} \mathbf{K}_{ij'} \cdot \mathbf{v}^{(\ell)}_{j's}}                       \\
     & = - \frac{\mathbf{a}_{is}}{\left(\mathbf{K}\mathbf{v}^{(\ell)}_s\right)_i^2} \cdot
    \frac{\partial }{\partial \mathbf{v}^{(\ell)}_{js}}
    \sum_{j'} \mathbf{K}_{ij'} \cdot \mathbf{v}^{(\ell)}_{j's}                                               \\
     & = - \frac{\mathbf{a}_{is}}{\left(\mathbf{K}\mathbf{v}^{(\ell)}_s\right)_i^2} \cdot \mathbf{K}_{ij}    \\
     & = - \mathbf{K}_{ij} \cdot \left(\frac{\mathbf{u}^{(\ell+1)}}{\mathbf{K}\mathbf{v}^{(\ell)}}\right)_i,
  \end{aligned}
\end{equation*}

and
\begin{equation*}
  \begin{aligned}
    \frac{\partial \mathcal{L}}{\partial \mathbf{v}^{(\ell)}_{js}}
    = - \sum_i \widebar{\mathbf{U}}^{(\ell+1)}_{is} \cdot
    \mathbf{K}_{ij} \cdot \left(\frac{\mathbf{U}^{(\ell+1)}}{\mathbf{K}\mathbf{V}^{(\ell)}}\right)_{is}
    = - \sum_i \mathbf{K}_{ij} \cdot
    \left(\frac{\widebar{\mathbf{U}}^{(\ell+1)} \odot \mathbf{U}^{(\ell+1)}}{\mathbf{K}\mathbf{V}^{(\ell)}}\right)_{is}.
  \end{aligned}
\end{equation*}

Thus the adjoint for $\mathbf{V}^{(\ell)}$ is
\begin{equation}\label{eqn:adjoint-of-V-ell}
  \begin{aligned}
    \widebar{\mathbf{V}}^{(\ell)}
    = - \mathbf{K}^\top \cdot
    \frac{\widebar{\mathbf{U}}^{(\ell+1)} \odot \mathbf{U}^{(\ell+1)}}{\mathbf{K} \mathbf{V}^{(\ell)}}.
  \end{aligned}
\end{equation}

The adjoint of $\mathbf{b}^{(\ell)}$ is given by
\begin{equation*}
  \begin{aligned}
    \frac{\partial \mathcal{L}}{\partial \mathbf{b}^{(\ell)}_j}
    = \sum_s \frac{\partial \mathcal{L}}{\partial \mathbf{v}^{(\ell)}_{js}} \cdot
    \frac{\partial \mathbf{v}^{(\ell)}_{js}}{\partial \mathbf{b}^{(\ell)}_j}
    = \sum_s \widebar{\mathbf{v}}^{(\ell)}_{js} \cdot
    \frac1{\mathbf{K}^\top \mathbf{u}^{(\ell)}_s}_j
    = \sum_s \left(\frac{\widebar{\mathbf{v}}^{(\ell)}}{\mathbf{K}^\top \mathbf{u}^{(\ell)}_s}\right)_j,
  \end{aligned}
\end{equation*}

as
\begin{equation*}
  \begin{aligned}
    \frac{\partial \mathbf{v}^{(\ell)}_{js}}{\partial \mathbf{b}^{(\ell)}_j}
    = \left(\frac1{\mathbf{K}^\top \mathbf{u}^{(\ell)}_s}\right)_j,
  \end{aligned}
\end{equation*}

and the adjoint of $\mathbf{b}^{(\ell)}$ for $\ell = 1, \ldots, L-1$ is thus
\begin{equation}\label{eqn:adjoint-of-b-ell}
  \begin{aligned}
    \widebar{\mathbf{b}}^{(\ell)}
    = \sum_{s=1}^{S} \frac{\widebar{\mathbf{v}}^{(\ell)}_s}{\mathbf{K}^\top \mathbf{u}^{(\ell)}_s}
    = \sum^{row} \frac{\widebar{\mathbf{V}}^{(\ell)}}{\mathbf{K}^\top \mathbf{U}^{(\ell)}}.
  \end{aligned}
\end{equation}

All that remains is the adjoint for $\mathbf{U}^{(\ell)}$ for $\ell = 1, \ldots, L-1$.
Notice that $\mathbf{U}^{(\ell)}$'s are the only variables in the computational graph
\cref{fig:parallel-barycenter-comp-graph} to have two downstream variables,
and our backward pass needs to take that into consideration too.
Thus,
\begin{equation*}
  \begin{aligned}
    \frac{\partial \mathcal{L}}{\partial \mathbf{u}^{(\ell)}_{is}}
    = \sum_j \frac{\partial \mathcal{L}}{\partial \mathbf{b}^{(\ell)}_j} \cdot
    \frac{\partial \mathbf{b}^{(\ell)}_j}{\partial \mathbf{u}^{(\ell)}_{is}}
    + \sum_j \frac{\partial \mathcal{L}}{\partial \mathbf{v}^{(\ell)}_{js}} \cdot
    \frac{\partial \mathbf{v}^{(\ell)}_{js}}{\partial \mathbf{u}^{(\ell)}_{is}}.
  \end{aligned}
\end{equation*}

We already have that
\begin{equation*}
  \begin{aligned}
    \frac{\partial \mathbf{b}^{(\ell)}_j}{\partial \mathbf{u}^{(\ell)}_{is}}
    = \mathbf{w}_s \cdot \mathbf{K}_{ij} \cdot \mathbf{v}^{(\ell)}_{js},
  \end{aligned}
\end{equation*}

and also
\begin{equation*}
  \begin{aligned}
    \frac{\partial \mathbf{v}^{(\ell)}_{js}}{\partial \mathbf{u}^{(\ell)}_{is}}
    = \frac{\partial }{\partial \mathbf{u}^{(\ell)}_{is}}
    \frac{\mathbf{b}^{(\ell)}_j}{\sum_{i'} \mathbf{K}_{i'j} \mathbf{u}^{(\ell)}_{i's}}
    = - \frac{\mathbf{b}^{(\ell)}_j}{\left(\mathbf{K}^\top \mathbf{u}^{(\ell)}_s\right)_j}
    = - \mathbf{K}_{ij} \cdot \left(
    \frac{\mathbf{V}^{(\ell)}}{\mathbf{K}^\top \mathbf{U}^{(\ell)}}
    \right)_{js}.
  \end{aligned}
\end{equation*}

Therefore,
\begin{equation*}
  \begin{aligned}
    \frac{\partial \mathcal{L}}{\partial \mathbf{u}^{(\ell)}_{is}}
    = \sum_j \mathbf{K}_{ij} \cdot
    \left[\left(\widebar{\mathbf{b}}^{(\ell)} \cdot \mathbf{w}^\top \right) \odot \mathbf{V}^{(\ell)}\right]_{js}
    - \sum_j \mathbf{K}_{ij} \cdot
    \left(
    \frac{\widebar{\mathbf{V}}^{(\ell)} \odot \mathbf{V}^{(\ell)}}{\mathbf{K}^\top \mathbf{U}^{(\ell)}}
    \right)_{js},
  \end{aligned}
\end{equation*}

and the adjoint of $\mathbf{U}^{(\ell)}$ is
\begin{equation}\label{eqn:adjoint-of-U-ell}
  \begin{aligned}
    \widebar{\mathbf{U}}^{(\ell)}
     & = \mathbf{K} \left[\left(\widebar{\mathbf{b}}^{(\ell)} \cdot \mathbf{w}^\top \right) \odot \mathbf{V}^{(\ell)}\right]
    - \mathbf{K} \cdot \frac{\widebar{\mathbf{V}}^{(\ell)} \odot \mathbf{V}^{(\ell)}}{\mathbf{K}^\top \mathbf{U}^{(\ell)}}   \\
     & = \mathbf{K} \cdot
    \left[
      \left(
      \widebar{\mathbf{b}}^{(\ell)} \cdot \mathbf{w}^\top -
      \frac{\widebar{\mathbf{V}}^{(\ell)}}{\mathbf{K}^\top \mathbf{U}^{(\ell)}}
      \right) \odot \mathbf{V}^{(\ell)}
      \right].
  \end{aligned}
\end{equation}

\begin{algorithm}[H]
  \caption{Parallel Wasserstein Barycenter Algorithm with Gradients}
  \begin{algorithmic}[1]\label{algo:parallel-barycenter-with-gradient-Aw}
    \Require $\mathbf{A} \in \Sigma_{M \times S}$, $\widetilde{\mathbf{b}} \in \Sigma_N$,
    $\mathbf{C} \in \mathbb{R}^{M \times N}$, $\varepsilon > 0$.
    \Initialize $\mathbf{U} = \mathbb{1}_{M \times S}$, $\mathbf{V}_{N \times S} = \mathbb{1}_{N \times S}$,
    $\mathbf{b} = \mathbb{0}_N$.
    \State $\mathbf{K} = \exp(-\frac{\mathbf{C}}{\varepsilon})$

    \State \# forward loop
    \While{$\ell = 1, \ldots, L$}
    \State $\mathbf{U} = \mathbf{A} \oslash (\mathbf{K} \mathbf{V})$
    \State $\mathbf{b} =
      \boldsymbol\Pi_{row}
      \left(\mathbf{K}^\top \mathbf{U}\right)^{\mathbb{1}_N \cdot \mathbf{w}^\top}
    $
    \State $\mathbf{V} = \mathbf{b} \oslash (\mathbf{K}^\top \mathbf{U})$
    \EndWhile

    \State \# backward loop
    \State $\widebar{\mathbf{b}} = 2 (\mathbf{b} - \widetilde{\mathbf{b}})$
    \State $\widebar{\mathbf{U}}
      = \mathbf{K} \left[\left(\widebar{\mathbf{b}}\cdot \mathbf{w}^\top\right) \odot \mathbf{V}\right]$
    \State $\widebar{\mathbf{A}} = \frac{\widebar{\mathbf{U}}}{\mathbf{K} \mathbf{V}^{(L-1)}}$
    \State $\widebar{\mathbf{w}}
      = \log \left(\mathbf{K}^\top \mathbf{U}\right)^\top \cdot \left(\widebar{\mathbf{b}} \odot \mathbf{b}\right)$
    \If{$\ell = L-1, \ldots, 1$}
    \State $\widebar{\mathbf{V}}$ from \cref{eqn:adjoint-of-V-ell}
    \State $\widebar{\mathbf{b}}$ from \cref{eqn:adjoint-of-b-ell}
    \State $\widebar{\mathbf{U}}$ from \cref{eqn:adjoint-of-U-ell}
    \State $\widebar{\mathbf{A}} = \widebar{\mathbf{A}} + \frac{\widebar{\mathbf{U}}}{\mathbf{K} \mathbf{V}^{(\ell-1)}}$
    \State $\widebar{\mathbf{w}} = \widebar{\mathbf{w}} +
      \log \left(\mathbf{K}^\top \mathbf{U}^{(\ell)}\right)^\top \cdot \left(\widebar{\mathbf{b}}^{(\ell)} \odot \mathbf{b}^{(\ell)}\right)$
    \EndIf

    \Ensure $\mathbf{b}$, $\widebar{\mathbf{A}}$, $\widebar{\mathbf{w}}$
  \end{algorithmic}
\end{algorithm}

\subsection{Gradients for Log-Stabilized Algorithm}\label{subsec:gradient-log-barycenter}

\begin{figure}[H]
  \centering
  \tikzset{every picture/.style={line width=0.75pt}} %

  \begin{tikzpicture}[x=0.75pt,y=0.75pt,yscale=-1,xscale=1]
    \draw    (230.11,83.06) -- (230.31,102) ;
    \draw [shift={(230.33,104)}, rotate = 269.39] [color={rgb, 255:red, 0; green, 0; blue, 0 }  ][line width=0.75]    (10.93,-3.29) .. controls (6.95,-1.4) and (3.31,-0.3) .. (0,0) .. controls (3.31,0.3) and (6.95,1.4) .. (10.93,3.29)   ;
    \draw    (424.78,81.72) -- (424.98,100.67) ;
    \draw [shift={(425,102.67)}, rotate = 269.39] [color={rgb, 255:red, 0; green, 0; blue, 0 }  ][line width=0.75]    (10.93,-3.29) .. controls (6.95,-1.4) and (3.31,-0.3) .. (0,0) .. controls (3.31,0.3) and (6.95,1.4) .. (10.93,3.29)   ;
    \draw    (230.78,129.72) -- (230.99,176) ;
    \draw [shift={(231,178)}, rotate = 269.74] [color={rgb, 255:red, 0; green, 0; blue, 0 }  ][line width=0.75]    (10.93,-3.29) .. controls (6.95,-1.4) and (3.31,-0.3) .. (0,0) .. controls (3.31,0.3) and (6.95,1.4) .. (10.93,3.29)   ;
    \draw    (195.44,175.06) -- (220.7,129.09) ;
    \draw [shift={(221.67,127.33)}, rotate = 118.79] [color={rgb, 255:red, 0; green, 0; blue, 0 }  ][line width=0.75]    (10.93,-3.29) .. controls (6.95,-1.4) and (3.31,-0.3) .. (0,0) .. controls (3.31,0.3) and (6.95,1.4) .. (10.93,3.29)   ;
    \draw    (390.78,171.72) -- (416.04,125.75) ;
    \draw [shift={(417,124)}, rotate = 118.79] [color={rgb, 255:red, 0; green, 0; blue, 0 }  ][line width=0.75]    (10.93,-3.29) .. controls (6.95,-1.4) and (3.31,-0.3) .. (0,0) .. controls (3.31,0.3) and (6.95,1.4) .. (10.93,3.29)   ;
    \draw    (242.11,123.72) -- (252.38,140.03) ;
    \draw [shift={(253.44,141.72)}, rotate = 237.8] [color={rgb, 255:red, 0; green, 0; blue, 0 }  ][line width=0.75]    (10.93,-3.29) .. controls (6.95,-1.4) and (3.31,-0.3) .. (0,0) .. controls (3.31,0.3) and (6.95,1.4) .. (10.93,3.29)   ;
    \draw    (436.78,123.72) -- (447.05,140.03) ;
    \draw [shift={(448.11,141.72)}, rotate = 237.8] [color={rgb, 255:red, 0; green, 0; blue, 0 }  ][line width=0.75]    (10.93,-3.29) .. controls (6.95,-1.4) and (3.31,-0.3) .. (0,0) .. controls (3.31,0.3) and (6.95,1.4) .. (10.93,3.29)   ;
    \draw    (477.78,164.06) -- (502.6,174.61) ;
    \draw [shift={(504.44,175.39)}, rotate = 203.03] [color={rgb, 255:red, 0; green, 0; blue, 0 }  ][line width=0.75]    (10.93,-3.29) .. controls (6.95,-1.4) and (3.31,-0.3) .. (0,0) .. controls (3.31,0.3) and (6.95,1.4) .. (10.93,3.29)   ;
    \draw    (512.44,115.06) -- (512.44,165.66) ;
    \draw [shift={(512.44,167.66)}, rotate = 270] [color={rgb, 255:red, 0; green, 0; blue, 0 }  ][line width=0.75]    (10.93,-3.29) .. controls (6.95,-1.4) and (3.31,-0.3) .. (0,0) .. controls (3.31,0.3) and (6.95,1.4) .. (10.93,3.29)   ;
    \draw    (254.78,161.99) -- (242.9,177.09) ;
    \draw [shift={(241.67,178.67)}, rotate = 308.17] [color={rgb, 255:red, 0; green, 0; blue, 0 }  ][line width=0.75]    (10.93,-3.29) .. controls (6.95,-1.4) and (3.31,-0.3) .. (0,0) .. controls (3.31,0.3) and (6.95,1.4) .. (10.93,3.29)   ;
    \draw    (262.44,90.99) -- (262.44,136.99) ;
    \draw [shift={(262.44,138.99)}, rotate = 270] [color={rgb, 255:red, 0; green, 0; blue, 0 }  ][line width=0.75]    (10.93,-3.29) .. controls (6.95,-1.4) and (3.31,-0.3) .. (0,0) .. controls (3.31,0.3) and (6.95,1.4) .. (10.93,3.29)   ;
    \draw    (456.44,90.99) -- (456.44,136.99) ;
    \draw [shift={(456.44,138.99)}, rotate = 270] [color={rgb, 255:red, 0; green, 0; blue, 0 }  ][line width=0.75]    (10.93,-3.29) .. controls (6.95,-1.4) and (3.31,-0.3) .. (0,0) .. controls (3.31,0.3) and (6.95,1.4) .. (10.93,3.29)   ;
    \draw    (196.11,118.06) -- (217.33,118) ;
    \draw [shift={(219.33,118)}, rotate = 179.86] [color={rgb, 255:red, 0; green, 0; blue, 0 }  ][line width=0.75]    (10.93,-3.29) .. controls (6.95,-1.4) and (3.31,-0.3) .. (0,0) .. controls (3.31,0.3) and (6.95,1.4) .. (10.93,3.29)   ;
    \draw    (201.11,194.06) -- (220.33,194.01) ;
    \draw [shift={(222.33,194)}, rotate = 179.85] [color={rgb, 255:red, 0; green, 0; blue, 0 }  ][line width=0.75]    (10.93,-3.29) .. controls (6.95,-1.4) and (3.31,-0.3) .. (0,0) .. controls (3.31,0.3) and (6.95,1.4) .. (10.93,3.29)   ;
    \draw    (202.44,181.06) -- (241.02,154.47) ;
    \draw [shift={(242.67,153.33)}, rotate = 145.42] [color={rgb, 255:red, 0; green, 0; blue, 0 }  ][line width=0.75]    (10.93,-3.29) .. controls (6.95,-1.4) and (3.31,-0.3) .. (0,0) .. controls (3.31,0.3) and (6.95,1.4) .. (10.93,3.29)   ;
    \draw    (395.44,174.06) -- (434.9,153.27) ;
    \draw [shift={(436.67,152.33)}, rotate = 152.21] [color={rgb, 255:red, 0; green, 0; blue, 0 }  ][line width=0.75]    (10.93,-3.29) .. controls (6.95,-1.4) and (3.31,-0.3) .. (0,0) .. controls (3.31,0.3) and (6.95,1.4) .. (10.93,3.29)   ;
    \draw    (397.11,115.06) -- (412.33,115.01) ;
    \draw [shift={(414.33,115)}, rotate = 179.82] [color={rgb, 255:red, 0; green, 0; blue, 0 }  ][line width=0.75]    (10.93,-3.29) .. controls (6.95,-1.4) and (3.31,-0.3) .. (0,0) .. controls (3.31,0.3) and (6.95,1.4) .. (10.93,3.29)   ;

    \draw (223,61) node [anchor=north west][inner sep=0.75pt]   [align=left] {$\displaystyle \mathbf{A}$};
    \draw (417.67,61) node [anchor=north west][inner sep=0.75pt]   [align=left] {$\displaystyle \mathbf{A}$};
    \draw (222.67,103.33) node [anchor=north west][inner sep=0.75pt]   [align=left] {$\displaystyle \mathbf{F}^{( 1)}$};
    \draw (417,101.67) node [anchor=north west][inner sep=0.75pt]   [align=left] {$\displaystyle \mathbf{F}^{( L)}$};
    \draw (173.83,181) node [anchor=north west][inner sep=0.75pt]   [align=left] {$\displaystyle \mathbf{G}^{( 0)}$};
    \draw (243,141.33) node [anchor=north west][inner sep=0.75pt]   [align=left] {$\displaystyle \log\mathbf{b}^{( 1)}$};
    \draw (438,142.33) node [anchor=north west][inner sep=0.75pt]   [align=left] {$\displaystyle \log\mathbf{b}^{( L)}$};
    \draw (506.67,95.67) node [anchor=north west][inner sep=0.75pt]   [align=left] {$\displaystyle \mathbf{b}$};
    \draw (358.5,177.17) node [anchor=north west][inner sep=0.75pt]   [align=left] {$\displaystyle \mathbf{G}^{( L-1)}$};
    \draw (505.33,168.67) node [anchor=north west][inner sep=0.75pt]   [align=left] {$\displaystyle \mathcal{L}$};
    \draw (223.33,180.67) node [anchor=north west][inner sep=0.75pt]   [align=left] {$\displaystyle \mathbf{G}^{( 1)}$};
    \draw (327,136.33) node [anchor=north west][inner sep=0.75pt]   [align=left] {$\displaystyle \dotsc \dotsc $};
    \draw (254.67,70.33) node [anchor=north west][inner sep=0.75pt]   [align=left] {$\displaystyle \mathbf{w}$};
    \draw (449,70.33) node [anchor=north west][inner sep=0.75pt]   [align=left] {$\displaystyle \mathbf{w}$};
    \draw (173.83,103.33) node [anchor=north west][inner sep=0.75pt]   [align=left] {$\displaystyle \mathbf{F}^{( 0)}$};
    \draw (358.5,100.67) node [anchor=north west][inner sep=0.75pt]   [align=left] {$\displaystyle \mathbf{F}^{( L-1)}$};

  \end{tikzpicture}

  \caption{Computational graph for the log-stabilized barycenter algorithm.}\label{fig:log-barycenter-comp-graph}
\end{figure}

Similar to \cref{subsec:gradient-parallel-barycenter}, we will need to calculate the adjoints for
$\mathbf{A}$ and $\mathbf{w}$ as the gradients of the loss w.r.t.~the variables.
To see $\widebar{\mathbf{A}}$, we have
\begin{equation*}
  \begin{aligned}
    \frac{\partial \mathcal{L}}{\partial \mathbf{a}_{is}}
    = \sum_{\ell = 1}^{L}
    \frac{\partial \mathcal{L}}{\partial \mathbf{f}^{(\ell)}_{is}} \cdot
    \frac{\partial \mathbf{f}^{(\ell)}_{is}}{\partial \mathbf{a}_{is}}
    = \sum_{\ell = 1}^{L}
    \widebar{\mathbf{f}}^{(\ell)}_{is}\cdot
    \frac{\partial \mathbf{f}^{(\ell)}_{is}}{\partial \mathbf{a}_{is}},
  \end{aligned}
\end{equation*}

and from \cref{eqn:update-fs-gs-by-minrow-mincol},
\begin{equation*}
  \begin{aligned}
    \frac{\partial \mathbf{f}^{(\ell)}_{is}}{\partial \mathbf{a}_{is}}
    = \frac{\partial }{\partial \mathbf{a}_{is}}
    \left[
      \mathbf{f}^{(\ell-1)}_{is} + \varepsilon \log \mathbf{a}_{is} +
      \min_\varepsilon \mathbf{R} \left(\mathbf{f}^{(\ell-1)}_s, \mathbf{g}^{(\ell)}_s\right)_{i,\boldsymbol\cdot}
      \right]
    = \frac{\varepsilon}{\mathbf{a}_{is}}.
  \end{aligned}
\end{equation*}

Therefore, we have $\widebar{\mathbf{a}}_{is}
  = \sum_{\ell = 1}^{L} \widebar{\mathbf{f}}^{(\ell)}_{is} \cdot \frac{\varepsilon}{\mathbf{a}_{is}}$, and
\begin{equation*}
  \begin{aligned}
    \widebar{\mathbf{A}} = \varepsilon \sum_{\ell = 1}^{L} \frac{\widebar{\mathbf{F}}^{(\ell)}}{\mathbf{A}}.
  \end{aligned}
\end{equation*}

To get $\widebar{\mathbf{w}}$, we have
\begin{equation*}
  \begin{aligned}
    \frac{\partial \mathcal{L}}{\partial \mathbf{w}_s}
    = \sum_{\ell = 1}^{L} \sum_j
    \frac{\partial \mathcal{L}}{\partial \log \mathbf{b}^{(\ell)}_j} \cdot
    \frac{\log \mathbf{b}^{(\ell)}_j}{\partial \mathbf{w}_s}
    = \sum_{\ell = 1}^{L} \sum_j
    \widebar{\log\mathbf{b}}^{(\ell)}_j \cdot
    \frac{\log \mathbf{b}^{(\ell)}_j}{\partial \mathbf{w}_s}.
  \end{aligned}
\end{equation*}

Then we have
\begin{equation*}
  \begin{aligned}
    \frac{\log \mathbf{b}^{(\ell)}_j}{\partial \mathbf{w}_s}
     & = \frac{\partial }{\partial \mathbf{w}_s}
    \left[
      - \frac1\varepsilon \sum_t \mathbf{g}^{(\ell-1)}_{jt} \cdot \mathbf{w}_t
      - \frac1\varepsilon \sum_t
      \min_\varepsilon \mathbf{R} \left(\mathbf{f}^{(\ell)}_t, \mathbf{g}^{(\ell-1)}_t\right)_{\boldsymbol\cdot,j}
      \cdot \mathbf{w}_t
    \right]                                             \\
     & = - \frac1\varepsilon \mathbf{g}^{(\ell-1)}_{js}
    - \frac1\varepsilon
    \min_\varepsilon \mathbf{R} \left(\mathbf{f}^{(\ell)}_s, \mathbf{g}^{(\ell-1)}_s\right)_{\boldsymbol\cdot,j},
  \end{aligned}
\end{equation*}

and also
\begin{equation*}
  \begin{aligned}
    \sum_j
    \widebar{\log\mathbf{b}}^{(\ell)}_j \cdot
    \frac{\log \mathbf{b}^{(\ell)}_j}{\partial \mathbf{w}_s}
     & = - \frac1\varepsilon \sum_j \widebar{\log\mathbf{b}}^{(\ell)}_j \cdot \mathbf{g}^{(\ell-1)}_{js}
    - \frac1\varepsilon \sum_j \widebar{\log\mathbf{b}}^{(\ell)}_j \cdot
    \min_\varepsilon \mathbf{R} \left(\mathbf{f}^{(\ell)}_s, \mathbf{g}^{(\ell-1)}_s\right)_{\boldsymbol\cdot,j} \\
     & = - \frac1\varepsilon \sum_j \widebar{\log\mathbf{b}}^{(\ell)}_j \cdot
    \left[
      \mathbf{G}^{(\ell-1)} +
      \text{Min}_{\varepsilon}^{col} \left(\mathbf{F}^{(\ell)}, \mathbf{G}^{(\ell-1)}\right)
      \right]_{js}.
  \end{aligned}
\end{equation*}

Therefore, the adjoint for $\mathbf{w}$ is
\begin{equation*}
  \begin{aligned}
    \widebar{\mathbf{w}}
    = - \frac1\varepsilon \sum_{\ell = 1}^{L}
    \left[
      \mathbf{G}^{(\ell-1)} +
      \text{Min}_{\varepsilon}^{col} \left(\mathbf{F}^{(\ell)}, \mathbf{G}^{(\ell-1)}\right)
      \right]^\top\cdot \widebar{\log\mathbf{b}}^{(\ell)}.
  \end{aligned}
\end{equation*}

Next, we will need the adjoints for all the intermediate variables
$\mathbf{F}^{(\ell)}$, $\mathbf{G}^{(\ell)}$, and $\log\mathbf{b}^{(\ell)}$, for all $\ell = 1, \ldots, L$.
Similar to \cref{subsec:gradient-parallel-barycenter},
we need to separate the cases for $\ell = L$ and $\ell = 1, \ldots, L-1$.
For $\ell = L$, we have from the quadratic loss that
\begin{equation}\label{eqn:adjoint-of-logb-L}
  \begin{aligned}
    \widebar{\log \mathbf{b}}^{(L)} = 2 \left(\mathbf{b}^{(L)}-\mathbf{b}\right) \odot \mathbf{b}^{(L)},
  \end{aligned}
\end{equation}

and for $\ell = 1, \ldots, L-1$, we have
\begin{equation*}
  \begin{aligned}
    \frac{\partial \mathcal{L}}{\partial \log \mathbf{b}^{(\ell)}}
    = \sum_s \frac{\partial \mathcal{L}}{\partial \mathbf{g}^{(\ell)}_{js}} \cdot
    \frac{\partial \mathbf{g}^{(\ell)}_{js}}{\partial \log \mathbf{b}^{(\ell)}_j}
    = \sum_s \widebar{\mathbf{g}}^{(\ell)}_{js} \cdot
    \frac{\partial \mathbf{g}^{(\ell)}_{js}}{\partial \log \mathbf{b}^{(\ell)}_j}
  \end{aligned}
\end{equation*}

and from \cref{eqn:update-fs-gs-by-minrow-mincol},
\begin{equation*}
  \begin{aligned}
    \frac{\partial \mathbf{g}^{(\ell)}_{js}}{\partial \log \mathbf{b}^{(\ell)}_j}
     & = \frac{\partial }{\partial \log \mathbf{b}^{(\ell)}_j}
    \left[
      \mathbf{g}^{(\ell-1)}_{js} + \varepsilon\log \mathbf{b}^{(\ell)}_j +
      \min_\varepsilon \mathbf{R} \left(\mathbf{f}^{(\ell)}_s, \mathbf{g}^{(\ell-1)}_s\right)_{\boldsymbol\cdot,j}
      \right]
     & = \varepsilon.
  \end{aligned}
\end{equation*}

Hence,
\begin{equation*}
  \begin{aligned}
    \frac{\partial \mathcal{L}}{\partial \log \mathbf{b}^{(\ell)}}
    = \sum_s \widebar{\mathbf{g}}^{(\ell)}_{js} \cdot \varepsilon
    = \varepsilon \sum_s \widebar{\mathbf{g}}^{(\ell)}_{js},
  \end{aligned}
\end{equation*}
and
\begin{equation}\label{eqn:adjoint-of-logb-ell}
  \begin{aligned}
    \widebar{\log \mathbf{b}}^{(\ell)} = \varepsilon \sum_{row} \widebar{\mathbf{G}}^{(\ell)}.
  \end{aligned}
\end{equation}

For $\ell = L$, to find the adjoint of $\mathbf{F}^{(L)}$, we have
\begin{equation*}
  \begin{aligned}
    \frac{\partial \mathcal{L}}{\partial \mathbf{f}^{(L)}_{is}}
    = \sum_j \frac{\partial \mathcal{L}}{\partial \log \mathbf{b}^{(L)}_j} \cdot
    \frac{\partial \log \mathbf{b}^{(L)}_j}{\partial \mathbf{f}^{(L)}_{is}}
    = \sum_j \widebar{\log \mathbf{b}}^{(L)}_j \cdot
    \frac{\partial \log \mathbf{b}^{(L)}_j}{\partial \mathbf{f}^{(L)}_{is}},
  \end{aligned}
\end{equation*}

and
\begin{equation*}
  \begin{aligned}
    \frac{\partial \log \mathbf{b}^{(L)}_j}{\partial \mathbf{f}^{(L)}_{is}}
     & = \frac{\partial }{\partial \mathbf{f}^{(L)}_{is}}
    \left[
      - \frac1\varepsilon \sum_t \mathbf{g}^{(L-1)}_{jt}
      - \frac1\varepsilon \sum_t \min_\varepsilon
      \mathbf{R} \left(\mathbf{f}^{(L)}_t, \mathbf{g}^{(L-1)}_t\right)_{\boldsymbol\cdot,j} \cdot \mathbf{w}_t
    \right]                                                                            \\
     & = - \frac1\varepsilon \mathbf{w}_s \cdot
    \frac{
      \partial \min_\varepsilon
      \mathbf{R} \left(\mathbf{f}^{(L)}_s, \mathbf{g}^{(L-1)}_s\right)_{\boldsymbol\cdot,j}
    }{\partial \mathbf{f}^{(L)}_{is}}                                                  \\
     & = - \frac1\varepsilon \mathbf{w}_s \cdot
    \sum_k \frac{
      \partial \min_\varepsilon \mathbf{R} \left(\mathbf{f}^{(L)}_s, \mathbf{g}^{(L-1)}_s\right)_{\boldsymbol\cdot,j}
    }{\partial \mathbf{R} \left(\mathbf{f}^{(L)}_s, \mathbf{g}^{(L-1)}_s\right)_{k,j}} \cdot
    \frac{
      \partial \mathbf{R} \left(\mathbf{f}^{(L)}_s, \mathbf{g}^{(L-1)}_s\right)_{k,j}
    }{\partial \mathbf{f}^{(L)}_{is}}                                                  \\
     & = \frac1\varepsilon \mathbf{w}_s \cdot
    \frac{
      \partial \min_\varepsilon \mathbf{R} \left(\mathbf{f}^{(L)}_s, \mathbf{g}^{(L-1)}_s\right)_{\boldsymbol\cdot,j}
    }{\partial \mathbf{R} \left(\mathbf{f}^{(L)}_s, \mathbf{g}^{(L-1)}_s\right)_{i,j}} \\
     & = \frac1\varepsilon \mathbf{w}_s \cdot
    \frac{
      \exp \left(
      - \frac{
        \mathbf{R} \left(\mathbf{f}^{(L)}_s, \mathbf{g}^{(L-1)}_s\right)_{i,j} -
        \min \mathbf{R} \left(\mathbf{f}^{(L)}_s, \mathbf{g}^{(L-1)}_s\right)_{\boldsymbol\cdot,j}
      }\varepsilon
      \right)
    }{
      \sum_{i'}
      \exp \left(
      - \frac{
        \mathbf{R} \left(\mathbf{f}^{(L)}_s, \mathbf{g}^{(L-1)}_s\right)_{i',j} -
        \min \mathbf{R} \left(\mathbf{f}^{(L)}_s, \mathbf{g}^{(L-1)}_s\right)_{\boldsymbol\cdot,j}
      }\varepsilon
      \right)
    }                                                                                  \\
     & = \frac1\varepsilon \mathbf{w}_s \cdot
    \mathbf{W} \left(\mathbf{f}^{(L)}_s, \mathbf{g}^{(L-1)}_s\right)_{i,j},
  \end{aligned}
\end{equation*}

where
\begin{equation*}
  \begin{aligned}
    \mathbf{W} \left(\mathbf{f}^{(L)}_s, \mathbf{g}^{(L-1)}_s\right)
    = \begin{bmatrix*}
        \ldots, &
        \nabla_{
          \mathbf{R} \left(\mathbf{f}^{(L)}_s, \mathbf{g}^{(L-1)}_s\right)_{\boldsymbol\cdot,j}
        } \min_\varepsilon \mathbf{R} \left(\mathbf{f}^{(L)}_s, \mathbf{g}^{(L-1)}_s\right)_{\boldsymbol\cdot,j},&
        \ldots
      \end{bmatrix*},
  \end{aligned}
\end{equation*}

which is identical to \cref{subsec:gradient-log-sinkhorn}.
Therefore,
\begin{equation*}
  \begin{aligned}
    \frac{\partial \mathcal{L}}{\partial \mathbf{f}^{(L)}_{is}}
     & = \sum_j \widebar{\log \mathbf{b}}^{(L)}_j \cdot
    \frac1\varepsilon \mathbf{w}_s \cdot
    \mathbf{W} \left(\mathbf{f}^{(L)}_s, \mathbf{g}^{(L-1)}_s\right)_{i,j} \\
     & = \frac1\varepsilon \mathbf{w}_s \cdot
    \left[
      \mathbf{W} \left(\mathbf{f}^{(L)}_s, \mathbf{g}^{(L-1)}_s\right) \cdot
      \widebar{\log \mathbf{b}}^{(L)}
      \right]_i,
  \end{aligned}
\end{equation*}

and the adjoint for each $\mathbf{f}^{(L)}_s$ is
\begin{equation*}
  \begin{aligned}
    \widebar{\mathbf{f}}^{(L)}_s
    = \frac1\varepsilon \mathbf{w}_s \cdot
    \mathbf{W} \left(\mathbf{f}^{(L)}_s, \mathbf{g}^{(L-1)}_s\right) \cdot
    \widebar{\log \mathbf{b}}^{(L)}.
  \end{aligned}
\end{equation*}

Finally, we can have the adjoint for $\mathbf{F}^{(L)}$
\begin{equation}\label{eqn:adjoint-of-F-L}
  \begin{aligned}
    \widebar{\mathbf{F}}^{(L)} =
    \begin{bmatrix*}
      \bar{\mathbf{f}}^{(L)}_1, &
      \ldots, &
      \bar{\mathbf{f}}^{(L)}_s, &
      \ldots, &
      \bar{\mathbf{f}}^{(L)}_{S}
    \end{bmatrix*}.
  \end{aligned}
\end{equation}

Next, we will also need the adjoint for $\mathbf{F}^{(\ell)}$, $\ell = 1, \ldots, L-1$.
Note that, from \cref{fig:log-barycenter-comp-graph}, we can see that $\mathbf{F}^{(\ell)}$ has three children nodes:
$\mathbf{F}^{(\ell+1)}$, $\log \mathbf{b}^{(\ell)}$, and $\mathbf{G}^{(\ell)}$. Then,
\begin{equation*}
  \begin{aligned}
    \frac{\partial \mathcal{L}}{\partial \mathbf{f}^{(\ell)}_{is}}
    = \sum_j \frac{\partial \mathcal{L}}{\partial \log \mathbf{b}^{(\ell)}_j} \cdot
    \frac{\partial \log \mathbf{b}^{(\ell)}_j}{\partial \mathbf{f}^{(\ell)}_{is}}
    + \frac{\partial \mathcal{L}}{\partial \mathbf{f}^{(\ell+1)}} \cdot
    \frac{\partial \mathbf{f}^{(\ell+1)}_{is}}{\partial \mathbf{f}^{(\ell)}_{is}}
    + \sum_j \frac{\partial \mathcal{L}}{\partial \mathbf{g}^{(\ell)}_{js}} \cdot
    \frac{\partial \mathbf{g}^{(\ell)}_{js}}{\partial \mathbf{f}^{(\ell)}_{is}}.
  \end{aligned}
\end{equation*}

Since from \cref{subsec:gradient-log-sinkhorn} we have
\begin{equation*}
  \begin{aligned}
    \frac{\partial \mathbf{f}^{(\ell+1)}_{is}}{\partial \mathbf{f}^{(\ell)}_{is}} = 0,
  \end{aligned}
\end{equation*}

and
\begin{equation*}
  \begin{aligned}
    \frac{\partial \mathbf{g}^{(\ell)}_{js}}{\partial \mathbf{f}^{(\ell)}_{is}}
     & = \frac{\partial }{\partial \mathbf{f}^{(\ell)}_{is}}
    \left[
      \mathbf{g}^{(\ell-1)}_{js} +
      \varepsilon \log \mathbf{b}^{(\ell)}_j +
      \min_\varepsilon \mathbf{R}\left(\mathbf{f}^{(\ell)}_s, \mathbf{g}^{(\ell-1)}_s\right)_{\boldsymbol\cdot,j}
    \right]                                                                                                     \\
     & = \frac{\partial }{\partial \mathbf{f}^{(\ell)}_{is}}
    \min_\varepsilon \mathbf{R}\left(\mathbf{f}^{(\ell)}_s, \mathbf{g}^{(\ell-1)}_s\right)_{\boldsymbol\cdot,j} \\
     & = \sum_k \frac{
      \min_\varepsilon \mathbf{R}\left(\mathbf{f}^{(\ell)}_s, \mathbf{g}^{(\ell-1)}_s\right)_{\boldsymbol\cdot,j}
    }{
      \mathbf{R}\left(\mathbf{f}^{(\ell)}_s, \mathbf{g}^{(\ell-1)}_s\right)_{k,j}
    } \cdot \frac{
      \mathbf{R}\left(\mathbf{f}^{(\ell)}_s, \mathbf{g}^{(\ell-1)}_s\right)_{k,j}
    }{\partial \mathbf{f}^{(\ell)}_{is}}                                                                        \\
     & = - \frac{
      \min_\varepsilon \mathbf{R}\left(\mathbf{f}^{(\ell)}_s, \mathbf{g}^{(\ell-1)}_s\right)_{\boldsymbol\cdot,j}
    }{
      \mathbf{R}\left(\mathbf{f}^{(\ell)}_s, \mathbf{g}^{(\ell-1)}_s\right)_{i,j}
    }                                                                                                           \\
     & = - \mathbf{W} \left(\mathbf{f}^{(\ell)}_s, \mathbf{g}^{(\ell-1)}_s\right)_{i,j},
  \end{aligned}
\end{equation*}

where
\begin{equation*}
  \begin{aligned}
    \mathbf{W} \left(\mathbf{f}^{(\ell)}_s, \mathbf{g}^{(\ell-1)}_s\right) =
    \begin{bmatrix*}
      \ldots, &
      \nabla_{
        \mathbf{R}\left(\mathbf{f}^{(\ell)}_s, \mathbf{g}^{(\ell-1)}_s\right)_{\boldsymbol\cdot,j}
      } \min_\varepsilon \mathbf{R}\left(\mathbf{f}^{(\ell)}_s, \mathbf{g}^{(\ell-1)}_s\right)_{\boldsymbol\cdot,j},&
      \ldots
    \end{bmatrix*}.
  \end{aligned}
\end{equation*}

Similar to the case of $\widebar{\mathbf{F}}^{(L)}$, we have
\begin{equation*}
  \begin{aligned}
    \frac{\partial \log \mathbf{b}^{(\ell)}_j}{\partial \mathbf{f}^{(\ell)}_{is}}
    = \frac1\varepsilon \mathbf{w}_s \cdot
    \mathbf{W} \left(\mathbf{f}^{(\ell)}_s, \mathbf{g}^{(\ell-1)}_s\right)_{i,j},
  \end{aligned}
\end{equation*}

thus, we have
\begin{equation*}
  \begin{aligned}
    \frac{\partial \mathcal{L}}{\partial \mathbf{f}^{(\ell)}_{is}}
    = \frac1\varepsilon \sum_j \mathbf{w}_s \cdot
    \mathbf{W} \left(\mathbf{f}^{(\ell)}_s, \mathbf{g}^{(\ell-1)}_s\right)_{i,j}\cdot
    \widebar{\log \mathbf{b}}^{(\ell)}_j -
    \sum_j
    \mathbf{W} \left(\mathbf{f}^{(\ell)}_s, \mathbf{g}^{(\ell-1)}_s\right)_{i,j}\cdot
    \widebar{\mathbf{g}}^{(\ell)}_{js}.
  \end{aligned}
\end{equation*}

Finally, we have for $\ell = 1, \ldots, L-1$ the adjoint of $\mathbf{f}^{(\ell)}_s$
\begin{equation*}
  \begin{aligned}
    \widebar{\mathbf{f}}^{(\ell)}_s
    = \frac1\varepsilon \mathbf{w}_s
    \cdot \mathbf{W} \left(\mathbf{f}^{(\ell)}_s, \mathbf{g}^{(\ell-1)}_s\right) \cdot
    \widebar{\log \mathbf{b}}^{(\ell)}
    - \mathbf{W} \left(\mathbf{f}^{(\ell)}_s, \mathbf{g}^{(\ell-1)}_s\right) \cdot
    \widebar{\mathbf{g}}^{(\ell)}_s,
  \end{aligned}
\end{equation*}

and
\begin{equation}\label{eqn:adjoint-of-F-ell}
  \begin{aligned}
    \widebar{\mathbf{F}}^{(\ell)} =
    \begin{bmatrix*}
      \bar{\mathbf{f}}^{(\ell)}_1,&
      \ldots,&
      \bar{\mathbf{f}}^{(\ell)}_s,&
      \ldots,&
      \bar{\mathbf{f}}^{(\ell)}_S
    \end{bmatrix*}.
  \end{aligned}
\end{equation}

Now what remain to be derived are the adjoints for the $\mathbf{G}^{(\ell)}$'s, for $\ell = 1, \ldots, L-1$,
as $\mathbf{G}^{(L)}$ is not needed for the computation of the gradient from \cref{fig:log-barycenter-comp-graph}.
Further, $\mathbf{G}^{(L-1)}$ has only two children nodes $\mathbf{F}^{(L)}$ and $\log \mathbf{b}^{(L)}$,
whereas $\mathbf{G}^{(\ell)}$ has three: $\mathbf{F}^{(\ell+1)}$, $\log\mathbf{b}^{(\ell+1)}$, and $\mathbf{G}^{(\ell+1)}$
for $\ell = 1, \ldots, L-2$.
Thus, we have
\begin{equation*}
  \begin{aligned}
    \frac{\partial \mathcal{L}}{\partial \mathbf{g}^{(L-1)}_{js}}
     & = \sum_i \frac{\partial \mathcal{L}}{\partial \mathbf{f}^{(L)}_{is}}\cdot
    \frac{\partial \mathbf{f}^{(L)}_{is}}{\partial \mathbf{g}^{(L-1)}_{js}}
    + \frac{\partial \mathcal{L}}{\partial \log \mathbf{b}^{(L)}_j}\cdot
    \frac{\partial \log \mathbf{b}^{(L)}_j}{\partial \mathbf{g}^{(L-1)}_{js}},        \\
    \frac{\partial \mathcal{L}}{\partial \mathbf{g}^{(\ell)}_{js}}
     & = \sum_i \frac{\partial \mathcal{L}}{\partial \mathbf{f}^{(\ell+1)}_{is}}\cdot
    \frac{\partial \mathbf{f}^{(\ell+1)}_{is}}{\partial \mathbf{g}^{(\ell)}_{js}}
    + \frac{\partial \mathcal{L}}{\partial \log \mathbf{b}^{(\ell+1)}_j}\cdot
    \frac{\partial \log \mathbf{b}^{(\ell+1)}_j}{\partial \mathbf{g}^{(\ell)}_{js}}
    + \frac{\partial \mathcal{L}}{\partial \mathbf{g}^{(\ell+1)}_{js}}\cdot
    \frac{\partial \mathbf{g}^{(\ell+1)}_{js}}{\partial \mathbf{g}^{(\ell)}_{js}},
  \end{aligned}
\end{equation*}

where $\ell = 1, \ldots, L-2$.
Then we can show that from \cref{subsec:gradient-log-sinkhorn},
\begin{equation*}
  \begin{aligned}
    \frac{\partial \mathbf{g}^{(\ell+1)}_{js}}{\partial \mathbf{g}^{(\ell)}_{js}}
     & = \frac{\partial }{\partial \mathbf{g}^{(\ell)}_{js}}
    \left[
      \mathbf{g}^{(\ell)}_{js} +
      \varepsilon\log \mathbf{b}^{(\ell+1)}_{j} +
      \min_\varepsilon \mathbf{R} \left(\mathbf{f}^{(\ell+1)}_s, \mathbf{g}^{(\ell)}_s\right)_{\boldsymbol\cdot,j}
    \right]                                                  \\
     & = 1 + \frac{
      \min_\varepsilon \mathbf{R} \left(\mathbf{f}^{(\ell+1)}_s, \mathbf{g}^{(\ell)}_s\right)_{\boldsymbol\cdot,j}
    }{\partial \mathbf{g}^{(\ell)}_{js}}                     \\
     & = 1 + \sum_k \frac{
      \min_\varepsilon \mathbf{R} \left(\mathbf{f}^{(\ell+1)}_s, \mathbf{g}^{(\ell)}_s\right)_{\boldsymbol\cdot,j}
    }{
      \mathbf{R} \left(\mathbf{f}^{(\ell+1)}_s, \mathbf{g}^{(\ell)}_s\right)_{k,j}
    }\cdot \frac{
      \mathbf{R} \left(\mathbf{f}^{(\ell+1)}_s, \mathbf{g}^{(\ell)}_s\right)_{k,j}
    }{\partial \mathbf{g}^{(\ell)}_{js}}                     \\
     & = 1 - \sum_k \frac{
      \min_\varepsilon \mathbf{R} \left(\mathbf{f}^{(\ell+1)}_s, \mathbf{g}^{(\ell)}_s\right)_{\boldsymbol\cdot,j}
    }{
      \mathbf{R} \left(\mathbf{f}^{(\ell+1)}_s, \mathbf{g}^{(\ell)}_s\right)_{k,j}
    }                                                        \\
     & = 0,
  \end{aligned}
\end{equation*}

and thus for all $\ell = 1, \ldots, L-1$ we have,
\begin{equation*}
  \begin{aligned}
    \frac{\partial \mathcal{L}}{\partial \mathbf{g}^{(\ell)}_{js}}
     & = \sum_i \frac{\partial \mathcal{L}}{\partial \mathbf{f}^{(\ell+1)}_{is}}\cdot
    \frac{\partial \mathbf{f}^{(\ell+1)}_{is}}{\partial \mathbf{g}^{(\ell)}_{js}}
    + \frac{\partial \mathcal{L}}{\partial \log \mathbf{b}^{(\ell+1)}_j}\cdot
    \frac{\partial \log \mathbf{b}^{(\ell+1)}_j}{\partial \mathbf{g}^{(\ell)}_{js}}.
  \end{aligned}
\end{equation*}

Now that we also have
\begin{equation*}
  \begin{aligned}
    \frac{\partial \log \mathbf{b}^{(\ell+1)}_j}{\partial \mathbf{g}^{(\ell)}_{js}}
     & = \frac{\partial }{\partial \mathbf{g}^{(\ell)}_{js}}
    \left[
      - \frac1\varepsilon \sum_t \mathbf{g}^{(\ell)}_{jt} \cdot \mathbf{w}_t
      - \frac1\varepsilon \sum_t
      \min_\varepsilon \mathbf{R} \left(\mathbf{f}^{(\ell+1)}_t, \mathbf{g}^{(\ell)}_t\right)_{\boldsymbol\cdot,j}
      \cdot \mathbf{w}_t
    \right]                                                  \\
     & = - \frac1\varepsilon \mathbf{w}_s
    - \frac1\varepsilon \mathbf{w}_s \cdot
    \frac{
      \partial \min_\varepsilon \mathbf{R} \left(\mathbf{f}^{(\ell+1)}_t, \mathbf{g}^{(\ell)}_t\right)_{\boldsymbol\cdot,j}
    }{\partial \mathbf{g}^{(\ell)}_{js}}                     \\
     & = - \frac1\varepsilon \mathbf{w}_s \left(
    1 + \frac{
      \partial \min_\varepsilon \mathbf{R} \left(\mathbf{f}^{(\ell+1)}_t, \mathbf{g}^{(\ell)}_t\right)_{\boldsymbol\cdot,j}
    }{\partial \mathbf{g}^{(\ell)}_{js}}
    \right)
     & = 0,
  \end{aligned}
\end{equation*}

and
\begin{equation*}
  \begin{aligned}
    \frac{\partial \mathbf{f}^{(\ell+1)}_{is}}{\partial \mathbf{g}^{(\ell)}_{js}}
     & = \frac{\partial }{\partial \mathbf{g}^{(\ell)}_{js}}
    \left[
      \mathbf{f}^{(\ell)}_{is} +
      \varepsilon \log \mathbf{a}_{is} +
      \min_\varepsilon \mathbf{R} \left(\mathbf{f}^{(\ell)}_s, \mathbf{g}^{(\ell)}_s\right)_{i,\boldsymbol\cdot}
    \right]                                                                            \\
     & = \frac{
      \partial \min_\varepsilon \mathbf{R} \left(\mathbf{f}^{(\ell)}_s, \mathbf{g}^{(\ell)}_s\right)_{i,\boldsymbol\cdot}
    }{\partial \mathbf{g}^{(\ell)}_{js}}                                               \\
     & = \sum_k \frac{
      \partial \min_\varepsilon \mathbf{R} \left(\mathbf{f}^{(\ell)}_s, \mathbf{g}^{(\ell)}_s\right)_{i,\boldsymbol\cdot}
    }{
      \mathbf{R} \left(\mathbf{f}^{(\ell)}_s, \mathbf{g}^{(\ell)}_s\right)_{i,k}
    } \cdot \frac{
      \mathbf{R} \left(\mathbf{f}^{(\ell)}_s, \mathbf{g}^{(\ell)}_s\right)_{i,k}
    }{\mathbf{g}^{(\ell)}_{js}}                                                        \\
     & = - \frac{
      \partial \min_\varepsilon \mathbf{R} \left(\mathbf{f}^{(\ell)}_s, \mathbf{g}^{(\ell)}_s\right)_{i,\boldsymbol\cdot}
    }{
      \mathbf{R} \left(\mathbf{f}^{(\ell)}_s, \mathbf{g}^{(\ell)}_s\right)_{i,j}
    }                                                                                  \\
     & = - \mathbf{X} \left(\mathbf{f}^{(\ell)}_s, \mathbf{g}^{(\ell)}_s\right)_{j,i},
  \end{aligned}
\end{equation*}
where
\begin{equation*}
  \begin{aligned}
    \mathbf{X} \left(\mathbf{f}^{(\ell)}_s, \mathbf{g}^{(\ell)}_s\right)
    = \begin{bmatrix*}
        \ldots, &
        \nabla_{
          \mathbf{R} \left(\mathbf{f}^{(\ell)}_s, \mathbf{g}^{(\ell)}_s\right)_{i,\boldsymbol\cdot}
        }\min_\varepsilon \mathbf{R} \left(\mathbf{f}^{(\ell)}_s, \mathbf{g}^{(\ell)}_s\right)_{i,\boldsymbol\cdot},&
        \ldots
      \end{bmatrix*}.
  \end{aligned}
\end{equation*}

Finally, for all $\ell = 1, \ldots, L-1$ we have
\begin{equation*}
  \begin{aligned}
    \frac{\partial \mathcal{L}}{\partial \mathbf{g}^{(\ell)}_{js}}
     & = \sum_i \frac{\partial \mathcal{L}}{\partial \mathbf{f}^{(\ell+1)}_{is}}\cdot
    \frac{\partial \mathbf{f}^{(\ell+1)}_{is}}{\partial \mathbf{g}^{(\ell)}_{js}}     \\
     & = - \sum_i \widebar{\mathbf{f}}^{(\ell+1)}_{is} \cdot
    \mathbf{X} \left(\mathbf{f}^{(\ell)}_s, \mathbf{g}^{(\ell)}_s\right)_{j,i},
  \end{aligned}
\end{equation*}

and the adjoint for $\mathbf{g}^{(\ell)}_s$, $\ell = 1, \ldots, L-1$, is
\begin{equation*}
  \begin{aligned}
    \widebar{\mathbf{g}}^{(\ell)}_s
    = - \mathbf{X} \left(\mathbf{f}^{(\ell)}_s, \mathbf{g}^{(\ell)}_s\right) \cdot
    \widebar{\mathbf{f}}^{(\ell+1)}_s,
  \end{aligned}
\end{equation*}

and the adjoint for $\mathbf{G}^{(\ell)}$, $\ell = 1, \ldots, L-1$, is
\begin{equation}\label{eqn:adjoint-of-G-ell}
  \begin{aligned}
    \widebar{\mathbf{G}}^{(\ell)} =
    \begin{bmatrix*}
      \bar{\mathbf{g}}^{(\ell)}_1, &
      \ldots,&
      \bar{\mathbf{g}}^{(\ell)}_s,&
      \ldots,&
      \bar{\mathbf{g}}^{(\ell)}_S
    \end{bmatrix*}.
  \end{aligned}
\end{equation}

\begin{algorithm}[H]
  \caption{Log-Stabilized Wasserstein Barycenter Algorithm with Gradients}
  \begin{algorithmic}[1]\label{algo:log-barycenter-with-gradient-Aw}
    \Require $\mathbf{A} \in \Sigma_{M \times S}$, $\widetilde{\mathbf{b}} \in \Sigma_N$,
    $\mathbf{C} \in \mathbb{R}^{M \times N}$, $\mathbf{w} \in \Sigma_S$,
    $\varepsilon > 0$.
    \Initialize $\mathbf{F} = \mathbb{0}_{M \times S}$, $\mathbf{G}_{N \times S} = \mathbb{1}_{N \times S}$,
    $\mathbf{logb} = \mathbb{0}_N$.

    \State \# forward loop
    \While{$\ell = 1, \ldots, L$}
    \State $\mathbf{R}_{min}^{row}$ from \cref{eqn:minrow-R-FG}
    \State $\mathbf{F} = \mathbf{F} + \varepsilon\log \mathbf{A} + \mathbf{R}_{min}^{row}$ \# update F
    \State $\mathbf{R}_{min}^{col}$ from \cref{eqn:mincol-R-FG}
    \State $\mathbf{logb} =
      -\frac1\varepsilon \mathbf{G}\cdot \mathbf{w} - \frac1\varepsilon \mathbf{R}_{min}^{col} \cdot \mathbf{w}$
    \# update logb
    \State $\mathbf{G} = \mathbf{G} + \varepsilon \mathbf{logb} \cdot \mathbb{1}_S^\top + \mathbf{R}_{min}^{col}$
    \# update G
    \EndWhile

    \State \# backward loop
    \State $\widebar{\log\mathbf{b}} = 2 (\mathbf{b} - \widetilde{\mathbf{b}}) \odot \mathbf{b}$
    \State $\widebar{\mathbf{F}}$ from \cref{eqn:adjoint-of-F-L}
    \State $\widebar{\mathbf{A}} = \varepsilon \widebar{\mathbf{F}} \oslash \mathbf{A}$
    \State $\widebar{\mathbf{w}} = (\mathbf{G}^{(L-1)} + \mathbf{R}^{col}_{min})^\top \cdot \widebar{\log\mathbf{b}}$
    \If{$\ell = L-1, \ldots, 1$}
    \State $\widebar{\mathbf{G}}$ from \cref{eqn:adjoint-of-G-ell}
    \State $\widebar{\log\mathbf{b}}$ from \cref{eqn:adjoint-of-logb-ell}
    \State $\widebar{\mathbf{F}}$ from \cref{eqn:adjoint-of-F-ell}
    \State $\widebar{\mathbf{A}} = \widebar{\mathbf{A}} + \varepsilon \widebar{\mathbf{F}} \oslash \mathbf{A}$
    \State $\widebar{\mathbf{w}}
      = \widebar{\mathbf{w}} + (\mathbf{G}^{(\ell-1)} + \mathbf{R}^{col, \ell}_{min})^\top \cdot \widebar{\log\mathbf{b}}$
    \EndIf

    \State $\mathbf{b} = \exp \mathbf{logb}$
    \State $\widebar{\mathbf{w}} = - \widebar{\mathbf{w}} / \varepsilon$

    \Ensure $\mathbf{b}$, $\widebar{\mathbf{A}}$, $\widebar{\mathbf{w}}$.
  \end{algorithmic}
\end{algorithm}

\section{Wasserstein Dictionary Learning}\label{sec:wdl}
\sectionmark{WDL}

\subsection{Numeric Optimization Methods}\label{subsec:numeric-optimization-methods}

Before moving on to discuss the Wasserstein Dictionary Learning algorithm,
let us briefly review some gradient-based optimization algorithms commonly employed in machine learning applications.
For a general treatment on this topic, please refer to \citet{nocedal2006}.

The most famous algorithm for solving an unconstrained optimization problem is probably gradient descent method.
This celebrated method can be attributed to \citet{cauchy1847}, even though he didn't
believe it will actually find the minimum nor did he show its convergence \citep{lemarechal2012}.
\citet{hadamard1908} also seemed to have proposed this idea independently, following the work of \citet{hilbert1900},
in the context of solving differential equations,
as an alternative to Rayleigh-Ritz method \citep{rayleigh1896,ritz1909,courant1943}.

Since gradient descent method only leverages gradient information, it is considered as first-order optimization method.
Though there exist other optimization methods in numeric optimization,
such as second-order methods\footnote{
  For example, Newton's method \citep{nocedal2006}.
} or global optimization methods\footnote{
  E.g. Differential Evolution algorithm \citep{storn1997} and Simulated Annealing \citep{kirkpatrick1983}.
} among others,f
it is the stochastic gradient descent (SGD) method that is \textit{de facto} the only method\footnote{
  One of the most popular neural network library, PyTorch \citep{paszke2017},
  only includes first-order optimizers,
  and they are all variants of SGD with the exception of the celebrated L-BFGS method \citep{liu1989},
  which is a variant of the seminal BFGS method \citep{broyden1970,fletcher1970,goldfarb1970,shanno1970}.
  Among all the SGD variants, Adam optimizer \citep{kingma2015} is the most common one,
  simply because of its empirical performance;
  in fact, all other SGD-based optimizers are either derived from Adam, or its predecessors.
  The list of optimizers implemented in PyTorch
  (\url{https://pytorch.org/docs/stable/optim.html\#algorithms}, version 2.4 as of this writing)
  are accessed on August 10, 2024.
}
being used in modern machine learning applications.
Its history can be traced back to the stochastic approximation method \citep{robbins1951},
its update \citep{kiefer1952},
and was used in the perceptron model \citep{rosenblatt1958} which is the early prototype of the modern
neural networks.

What is now commonly employed is the mini-batch SGD, where a batch of sample are evaluated and their gradients averaged
as the current iteration's estimate for the gradient of the objective function.
Then this average gradient is used in some first-order optimization algorithm updating rule
to update the parameters for the model.
The procedure repeats until convergence.

Here, I list the parameter updating rules for several commonly employed optimizers\footnote{
  Interested readers should refer to \citet{nocedal2006} for a general treatment on the field of numeric optimization.
},
which are implemented in the \textbf{\textit{wig}} package.

\begin{definition}[{\cite[Mini-Batch Gradient]{dekel2012}}]
  For some batch size $B$ and objective function to be minimized
  $f \left(\mathbf{x} ; \boldsymbol\theta_t\right)$ with
  data $\mathbf{x}$ and
  parameters $\boldsymbol\theta_t$ at timestamp $t$,
  define the mini-batch gradient as

  \begin{equation}
    \begin{aligned}
      \mathbf{g}_t =
      \frac1B \sum_{b=1}^B \nabla_\theta f(\mathbf{x}_b; \boldsymbol\theta_{t}),
    \end{aligned}
  \end{equation}

  where $\left\{\mathbf{x}_b\right\}_{b = 1}^{B}$ denote the sequence of data.
\end{definition}

\begin{update}[Stochastic Gradient Descent]
  Given a function to be minimized $f \left(\mathbf{x} ; \boldsymbol\theta_t\right)$
  and a mini-batched gradient $\mathbf{g}_{t-1}$
  with current parameters $\boldsymbol\theta_{t-1}$
  and learning rate $\eta$,
  \begin{equation*}
    \begin{aligned}
      \boldsymbol\theta_t = \boldsymbol\theta_{t-1} - \eta \cdot \mathbf{g}_{t-1}.
    \end{aligned}
  \end{equation*}
\end{update}

\begin{update}[{\cite[Adam]{kingma2015}}]
  Given a function to be minimized $f \left(\mathbf{x} ; \boldsymbol\theta_t\right)$
  and a mini-batched gradient $\mathbf{g}_{t-1}$
  with current parameters $\boldsymbol\theta_{t-1}$
  and learning rate $\eta$,
  \begin{equation*}
    \begin{aligned}
      \mathbf{m}_{t}           & = \beta_1 \cdot \mathbf{m}_{t-1} + (1 - \beta_1)\cdot \mathbf{g}_{t-1},                                       \\
      \mathbf{v}_{t}           & = \beta_2 \cdot \mathbf{v}_{t-1} + (1 - \beta_2)\cdot \mathbf{g}_{t-1}^2,                                     \\
      \widehat{\mathbf{m}}_{t} & = \frac{\mathbf{m}_{t}}{1 - \beta_1^t},                                                                       \\
      \widehat{\mathbf{v}}_{t} & = \frac{\mathbf{v}_{t}}{1 - \beta_2^t},                                                                       \\
      \boldsymbol\theta_{t}    & = \boldsymbol\theta_{t-1} - \eta \cdot \frac{\widehat{\mathbf{m}}_t}{\sqrt{\widehat{\mathbf{v}}_t}+\epsilon},
    \end{aligned}
  \end{equation*}
  where $\left(\cdot\right)^2$ denote element-wise power of 2.
\end{update}

\begin{update}[{\cite[AdamW]{loshchilov2019}}]
  Given a function to be minimized $f \left(\mathbf{x} ; \boldsymbol\theta_t\right)$
  and a mini-batched gradient $\mathbf{g}_{t-1}$
  with current parameters $\boldsymbol\theta_{t-1}$
  decay parameter $\gamma$,
  and learning rate $\eta$,
  \begin{equation*}
    \begin{aligned}
      \mathbf{m}_{t}           & = \beta_1 \cdot \mathbf{m}_{t-1} + (1 - \beta_1)\cdot \widehat{\mathbf{g}}_t,   \\
      \mathbf{v}_{t}           & = \beta_2 \cdot \mathbf{v}_{t-1} + (1 - \beta_2)\cdot \widehat{\mathbf{g}}_t^2, \\
      \widehat{\mathbf{m}}_{t} & = \frac{\mathbf{m}_t}{1 - \beta_1^t},                                           \\
      \widehat{\mathbf{v}}_{t} & = \frac{\mathbf{v}_t}{1 - \beta_2^t}                                            \\
      \boldsymbol\theta_{t}    & =
      \left(1 - \eta \gamma\right) \boldsymbol\theta_{t-1} -
      \eta \frac{\widehat{\mathbf{m}}_t}{\sqrt{\widehat{\mathbf{v}}_t}+\epsilon}.
    \end{aligned}
  \end{equation*}
\end{update}

\subsection{Wasserstein Dictionary Learning}\label{subsec:wasserstein-dictionary-learning}

In the Wasserstein barycenter problem, we have the data as sequence of probability vectors
$\left\{\mathbf{a}_s\right\}_{s = 1}^S$,
and would like to obtain a ``centroid'' of the data $\mathbf{b}$.
Here, $\mathbf{a}_s$'s are known as data, but $\mathbf{b}$ is unknown and thus need to be estimated or computed by optimization.
The so-called Wasserstein Dictionary Learning \citep{schmitz2018} is the almost converse problem:
suppose we have the observed $\mathbf{b}$'s, how to find the latent and thus unknown $\mathbf{a}_s$ such that the barycenter from which would reconstruct the data $\mathbf{b}$.
In this case, it is the data $\mathbf{b}$ that is known, but not $\mathbf{a}_s$'s.

The Wasserstein Dictionary Learning can be considered broadly as one of the (nonlinear) topic models \citep{blei2009},
the most famous of which would probably be Latent Dirichlet Allocation model \citep[LDA]{blei2003}.
As is the case with all other topic models,
we have a sequence of data, and we need to discover the common ``factors'' or ``topics'' among them,
either in natural language processing \citep{xu2018} or image processing \citep{schmitz2018} settings.
In this note, however, I will only be discussing WDL for NLP applications.

In order to apply the Wasserstein Barycenter algorithm as in \cref{sec:wasserstein-barycenter},
we will need the distance matrix ($\mathbf{C}$ in \cref{sec:wasserstein-barycenter}),
but not anymore of $\mathbf{A}$ and $\mathbf{w}$ as they are considered ``trainable'' parameters now
with the gradients from \cref{sec:barycenter-gradient}.

To state the problem,
first suppose that we have $M$ documents with a total of $N$ tokens,
then we can let $\mathbf{Y} \in \mathbb{R}_+^{N\times M}$ be the matrix containing the data,
i.e. $M$ probability vectors each of size $N$,
and thus for each column $\mathbf{y}_m$ we have $\mathbf{y}_m \in \Sigma_N$
for $m = 1, \ldots, M$.
To obtain the only external variable, distance matrix $\mathbf{C} \in \mathbb{R}^{N \times N}$,
we will need the Euclidean distance between the word embedding vectors for any pair of tokens.
In other words, we need to train the Word2Vec \citep{mikolov2013} network
and obtain the embedding vectors $\mathbf{x}_i$ for each of the token in the dataset, $i= 1, \ldots, N$.
Then the entries of the distance matrix $\mathbf{C}_{i,j} = d\,(\mathbf{x}_i, \mathbf{x}_j)$
where $i = 1, \ldots, N$, $j = 1, \ldots, N$,
and $d\, \left(\cdot, \cdot\right)$ is usually considered Euclidean \citep{xu2018,xie2020a}.

For the latent topic\footnote{
  Sometimes also called atom or factor.
} matrix $\mathbf{A} \in \mathbb{R}_+^{N \times S}$ where each column $\mathbf{a}_s \in \Sigma_N$,
and the weight matrix $\mathbf{W} \in \mathbb{R}_+^{S \times M}$ where each column $\mathbf{w}_m \in \Sigma_S$,
we have

\begin{equation}\label{eqn:wasserstein-dictionary-problem-constrained}
  \begin{aligned}
    \min_{\mathbf{A}, \mathbf{W}}
     &
    \sum_{m = 1}^M \mathcal{L} \left(\widehat{\mathbf{b}}_m, \mathbf{b}_m\right),                          \\
    \text{s.t.}
     & \widehat{\mathbf{b}}_m =
    \argmin_{\mathbf{b} \in \Sigma_N}
    \sum_{s = 1}^S \mathbf{W}_{sm} \, \ell_{\mathbf{C}}^\varepsilon \left(\mathbf{a}_s, \mathbf{b}\right), \\
     & \mathbf{A} = \left[\mathbf{a}_1, \ldots, \mathbf{a}_S\right] \in \mathbb{R}_+^{N \times S},
    \text{where } \mathbf{a}_s \in \Sigma_N,                                                               \\
     & \mathbf{W} = \left[\mathbf{w}_1, \ldots, \mathbf{w}_M\right] \in \mathbb{R}_+^{S \times M},
    \text{where } \mathbf{w}_m \in \Sigma_S.
  \end{aligned}
\end{equation}

Essentially, we are considering the $\widehat{\mathbf{b}}_m$,
the computed barycenter from the current $\mathbf{A}$,
as the reconstruction for the true $\mathbf{b}_m$ that we observe in the data.
Then we will need to optimize $\mathbf{A}$ (and also implicitly $\mathbf{W}$)
such that the reconstructions approximate the original data well.
This problem, however, is not guaranteed to be convex \citep{schmitz2018},
whether jointly in $\mathbf{A}$ and $\mathbf{W}$ or for each separately,
unlike the original Sinkhorn problem discussed in \cref{sec:review}.
Thus, we aim to solve this problem by a gradient descent approach to find a local minimum,
as is the usual case in machine learning literature.
In order to do so, we will take the gradient of the loss with respect to the parameters
(here $\mathbf{A}$ and $\mathbf{W}$), and then optimize them by some learning rate,
and then repeat until some stopping criterion.
The gradient computation can be easily carried out by any modern Automatic Differentiation system\footnote{
  They are usually embedded within neural network libraries, e.g. PyTorch \citep{paszke2017}, TensorFlow \citep{abadi2016}, etc.,
  but can also be found as stand-alone libraries, e.g. autodiff in C++ \citep{leal2018} or Zygote in Julia \citep{innes2019}.
  Unlike Symbolic or Numeric Differentiation, Automatic Differentiation (AD) essentially relies on chain's rule
  and could provide accurate gradient computation.
  There are usually two modes for AD systems: forward mode or reverse mode, and each has their own advantages.
  To see this,
  let a function $f: \mathbb{R}^n \to \mathbb{R}^m$, when $n \ll m$, then forward mode is more efficient;
  if, however, $n \gg m$, then reverse mode is more efficient.
  Neural Network libraries are usually implemented in reverse mode,
  since there is usually a scalar loss but with many parameters of the neural network model.
  Interested readers could refer to \citet[Chapter 8]{nocedal2006} for details of AD.
},
but in this note, I manually derive the gradients of all algorithms
(in \cref{sec:sinkhorn-gradient,sec:barycenter-gradient}) to achieve memory-efficient and faster code.

Note, however, we cannot directly optimize the parameters from \cref{eqn:wasserstein-dictionary-problem-constrained},
since the problem is a constrained optimization problem with two constraints on the $\mathbf{a}_s$ and $\mathbf{w}_m$
being in their respective simplices.
Therefore, we need to convert this constrained optimization problem into an unconstrained one \citep{schmitz2018,xu2018}:
that is, to let $\mathbf{a}_s = \text{softmax} \left(\boldsymbol\upalpha_s\right)$
and $\mathbf{w}_m = \text{softmax}\left(\mathbf{w}_s\right)$,
as the $\text{softmax}$ function will ensure its output to be a probability vector,
i.e. $\mathbf{a}_s \in \Sigma_N$ and $\mathbf{w}_m \in \Sigma_S$.
Thus, we have transformed the problem as an optimization problem with parameters
$\boldsymbol\upalpha$ and $\boldsymbol\uplambda$.
Now we need to define the $\text{softmax}$ function.

From the \textit{softmax} functions defined in \cref{lemma:jacobian-softmax-vec,lemma:jacobian-softmax-mat},
we can enforce the constraints $\mathbf{a}_s \in \Sigma_N$ and $\mathbf{w}_m \in \Sigma_S$
by the following change of variables

\begin{equation}
  \begin{aligned}
    \mathbf{A} = \text{softmax} \left(\boldsymbol\upalpha\right),
    \quad\text{ and }\quad
    \mathbf{W} = \text{softmax} \left(\boldsymbol\uplambda\right),
  \end{aligned}
\end{equation}

where $\boldsymbol\upalpha = \left[\boldsymbol\upalpha_1, \ldots, \boldsymbol\upalpha_S\right] \in \mathbb{R}^{N\times S}$,
and $\boldsymbol\uplambda = \left[\boldsymbol\uplambda_1, \ldots, \boldsymbol\uplambda_M\right] \in \mathbb{R}^{S\times M}$.
Let us slightly abuse the notation to denote the vector of losses across the topics

\begin{equation*}
  \begin{aligned}
    \ell_{\mathbf{C}}^\varepsilon \left(\mathbf{A}, \mathbf{B}\right)
    = \left[
      \ell_{\mathbf{C}}^\varepsilon \left(\mathbf{a}_1, \mathbf{b}\right),
      \,\ldots\,,
      \ell_{\mathbf{C}}^\varepsilon \left(\mathbf{a}_S, \mathbf{b}\right)
      \right]^\top \in \mathbb{R}^N.
  \end{aligned}
\end{equation*}

Note that $\ell_{\mathbf{C}}^\varepsilon \left(\mathbf{A}, \mathbf{B}\right)$
is essentially the same as $\ell_{\mathbf{C}}^\varepsilon \left(\mathbf{A}, \mathbf{b}\right)$,
where $\mathbf{B} = \mathbf{b}\cdot \mathbb{1}_S^\top \in \mathbb{R}^{N\times S}$
as the ``many-to-one'' transport problem in \cref{subsec:parallel-sinkhorn}.
Therefore, we can rewrite the problem in \cref{eqn:wasserstein-dictionary-problem-constrained} as
an unconstrained one with respect to parameters $\boldsymbol\upalpha$ and $\boldsymbol\uplambda$

\begin{equation}\label{eqn:wasserstein-dictionary-problem-unconstrained}
  \begin{aligned}
    \min_{\boldsymbol\upalpha, \boldsymbol\uplambda}
    \sum_{m = 1}^M \mathcal{L} \left(\widehat{\mathbf{b}}_m, \mathbf{b}_m\right),
  \end{aligned}
\end{equation}

where $\widehat{\mathbf{b}}_m$, $\mathbf{A}$, and $\mathbf{W}$ satisfy the following

\begin{equation*}
  \begin{aligned}
    \widehat{\mathbf{b}}_m
     & =
    \argmin_{\mathbf{b} \in \Sigma_N}\,
    \mathbf{w}_{m}^\top \cdot \ell_{\mathbf{C}}^\varepsilon \left(\mathbf{A}, \mathbf{b}\right), \\
     & \mathbf{A} = \text{softmax} \left(\boldsymbol\upalpha\right),                             \\
     & \mathbf{W} = \text{softmax} \left(\boldsymbol\uplambda\right).
  \end{aligned}
\end{equation*}

It is obvious that \cref{eqn:wasserstein-dictionary-problem-unconstrained}
is an unconstrained optimization problem,
compared to the original constrained formulation in \cref{eqn:wasserstein-dictionary-problem-constrained},
since the parameters need optimizing are the auxiliary variables
$\boldsymbol\upalpha$ and $\boldsymbol\uplambda$,
and all we need to do is to obtain the gradients through the loss,
and ``backproporgate'' to update $\boldsymbol\upalpha$ and $\boldsymbol\uplambda$.
Thus, since we know that
$\mathbf{A} = \text{softmax} \left(\boldsymbol\upalpha\right)$
and $\mathbf{W} = \text{softmax} \left(\boldsymbol\uplambda\right)$,
we can retrieve the topics' distribution $\mathbf{A}$ and the weights $\mathbf{W}$.

\subsection{WDL Algorithm}\label{subsec:wdl-algorithm}

To leverage the idea of Wasserstein Barycenter and use it in a ``learning'' environment,
i.e.~Wasserstein Dictionary Learning (WDL),
we will take some random initialized parameters $\boldsymbol\upalpha$ and $\boldsymbol\uplambda$,
calculate their $softmax$ column-wise to get $\mathbf{A}$ and $\mathbf{W}$,
calculate their optimal barycenter,
calculate the loss between the barycenter and the original data,
take the gradients of the loss with respect to the parameters,
and then update the parameters based on the gradients.
The computation of the optimal barycenter is based on the Wasserstein barycenter problem,
detailed in \cref{sec:wasserstein-barycenter} with its gradient calculation in \cref{sec:barycenter-gradient},
and the parameter-updating routine can be chosen to be any optimization technique,
for example \citet{schmitz2018} chooses L-BFGS,
and \citet{xie2020} uses Adam.

To lay out the algorithm of WDL, we are taking what we have built so far in this entire note,
but there are a few more things remain to be derived.
Notice that the gradient
\begin{equation*}
  \begin{aligned}
    \nabla_{\boldsymbol\upalpha}\, \mathcal{L}(\widehat{\mathbf{b}}, \mathbf{b})
    = \begin{bmatrix*}
        \nabla_{\boldsymbol\upalpha_1}\, \mathcal{L}(\widehat{\mathbf{b}}, \mathbf{b}) &
        \cdots &
        \nabla_{\boldsymbol\upalpha_S}\, \mathcal{L}(\widehat{\mathbf{b}}, \mathbf{b})
      \end{bmatrix*},
  \end{aligned}
\end{equation*}
where each $\nabla_{\boldsymbol\upalpha_s}\, \mathcal{L}(\widehat{\mathbf{b}}, \mathbf{b})
  = \left(D_{\boldsymbol\upalpha_s}\, \mathcal{L}\right)^\top$,
$\boldsymbol\upalpha = \left[\boldsymbol\upalpha_1, \ldots, \boldsymbol\upalpha_s, \ldots\boldsymbol\upalpha_S\right]$,
and $\boldsymbol\upalpha_s \in \mathbb{R}^N$.
Then we have
\begin{equation*}
  \begin{aligned}
    \nabla_{\boldsymbol\upalpha_s}\, \mathcal{L}
    = \left(D_{\boldsymbol\upalpha_s}\, \mathcal{L}\right)^\top
    = \left(
    D_{\mathbf{A}} \mathcal{L} \cdot
    D_{\boldsymbol\upalpha_s} \mathbf{A}
    \right)^\top
    = \left(D_{\boldsymbol\upalpha_s} \mathbf{A}\right)^\top \cdot
    \left(D_{\mathbf{A}} \mathcal{L}\right)^\top
    = \left(D_{\boldsymbol\upalpha_s} \mathbf{A}\right)^\top \cdot \nabla_{\mathbf{A}} \mathcal{L}.
  \end{aligned}
\end{equation*}

Therefore, to derive the above gradient, we will need the Jacobian of the \textit{softmax} function,
which is given by \cref{lemma:jacobian-softmax-vec}.
The other term $\nabla_{\mathbf{A}} \mathcal{L}$ is derived in \cref{sec:barycenter-gradient},
with the quadratic loss between the computed barycenter and some ``data vector'' $\mathbf{b}$.

\begin{algorithm}[H]
  \caption{Wasserstein Dictionary Learning Algorithm}
  \begin{algorithmic}[1]\label{algo:wdl-algorithm}
    \Require $\mathbf{Y} \in \mathbb{R}^{N \times M}$, $\mathbf{C} \in \mathbb{R}^{N \times N}$,
    $\varepsilon > 0$, $S \in \mathbb{N}_+$
    \Initialize $\boldsymbol\upalpha \in \mathbb{R}^{N \times S}$,
    $\boldsymbol\uplambda \in \mathbb{R}^{S \times M}$
    \State $\mathbf{K} = \exp(-\frac{\mathbf{C}}{\varepsilon})$
    \# only used with \cref{algo:parallel-barycenter-with-gradient-Aw}
    \State $\mathbf{A} = \text{softmax} (\boldsymbol\upalpha)$, $\mathbf{W} = \text{softmax} (\boldsymbol\uplambda)$

    \While{Not Convergence}
    \State $\mathbf{g_{\boldsymbol\upalpha}} = \mathbb{0}_{N \times S}$,
    $\mathbf{g_{\boldsymbol\uplambda}} = \mathbb{0}_{S}$

    \For{Each $\ell = 1, \ldots, B$ in Batch B}
    \State get $\mathbf{y}^{\ell} \in \mathbb{R}^N$ from $\mathbf{Y}$,
    and get $\mathbf{l} \in \mathbb{R}^S$ from $\boldsymbol\uplambda$

    \State $\mathbf{b}^\ell$, $\widebar{\mathbf{A}}^\ell$, $\widebar{\mathbf{w}}^\ell$
    from \cref{algo:parallel-barycenter-with-gradient-Aw} or \cref{algo:log-barycenter-with-gradient-Aw}

    \State $\mathbf{J_A}$, $\mathbf{J_w}$ from \cref{eqn:jacobian-softmax-vec}

    \State $\mathbf{g_{\boldsymbol\upalpha}}
      = \mathbf{g_{\boldsymbol\upalpha}} + \mathbf{J_A}^\top \cdot \widebar{\mathbf{A}}^\ell$,
    $\mathbf{g_{\boldsymbol\uplambda}}
      = \mathbf{g_{\boldsymbol\uplambda}} + \mathbf{J_w}^\top \cdot \widebar{\mathbf{w}}^\ell$
    \EndFor
    \State $\mathbf{g_{\boldsymbol\upalpha}} = \mathbf{g_{\boldsymbol\upalpha}} / B$,
    $\mathbf{g_{\boldsymbol\uplambda}} = \mathbf{g_{\boldsymbol\uplambda}} / B$ \# average gradient
    \State Update $\boldsymbol\upalpha$ and $\boldsymbol\uplambda$ by \cref{subsec:numeric-optimization-methods}
    with $\mathbf{g}_{\boldsymbol\upalpha}$ and $\mathbf{g}_{\boldsymbol\uplambda} \cdot \mathbb{1}_M^\top$
    \State $\mathbf{A} = \text{softmax} (\boldsymbol\upalpha)$, $\mathbf{W} = \text{softmax} (\boldsymbol\uplambda)$
    \EndWhile
    \Ensure $\mathbf{A}$, $\mathbf{W}$
  \end{algorithmic}
\end{algorithm}

\newpage
\addcontentsline{toc}{section}{References}
\printbibliography

\newpage

\begin{appendices}\label{appendix}

\section{Mathematical Notations}\label{appendix:math-notation}

Mathematical notations used in this note:

\begin{itemize}
  \item $\mathbf{a}$: lower case letters in bold denote vectors
  \item $\mathbf{A}$: upper case letters in bold denote matrices
  \item $\mathbb{1}_n$, $\mathbb{0}_n$: a vector of ones (or zeros) of length $n$
  \item $\Sigma_n$: probability simplex $\Sigma_n \equiv \left\{
          \mathbf{a} \in \mathbb{R}_+^n: \sum_{i=1}^{n} \mathbf{a}_i = 1
          \right\}.$
  \item $\langle \cdot, \cdot\rangle $: inner (dot) product of two matrices of same size
  \item $\odot$: element-wise (Hadamard) multiplication
  \item $\oslash$ or $\frac1{\mathbf{X}}$: element-wise (Hadamard) division
  \item $\log(\cdot)$, $\exp(\cdot)$: element-wise logarithm and exponential functions
  \item $\diag(\mathbf{x})$: create a diagonal matrix $n \times n$ from vector $\mathbf{x} \in \mathbb{R}^n$
  \item $\vec \mathbf{A}\mathbf{B}
          = \left(I_m \otimes \mathbf{A}\right)\vec \mathbf{B}
          = \left(\mathbf{B}^\top \otimes I_k \right) \vec \mathbf{A}$: for $\mathbf{A}: k\times l$ and $\mathbf{B}:l\times m$
  \item $\mathcal{K}^{(m,n)}$: Commutation matrix such that $\mathcal{K}^{(m,n)}\vec \mathbf{A} = \vec \left(\mathbf{A}^\top\right)$ for a matrix $\mathbf{A} \in \mathbb{R}^{m \times n}$
  \item $\left[P\right]$: Iverson bracket, i.e.~$\left[P\right] = 1$ if $P$ is true, and $\left[P\right] = 0$ otherwise.
  \item Matrix calculus notations used throughout this note are based on \citet{magnus2019}
\end{itemize}

\newpage
\section{Lemmas}\label{appendix:lemmas}

\begin{lemma}[]\label{lemma:ones-otimes-b-transpose}
  For any column vectors $\mathbf{a} \in \mathbb{R}^m$ and $\mathbf{b} \in \mathbb{R}^n$,
  $\mathbf{a} \otimes \mathbf{b}^\top = \mathbf{a} \cdot \mathbf{b}^\top$.
\end{lemma}

\begin{proof}[Proof of \cref{lemma:ones-otimes-b-transpose}]
  To see this, we use the property of Kronecker product,

  \begin{equation*}
    \begin{aligned}
      \mathbf{a} \otimes \mathbf{b}^\top
       & = \left(\mathbf{a} \cdot I_1\right) \otimes \left(I_1 \cdot \mathbf{b}^\top\right)
       & = \left(\mathbf{a} \otimes I_1\right) \cdot \left(I_1 \otimes \mathbf{b}^\top\right)
       & = \mathbf{a} \cdot \mathbf{b}^\top.
    \end{aligned}
  \end{equation*}
\end{proof}

\begin{lemma}[]\label{lemma:ones-transpose-otimes-b}
  For any column vectors $\mathbf{a} \in \mathbb{R}^m$ and $\mathbf{b} \in \mathbb{R}^n$,
  $\mathbf{a}^\top \otimes \mathbf{b} = \mathbf{b} \cdot \mathbf{a}^\top$.
\end{lemma}

\begin{proof}[Proof of \cref{lemma:ones-transpose-otimes-b}]
  To see this, we use the same property as in the previous lemma,

  \begin{equation*}
    \begin{aligned}
      \mathbf{a}^\top \otimes \mathbf{b}
       & = \left(I_1  \cdot \mathbf{a}^\top\right) \otimes \left(\mathbf{b} \cdot I_1\right)
       & = \left(I_1 \otimes \mathbf{b}\right) \cdot \left(\mathbf{a}^\top \otimes I_1\right)
       & = \mathbf{b} \cdot \mathbf{a}^\top.
    \end{aligned}
  \end{equation*}
\end{proof}

\begin{lemma}[]\label{lemma:diagx-A-diagy}
  Given a matrix $\mathbf{A} \in \mathbb{R}^{m\times n}$ and two vectors
  $\mathbf{x} \in \mathbb{R}^m$, $\mathbf{y} \in \mathbb{R}^n$,
  $\diag \mathbf{x} \cdot \mathbf{A} \cdot \diag \mathbf{y} = (\mathbf{x} \mathbf{y}^\top) \odot \mathbf{A}$.
\end{lemma}

\begin{proof}[Proof of \cref{lemma:diagx-A-diagy}]
  In the derivation of the Jacobians, we frequently encounter matrix computations of the form
  $\diag \mathbf{x} \cdot \mathbf{A} \cdot \diag \mathbf{y}$.
  To avoid the creation of highly sparse diagonal matrices by the $\diag(\cdot)$ function from the vectors,
  we can compute the outer product by the two vectors first and then compute its Hadamard product
  with the internal matrix.
  Not only will this give us a formulation that is more compact,
  but also more efficient in the computation.

  To see this, note that $\mathbf{x} \mathbf{y}^T$ gives us an $m \times n$ matrix with each element being
  $\left(\mathbf{x} \mathbf{y}^\top\right)_{ij} = \mathbf{x}_i \mathbf{y}_j$,
  then the elements of the final Hadamard product would be
  $\left(\left(\mathbf{x} \mathbf{y}^\top\right)\odot \mathbf{A}\right)_{ij} = \mathbf{A}_{ij} \mathbf{x}_i \mathbf{y}_j$.

  Now, on the other hand, if we look at the form $\diag \mathbf{x} \cdot \mathbf{A} \cdot \diag \mathbf{y}$,
  the element is also
  $\left(\diag \mathbf{x} \cdot \mathbf{A} \cdot \diag \mathbf{y}\right)_{ij} = \mathbf{x}_i \mathbf{A}_{ij} \mathbf{y}_j$,
  since for each row $i$ in $\mathbf{A}$, we need to multiply by $\mathbf{x}_i$ and
  for each column $j$ in $\mathbf{A}$, we need to multiply by $\mathbf{y}_j$.
  This concludes the proof as the two forms are identical.
\end{proof}

\begin{lemma}[Gradient of \textit{soft-minimum} function]\label{lemma:grad-soft-minimum}
  For the \textit{soft-minimum} function defined in \cref{eqn:soft-minimum-stable},
  a vector $\mathbf{z} \in \mathbb{R}^n$ with its minimum as $\ubar{\mathbf{z}}$, we have its gradient as
  \begin{equation*}
    \begin{aligned}
      \nabla_{\mathbf{z}} \text{min}_{\varepsilon} \mathbf{z}
      = \frac1{\mathbb{1}_n^\top \cdot \exp \left(-\frac{\mathbf{z} - \ubar{\mathbf{z}}}{\varepsilon}\right)}
      \cdot
      \exp \left(-\frac{\mathbf{z} - \ubar{\mathbf{z}}}{\varepsilon}\right).
    \end{aligned}
  \end{equation*}
\end{lemma}

\begin{proof}[Proof of \cref{lemma:grad-soft-minimum}]
  Given the \textit{soft-minimum} function, we can take the differential
  \begin{equation*}
    \begin{aligned}
      \frac{\partial \min_{\varepsilon} \mathbf{z}}{\partial \mathbf{z}_i}
      = \frac{\partial \ubar{\mathbf{z}}}{\partial \mathbf{z}_i}
      - \varepsilon \, \frac1{
        \sum_j \exp \left(-\frac{\mathbf{z}_j - \ubar{\mathbf{z}}}{\varepsilon}\right)
      }\cdot
      \frac{\partial}{\partial \mathbf{z}_i}
      \left[
        \sum_j \exp \left(-\frac{\mathbf{z}_j - \ubar{\mathbf{z}}}{\varepsilon}\right)
        \right].
    \end{aligned}
  \end{equation*}

  Then we can see that

  \begin{dmath*}
    \frac{\partial}{\partial \mathbf{z}_i}
    \left[
      \sum_j \exp \left(-\frac{\mathbf{z}_j - \ubar{\mathbf{z}}}{\varepsilon}\right)
      \right]
    = \sum_j \frac{\partial }{\partial \mathbf{z}_i} \left[
      \exp \left(-\frac{\mathbf{z}_j - \ubar{\mathbf{z}}}{\varepsilon}\right)
      \right]
    = \frac{\partial }{\partial \mathbf{z}_i}
    \exp \left(-\frac{\mathbf{z}_i - \ubar{\mathbf{z}}}{\varepsilon}\right)
    +
    \sum_{j \not= i}
    \frac{\partial }{\partial \mathbf{z}_i} \left[
      \exp \left(-\frac{\mathbf{z}_j - \ubar{\mathbf{z}}}{\varepsilon}\right)
      \right]
    = -\frac1\varepsilon
    \exp \left(-\frac{\mathbf{z}_i - \ubar{\mathbf{z}}}{\varepsilon}\right)
    \left(1 - \frac{\partial \ubar{\mathbf{z}}}{\partial \mathbf{z}_i}\right)
    + \frac1\varepsilon \sum_{j \not= i}
    \exp \left(-\frac{\mathbf{z}_j - \ubar{\mathbf{z}}}{\varepsilon}\right)
    \cdot \frac{\partial \ubar{\mathbf{z}}}{\partial \mathbf{z}_i}
    = -\frac1\varepsilon \exp \left(-\frac{\mathbf{z}_i - \ubar{\mathbf{z}}}{\varepsilon}\right)
    + \frac1\varepsilon \sum_j \exp \left(-\frac{\mathbf{z}_j - \ubar{\mathbf{z}}}{\varepsilon}\right)
    \cdot \frac{\partial \ubar{\mathbf{z}}}{\partial \mathbf{z}_i}.
  \end{dmath*}

  Putting everything together, we have
  \begin{equation*}
    \begin{aligned}
      \frac{\partial \min_{\varepsilon} \mathbf{z}}{\partial \mathbf{z}_i}
      =
      \frac1{\sum_j \exp \left(-\frac{\mathbf{z}_j - \ubar{\mathbf{z}}}{\varepsilon}\right)}
      \cdot \exp \left(-\frac{\mathbf{z}_j - \ubar{\mathbf{z}}}{\varepsilon}\right).
    \end{aligned}
  \end{equation*}

  One thing that worth mentioning is about $\frac{\partial \ubar{\mathbf{z}}}{\partial \mathbf{z}_i}$.
  Even though this term was cancelled out in the final expression,
  it was handled differently in different AD implementations.
  Julia Zygote implements it with ``left-most'' or ``first-occurrence'' rule,
  while PyTorch implements it as ``subgradient averaging'' approach.
\end{proof}

\begin{lemma}[Jacobian of \textit{softmax} function for a vector]\label{lemma:jacobian-softmax-vec}
  For the \textit{softmax} function and a vector $\mathbf{x} \in \mathbb{R}^n$ with its maximum as $\bar{\mathbf{x}}$,
  we have
  \begin{equation*}
    \begin{aligned}
      \text{softmax} \left(\mathbf{x}\right)
      = \frac{\exp (\mathbf{x})}{\mathbb{1}_n^\top \cdot \exp (\mathbf{x})}
      = \frac{\exp (\mathbf{x}- \bar{\mathbf{x}})}{\mathbb{1}_n^\top \cdot \exp (\mathbf{x} - \bar{\mathbf{x}})},
    \end{aligned}
  \end{equation*}
  where the second expression would be numerically stable.
  Moreover, we can derive its Jacobian w.r.t.~$\mathbf{x}$ as
  \begin{equation}\label{eqn:jacobian-softmax-vec}
    \begin{aligned}
      D_{\mathbf{x}}\, \text{softmax} (\mathbf{x})
      = \diag \text{softmax} (\mathbf{x}) - \text{softmax} (\mathbf{x}) \cdot \text{softmax} \left(\mathbf{x}\right)^\top.
    \end{aligned}
  \end{equation}
\end{lemma}

\begin{proof}[Proof of \cref{lemma:jacobian-softmax-vec}]
  To make the computation of the \textit{softmax} function stable, we will need to subtract the maximum value
  on both the numerator and denominator and the result would be the same.
  Therefore, its derivative is also the same.
  To get the Jacobian of $\text{softmax}$, we have
  \begin{equation*}
    \begin{aligned}
      \frac{\partial\, \text{softmax}\left(\mathbf{x}\right)_i}{\partial \mathbf{x}_j}
       & = \frac{
        \frac{\partial}{\partial \mathbf{x}_j} \exp \mathbf{x}_i \cdot
        \sum_j \exp \mathbf{x}_j
        - \exp \mathbf{x}_i \cdot \frac{\partial}{\partial \mathbf{x}_j} \sum_j \exp \mathbf{x}_j
      }{\left(\sum_j \exp \mathbf{x}_j\right)^2}              \\
       & = \frac{\exp \mathbf{x}_i}{\sum_j \exp \mathbf{x}_j}
      \cdot \frac{\partial \mathbf{x}_i}{\partial \mathbf{x}_j} -
      \frac{\exp \mathbf{x}_i}{\sum_j \exp \mathbf{x}_j} \cdot
      \frac{\exp \mathbf{x}_j}{\sum_j \exp \mathbf{x}_j}      \\
       & =
      \text{softmax} \left(\mathbf{x}\right)_i \cdot \frac{\partial \mathbf{x}_i}{\partial \mathbf{x}_j}
      - \text{softmax} \left(\mathbf{x}\right)_i \cdot \text{softmax} \left(\mathbf{x}\right)_j.
    \end{aligned}
  \end{equation*}

  Therefore, the Jacobian matrix is
  \begin{equation*}
    \begin{aligned}
      D_{\mathbf{x}}\, \text{softmax} \left(\mathbf{x}\right)
      = \diag\text{softmax} \left(\mathbf{x}\right)
      - \text{softmax}\left(\mathbf{x}\right) \cdot \text{softmax}\left(\mathbf{x}\right)^\top.
    \end{aligned}
  \end{equation*}
\end{proof}

\begin{lemma}[Jacobian of \textit{softmax} function for a matrix]\label{lemma:jacobian-softmax-mat}
  Suppose we have a matrix $\mathbf{X} \in \mathbb{R}^{m\times n}$,
  and let
  \begin{equation*}
    \begin{aligned}
      \text{softmax}\left(\mathbf{X}\right) = \left[
        \text{softmax}\left(\mathbf{x}_1\right),
        \ldots,
        \text{softmax}\left(\mathbf{x}_n\right)
        \right],
    \end{aligned}
  \end{equation*}
  where $\text{softmax}(\mathbf{x}_i)$ for $i = 1, \ldots, n$ is from \cref{lemma:jacobian-softmax-vec}.
  Moreover, we have
  \begin{equation}\label{eqn:jacobian-softmax-mat}
    \begin{aligned}
      D_{\mathbf{X}} \, \text{softmax} \left(\mathbf{X}\right)
      = \begin{bmatrix*}
          \mathbf{J}_1& 0& \ldots& 0\\
          0 &\mathbf{J}_2 &\ldots &0\\
          \vdots & \vdots & \ddots &\vdots\\
          0 & 0 & \ldots & \mathbf{J}_n
        \end{bmatrix*} \in \mathbb{R}^{mn \times mn},
    \end{aligned}
  \end{equation}
  where $\mathbf{J}_j \in \mathbb{R}^{m\times m}$ is the Jacobian of \textit{softmax} of $\mathbf{X}_{\boldsymbol\cdot, j}$.
\end{lemma}

\begin{proof}[Proof of \cref{lemma:jacobian-softmax-mat}]
  Let $\mathbf{Y} = \text{softmax}\left(\mathbf{X}\right) \in \mathbb{R}^{m\times n}$,
  with each column $\mathbf{Y}_{\boldsymbol\cdot, j} = \text{softmax} \left(\mathbf{X}_{\boldsymbol\cdot, j}\right)$,
  and
  \begin{equation*}
    \begin{aligned}
      \text{softmax}\left(\mathbf{X}_{\boldsymbol\cdot,j}\right)_i
      = \frac{\exp\mathbf{X}_{ij}}{\sum_{k=1}^m \exp \mathbf{X}_{kj}},
    \end{aligned}
  \end{equation*}
  for $i = 1, \ldots, m$ and $j = 1, \ldots, n$.
  Since the \textit{softmax} function is applied independently to each column, the Jacobian of $\mathbf{Y}$
  w.r.t.~$\mathbf{X}$ is a block-diagonal matrix of the individual \textit{softmax} Jacobians.
  Therefore,
  the Jacobian of \textit{softmax}$\mathbf{X}$ takes the form of \cref{eqn:jacobian-softmax-mat},
  i.e.
  \begin{equation*}
    \begin{aligned}
      D_{\mathbf{X}} \, \mathbf{Y}
      = D_{\mathbf{X}} \, \text{softmax} \left(\mathbf{X}\right),
    \end{aligned}
  \end{equation*}
  where $\mathbf{J}_j$ is given by \cref{lemma:jacobian-softmax-vec},
  i.e.
  \begin{equation*}
    \begin{aligned}
      \mathbf{J}_j
       & = D_{\mathbf{X}_{\boldsymbol\cdot, j}}\, \text{softmax} \left(\mathbf{X}_{\boldsymbol\cdot, j}\right) \\
       & = \diag \text{softmax} \left(\mathbf{X}_{\boldsymbol\cdot, j}\right)
      - \text{softmax} \left(\mathbf{X}_{\boldsymbol\cdot, j}\right) \cdot
      \text{softmax} \left(\mathbf{X}_{\boldsymbol\cdot, j}\right).
    \end{aligned}
  \end{equation*}
\end{proof}

\end{appendices}

\end{document}